\def\no{\if01}
\def\iftwelvept{\no}

\def\ifusepdf{\no}
\def\ifpsfont{\no}

\iftwelvept
\documentclass[leqno,12pt]{amsart}
\else
\documentclass[leqno,11pt]{amsart}
\fi
\usepackage{amssymb}
\usepackage{amscd}
\usepackage{latexsym}
\usepackage{verbatim}
\usepackage[all]{xy}

\ifusepdf
\usepackage{hyperref}
\else\fi
\ifpsfont
\usepackage[T1]{fontenc}
\usepackage{times}
\else\fi

\iftwelvept
\else\fi

\theoremstyle{plain}
\newtheorem{Theorem}{Theorem}[section]
\newtheorem{Warning}{Warning}[Theorem]
\newtheorem{TTheorem}{Theorem}

\newtheorem{Proposition}[Theorem]{Proposition}
\newtheorem{Lemma}[Theorem]{Lemma}
\newtheorem{Corollary}[Theorem]{Corollary}

\theoremstyle{definition}

\newtheorem{Definition}[Theorem]{Definition}
\newtheorem{Remark}[Theorem]{Remark}
\newtheorem{Example}[Theorem]{Example}

\newcommand{\ZZ}{{\mathbb{Z}}}
\newcommand{\QQ}{{\mathbb{Q}}}

\newcommand{\NN}{{\mathbb{N}}}
\newcommand{\NNNN}{\operatorname{N}}

\newcommand{\uCat}{\operatorname{Cat}_\infty}

\newcommand{\Sh}{\operatorname{Sh}}

\newcommand{\KKK}{\mathbf{K}}

\newcommand{\Rep}{\operatorname{Rep}}
\newcommand{\PRep}{\operatorname{PRep}}

\newcommand{\oDTM}{\overline{\mathsf{DTM}}}
\newcommand{\MMM}{\mathcal{M}}

\newcommand{\dga}{\mathsf{dga}}

\newcommand{\DM}{\mathsf{DM}}
\newcommand{\DTM}{\mathsf{DTM}}
\newcommand{\Art}{\mathsf{Art}}
\newcommand{\Hom}{\operatorname{Hom}}

\newcommand{\Spec}{\operatorname{Spec}}

\newcommand{\Gal}{\operatorname{Gal}}

\newcommand{\Cor}{\operatorname{Cor}}
\newcommand{\SP}{\operatorname{Sp}}

\newcommand{\PSP}{\operatorname{PSp}}

\newcommand{\Mod}{\operatorname{Mod}}
\newcommand{\PMod}{\operatorname{PMod}}

\newcommand{\SSS}{\mathcal{S}}
\newcommand{\SEC}{\operatorname{Sec}}

\newcommand{\Map}{\operatorname{Map}}

\newcommand{\Fun}{\operatorname{Fun}}
\newcommand{\Alg}{\operatorname{Alg}}

\newcommand{\PSh}{\operatorname{PSh}}
\newcommand{\Aff}{\operatorname{Aff}}
\newcommand{\Grp}{\operatorname{Grp}}

\newcommand{\FIN}{\operatorname{Fin}_\ast}

\newcommand{\wCat}{\widehat{\textup{Cat}}_{\infty}}

\newcommand{\CAlg}{\operatorname{CAlg}}

\newcommand{\Ind}{\operatorname{Ind}}
\newcommand{\sCat}{\operatorname{Cat}^{\textup{st}}_{\infty}}

\newcommand{\Aut}{\operatorname{Aut}}

\newcommand{\Sym}{\operatorname{Sym}}
\newcommand{\Proof}{{\sl Proof.}\quad}
\newcommand{\QED}{{\unskip\nobreak\hfil\penalty50\quad\null\nobreak\hfil
{$\Box$}\parfillskip0pt\finalhyphendemerits0\par\medskip}}

\begin{document}

\title{Bar construction and Tannakization}

\begin{abstract}
In this note we continue our development of tannakizations of
symmetric monoidal $\infty$-categories, begun in \cite{Tan}.
The issue treated in this note is the calculation of tannakizations of
 examples of symmetric
monoidal stable $\infty$-categories with fiber functors.
We consider the case of symmetric monoidal $\infty$-categories of
perfect complexes on perfect derived stacks.
The first main result especially says that
our tannakization includes the bar construction for an augmented commutative ring
spectrum and its equivariant version as a special case.
We apply it to the study of the
tannakization of the stable infinity-category of
mixed Tate motives over a perfect field.
We prove that its tannakization can be obtained from the
$\mathbb{G}_m$-equivariant bar construction of a commutative differential
graded algebra equipped with $\mathbb{G}_m$-action.
Moreover, under Beilinson-Soul\'e vanishing conjecture,
we prove that the underlying group scheme of the tannakization
is the motivic Galois group for mixed Tate motives, constructed in \cite{BK}, \cite{KM}, \cite{Lev}. The case of Artin motives is also included.
\end{abstract}

\author{Isamu Iwanari}

\thanks{The author is partially supported by Grant-in-Aid for Scientific Research,
Japan Society for the Promotion of Science.}

\address{Mathematical Institute, Tohoku University, Sendai, Miyagi, 980-8578,
Japan}

\email{iwanari@math.tohoku.ac.jp}

\maketitle

\section{Introduction}

In \cite{Tan} we have constructed tannakizations of symmetric monoidal $\infty$-categories.
Let $R$ be a commutative ring spectrum.
Let $\mathcal{C}^\otimes$ be a small symmetric monoidal
$\infty$-category,
equipped with a symmetric monoidal functor
$F:\mathcal{C}^\otimes \to \PMod_R^\otimes$ where
$\PMod_R^\otimes$ denotes the symmetric monoidal
$\infty$-category of compact $R$-spectra.
(Despite we use the machinery of quasi-categories in the text,
by an $\infty$-category we informally mean an $(\infty,1)$-category
in this introduction.)
In loc. cit., given $F:\mathcal{C}^\otimes \to \PMod_R^\otimes$
we construct a derived affine group scheme $G$
over $R$ which
represents the automorphism group of $F$
and has a certain universality (see Theorem~\ref{premain}).
A derived affine group schme
is an analogue of an affine group scheme in derived
algebraic geometry \cite{HAG2}, \cite{DAGn}.
For simplicity,
we shall call it the tannakization of $F:\mathcal{C}^\otimes \to \PMod_R^\otimes$. This construction was applied
to the stable $\infty$-category of mixed motives to construct derived motivic
Galois groups and underived motivic Galois groups.
For the reader who is not familiar with
the higher category theory
it is worth emphasizing that the $\infty$-categorical framework
is crucial for the nice representability of automorphism groups, whereas
coarse machineries, such as triangulated categories, prevent
us from getting it.

The purpose of this note is to calculate
tannakizations of some examples of
$F:\mathcal{C}^\otimes \to \PMod_R^\otimes$; our principal interest here is
the case when $\mathcal{C}^\otimes$ is the symmetric monoidal
$\infty$-category $\PMod_Y^\otimes$
of perfect complexes on a derived stack $Y$
and $F$ is induced by $\Spec R\to Y$.
We will study the tannakization
under the assumption of perfectness on derived
stacks, introduced in \cite{BFN},
which particularly includes two cases:
\begin{enumerate}
\renewcommand{\labelenumi}{(\roman{enumi})}
\item $Y$ is an affine  derived
scheme over $R$, that is, $Y=\Spec A$ over $\Spec R$ with $A$ a commutative ring spectrum,
\item $Y$ is
the quotient stack $[X/G]$
where $X$ is an affine derived scheme $X=\Spec A$ and $G$ is an algebraic
group in characteristic zero.
\end{enumerate}
We note that for our purpose
the assumption of affineness on $Y$ in (i) and $X$ in (ii) is
not essential since $\PMod^\otimes_Y\to \PMod_R^\otimes$
depends only on a Zariski neighborhood of the image of $\Spec R\to Y$.
Also, we remark that $A$ in (i) and (ii) can be nonconnective.
Our result may be expressed as follows (cf. Theorem~\ref{perfectkan},
Corollary~\ref{main}):

\begin{TTheorem}
\label{Thm1}
Let $Y$ be a derived stack over $R$ and $\Spec R\to Y$ a section of
the structure map $Y\to \Spec R$. Let $\PMod_{Y}^\otimes\to \PMod_R^\otimes$
be the associated pullback symmetric monoidal functor.
Suppose that $Y$ is
perfect (the cases (i) and (ii) satisfy this property).
Let $G$ be the derived affine group scheme
arising from \v{C}ech nerve associated to $\Spec R\to Y$.
Then the tannakization of
$\PMod_Y^\otimes\to \PMod_R^\otimes$ is equivalent to $G$.
\end{TTheorem}

{\it Bar construction and equivariant bar construction.}
One of our motivations of this note
arises from comparison between
derived group schemes obtained by tannakization
and bar constructions and its variants.
Bar construction has been an important device
in various contexts of homotopy theory, mixed Tate motives
and non-abelian Hodge theory, etc.
In the case (i), \v{C}ech nerve in $\Aff_R$
associated to $\Spec R\to Y=\Spec A$, which we can regard as a
derived affine
group scheme over $R$,
is known as the bar construction of an augmented commutative
ring spectrum (or commutative differential graded algebra)
whose explicit construction can be given by bar resolutions.
In the case (ii), we can think of the \v{C}ech nerve
as the $G$-equivariant version of the bar construction.
As a matter of fact, our actual aim is to study a relationship
between our tannakization and bar constructions and its equivariant
versions; Theorem~\ref{Thm1} especially means that
our method of tannakizations includes bar constructions and the equivariant
versions as a special case.
This allows one to link bar constructions and the variants
to more general method of tannakizations.
It will be of use in the subsequent works
as a key ingredient (see the end of this introduction).

{\it Mixed Tate motives.}
It would be worth mentioning that
the equivariant versions are also important
to applications to the motivic contexts: for instance,
in order to take weight structures into account,
one often uses $\mathbb{G}_m$-equivariant version of
bar construction.
Our results fit very naturally in with the structure of
mixed Tate motives.
In Section 6 and 7, we will study the applications to mixed Tate motives.
Let $\DM^\otimes:=\DM^\otimes(k)$ be the symmetric monoidal stable $\infty$-category of mixed motives over a base scheme $\Spec k$, where $k$ is a perfect field
 (see Section 6.1 for our convention).
We work with coefficients of a field $\KKK$ of characteristic zero;
all stable $\infty$-categories are $H\KKK$-linear, where $H\KKK$ denotes
the Eilenberg-MacLane spectrum.
Let $\DTM_{\vee}^\otimes\subset \DM^\otimes$ be the small symmetric monoidal
stable $\infty$-category of mixed Tate motives which admit duals (see Section 6.2).
For a mixed Weil cohomology theory (such as \'etale cohomology, de Rham cohomology), there exists a homological realization functor
$\mathsf{R}_T:\DTM^\otimes_\vee\to \PMod_{H\KKK}^\otimes$, that is a
$H\KKK$-linear symmetric monoidal exact functor (the field of coefficients
$\KKK$ depends on the choice of a mixed Weil cohomology theory).
By applying the above theorem, we deduce Theorem~\ref{tatemain} which informally says:

\begin{TTheorem}
Let $\mathsf{MTG}=\Spec B$ be the derived affine group scheme obtained
as the tannakization of $\mathsf{R}_T:\DTM^\otimes_\vee\to \PMod_{H\KKK}^\otimes$. (Here $B$ is a commutative differential graded $\KKK$-algebra.)
Then $\mathsf{MTG}$ is obtained from the 
$\mathbb{G}_m$-equivariant bar construction of a commutative
differential graded $\KKK$-algebra $\overline{Q}$ equipped with
$\mathbb{G}_m$-action.
Namely, it is the \v{C}ech nerve of a morphism of derived
stacks $\Spec H\KKK\to [\Spec \overline{Q}/\mathbb{G}_m]$.
\end{TTheorem}

We remark that the underlying complex $\overline{Q}$
can be described in terms of Bloch's cycle complexes.
The proof of Theorem 2 requires two ingredients; one is Theorem 1, and
another is to
identify $\mathsf{R}_T:\DTM^\otimes_\vee\to \PMod_{H\KKK}^\otimes$
with a certain pullback functor between $\infty$-categories of perfect
complexes on derived stacks, which makes use of
the module-theoretic (i.e. Morita-theoretic)
presentation theorem of the stable $\infty$-category $\DTM_\vee^\otimes$, see
\cite{Sp1}.

If Beilinson-Soul\'e vanishing conjecture holds for the base field $k$
(e.g. $k$ is a number field), there is a traditional line
passing to a group scheme.
Under the vanishing conjecture, one can define
the motivic $t$-structure on $\DTM_\vee$.
The heart of this $t$-structure is a neutral Tannakian category
(cf. \cite{Sa}, \cite{DM}), and we can extract an affine group
scheme $MTG$ over $\KKK$ from it.
The so-called motivic Galois group for mixed Tate motives $MTG$ is constructed notably by Bloch-Kriz, Kriz-May, Levine \cite{BK}, \cite{KM}, \cite{Lev}.
The vanishing conjecture does not imply that the stable
$\infty$-category of complexes of the heart recovers the original
$\infty$-category $\DTM_\vee$. However,
we can describe a quite nice relation between $\mathsf{MTG}$ and $MTG$:

\begin{TTheorem}
Suppose that Beilinson-Soul\'e vanishing conjecture holds for $k$.
Then the group scheme $MTG$ is an excellent coarse moduli space
(cf. Definition~\ref{under}) of $\mathsf{MTG}$.
\end{TTheorem}

This result is proved in Section 7; Theorem~\ref{comparison}.
Roughly speaking, in this case the coarse moduli space
of $\mathsf{MTG}$
is obtained by truncating higher homotopy groups
of valued points of $\mathsf{MTG}$.
In view of Theorem 2 and 3, we can say that
the derived motivic Galois group constructed from $\DM^\otimes$
in \cite{Tan}
is a natural generalization of $MTG$ to the whole mixed motives.

In the final Section, for the sake of completeness
we will also treat the stable subcategory
of Artin motives in $\mathsf{DM}$, that is generated by motives
of smooth $0$-dimensional varieties.
We show that
the tannakization of the stable $\infty$-category of Artin motives
is the absolute Galois group $\Gal(\bar{k}/k)$ (see Proposition~\ref{toabsolute}).

\vspace{1mm}

{\it Application to future work.}
The result (Theorem~\ref{Thm1}) has already found nice applications
and has been applied in the work of
the study of motivic Galois group of mixed motives
generated by an abelian varieity (do not confuse it with 1-motives); see \cite{PM}.
In {\it loc. cit.}, it connects certain based loop stacks
with the representability
of automorphisms, which allows one to use various techniques such as
 Galois representations, rational homotopy theory, etc.
Combining with the general method of perfect adjoint pairs discussed
in \cite[Section 3]{PM} one may expect more to this and other directions.

\vspace{2mm}

This article is organized as follows:
In Section 2, we will review some of notions
and notation which we need
in this note.
In Section 3, we recall the
definitions of representation of
derived affine group schemes, automorphism
group of symmetric functors, etc.
Section 4 contains the proof of Theorem 1.
In Section 5, we give a brief exposition of bar constructions
from our viewpoint.
In Sections 6,7,8, we give applications to examples.
Sections 6 and 7 are devoted to the study of the tannakization of
stable $\infty$-category of mixed Tate motives; we prove Theorem 2 and 3.
In Section 8, we prove that the tannakization of
stable $\infty$-category of Artin motives
endowed with a realization functor
is the absolute Galois group.

\section{Notation and Convention}

We fix notation and convention.

{\it $\infty$-categories}.
In this note we use theory of quasi-categories as in \cite{Tan}.
A quasi-category is a simplicial set which
satisfies the weak Kan condition of Boardman-Vogt:
A quasi-category $S$ is a  simplicial set
such that for any $0< i< n$ and any diagram
\[
\xymatrix{
\Lambda^{n}_i \ar[r] \ar[d] & S \\
\Delta^n \ar@{..>}[ru] & 
}
\]
of solid arrows, there exists a dotted arrow filling the diagram.
Here $\Lambda^{n}_i$ is the $i$-th horn and $\Delta^n$
is the standard $n$-simplex.
Following \cite{HTT} we shall refer to quasi-categories
as $\infty$-categories.
Our main references are \cite{HTT}
 and \cite{HA}
(see also \cite{Jo}, \cite{DAGn}).
We often refer to a map $S\to T$ of $\infty$-categories
as a functor. We call a vertex in an $\infty$-category $S$
(resp. an edge) an object (resp. a morphism).
For the rapid introduction to $\infty$-categories, we
refer to \cite[Chapter 1]{HTT}, \cite{Gro}, \cite[Section 2]{FI}.
For the quick survey on various approaches to $(\infty,1)$-categories
and their relations,
we refer to \cite{Be}.

\begin{itemize}

\item $\Delta$: the category of linearly ordered finite sets (consisting of $[0], [1], \ldots, [n]=\{0,\ldots,n\}, \ldots$)

\item $\Delta^n$: the standard $n$-simplex

\item $\textup{N}$: the simplicial nerve functor (cf. \cite[1.1.5]{HTT})

\item $\mathcal{C}^{op}$: the opposite $\infty$-category of an $\infty$-category $\mathcal{C}$

\item Let $\mathcal{C}$ be an $\infty$-category and suppose that
we are given an object $c$. Then $\mathcal{C}_{c/}$ and $\mathcal{C}_{/c}$
denote the undercategory and overcategory respectively (cf. \cite[1.2.9]{HTT}).

\item $\operatorname{Cat}_\infty$: the $\infty$-category of small $\infty$-categories  in a fixed universe (cf. \cite[3.0.0.1]{HTT})

\item $\wCat$: $\infty$-category of $\infty$-categories

\item $\SSS$: $\infty$-category of small spaces (cf. \cite[1.2.16]{HTT})

\item $\textup{h}(\mathcal{C})$: homotopy category of an $\infty$-category (cf. \cite[1.2.3.1]{HTT})

\item $\Fun(A,B)$: the function complex for simplicial sets $A$ and $B$

\item $\Fun_C(A,B)$: the simplicial subset of $\Fun(A,B)$ classifying
maps which are compatible with
given projections $A\to C$ and $B\to C$.

\item $\Map(A,B)$: the largest Kan complex of $\Fun(A,B)$ when $A$ and $B$ are $\infty$-categories,

\item $\Map_{\mathcal{C}}(C,C')$: the mapping space from an object $C\in\mathcal{C}$ to $C'\in \mathcal{C}$ where $\mathcal{C}$ is an $\infty$-category.
We usually view it as an object in $\mathcal{S}$ (cf. \cite[1.2.2]{HTT}).

\end{itemize}

{\it Stable $\infty$-categories, symmetric monoidal $\infty$-categories and spectra}.
For the definitions of (symmetric) monoidal $\infty$-categories
and $\infty$-operads,
their algebra objects, we shall refer to \cite{HA}.
The theory of stable $\infty$-categories is developed in
\cite[Chapter 1]{HA}. We list some of notation.

\begin{itemize}

\item $\mathbb{S}$: the sphere spectrum

\item $\SP$: $\infty$-category of spectra, we denote the smash product by $\otimes$

\item $\PSP$ the full subcategory of $\SP$ spanned by compact spectra

\item $\Mod_A$: $\infty$-category of $A$-module spectra
for a commutative ring spectrum $A$

\item $\PMod_A$: the full subcategory of $\Mod_A$ spanned by compact objects
(in $\Mod_A$, an object is compact if and only if it is dualizable, see \cite{BFN}) .
We refer to objects in $\PMod_A$ as perfect $A$-module (spectra).

\item $\FIN$: the category of pointed finite sets $\langle 0\rangle_\ast=\{\ast\}, \langle 1\rangle_\ast=\{1, \ast\},\ldots,\langle n\rangle_\ast=\{1\ldots,n, \ast\},\ldots$. A morphism
is a map $f:\langle n\rangle_\ast\to \langle m \rangle_\ast$ such that
$f(\ast)=\ast$. Note that $f$ is not assumed to be order-preserving.

\item Let $\mathcal{M}^\otimes\to \mathcal{O}^\otimes$ be a fibration of $\infty$-operads. We denote by
$\Alg_{/\mathcal{O}^\otimes}(\mathcal{M}^\otimes)$ the $\infty$-category of algebra objects (cf. \cite[2.1.3.1]{HA}).  We often write $\Alg(\MMM^\otimes )$ or 
$\Alg(\MMM)$ for $\Alg_{/\mathcal{O}^\otimes}(\MMM^\otimes)$.
Suppose that $\mathcal{P}^\otimes\to \mathcal{O}^\otimes$ is a map of $\infty$-operads. $\Alg_{\mathcal{P}^\otimes/\mathcal{O}^\otimes}(\mathcal{M}^\otimes)$: $\infty$-category of $\mathcal{P}$-algebra objects.

\item $\CAlg(\mathcal{M}^\otimes)$: $\infty$-category of commutative
algebra objects in a symmetric
monoidal $\infty$-category $\mathcal{M}^\otimes\to \NNNN(\FIN)$.

\item $\CAlg_R$: $\infty$-category of commutative
algebra objects in the symmetric monoidal $\infty$-category $\Mod_R^\otimes$
where $R$ is a commutative ring spectrum. When $R=\mathbb{S}$, we set
$\CAlg=\CAlg_{\mathbb{S}}$.

\item $\Mod_A^\otimes(\mathcal{M}^\otimes)\to \NNNN(\FIN)$: symmetric monoidal
$\infty$-category of
$A$-module objects,
where $\mathcal{M}^\otimes$
is a symmetric monoidal $\infty$-category such that (1)
the underlying $\infty$-category admits a colimit for any simplicial diagram, and (2)
its tensor product functor $\mathcal{M}\times\mathcal{M}\to \mathcal{M}$
preserves
colimits of simplicial diagrams separately in each variable.
Here $A$ belongs to $\CAlg(\mathcal{M}^\otimes)$.
cf. \cite[3.3.3, 4.4.2]{HA}.

\end{itemize}

 Let $\mathcal{C}^\otimes$ be the symmetric monoidal $\infty$-category.
We usually denote, dropping the subscript $\otimes$,
by $\mathcal{C}$ its underlying $\infty$-category.
We say that an object $X$ in $\mathcal{C}$
is dualizable if there exist an object
$X^\vee$ and two morphisms $e:X\otimes X^\vee\to \mathsf{1}$ and
$c:\mathsf{1} \to X\otimes X^\vee$
with $\mathsf{1}$ a unit such that the composition
\[
X \stackrel{\textup{Id}_X\otimes c}{\longrightarrow} X\otimes X^\vee\otimes X \stackrel{e\otimes\textup{Id}_X}{\longrightarrow} X
\]
is equivalent to the identity, and 
\[
X^\vee \stackrel{c\otimes\textup{Id}_{X^\vee}}{\longrightarrow} X^\vee \otimes X\otimes X^\vee \stackrel{\textup{Id}_{X^\vee}\otimes e}{\longrightarrow} X^\vee
\]
is equivalent to the identity.
The symmetric monoidal structure of $\mathcal{C}$ induces
that of the homotopy category
$\textup{h}(\mathcal{C})$.
If we consider $X$ to be an object also
in $\textup{h}(\mathcal{C})$,
then $X$ is dualizable in $\mathcal{C}$ if and only if $X$
is dualizable in $\textup{h}(\mathcal{C})$.
For example, for $R\in \CAlg$, compact and dualizable objects
coincide in the symmetric monoidal $\infty$-category
$\Mod_R^\otimes$ (cf. \cite{BFN}).

Let us recall the symmetric monoidal $\infty$-categories $\widehat{\textup{Cat}}_\infty^{\textup{L,st}}$ and $\sCat$ (see \cite[Section 4]{BFN},
\cite[6.3]{HA}).
Let $\widehat{\textup{Cat}}_\infty^{\textup{L,st}}$ be the
subcategory of $\wCat$ spanned by stable presentable $\infty$-categories, in which morphisms are functors which preserves small colimits.
For $\mathcal{C},\mathcal{D}\in \widehat{\textup{Cat}}_\infty^{\textup{L,st}}$,
$\Fun^{\textup{L}}(\mathcal{C},\mathcal{D})$ is defined to be the
full subcategory of $\Fun(\mathcal{C},\mathcal{D})$
spanned by functors which preserves small colimits.
Then $\widehat{\textup{Cat}}_\infty^{\textup{L,st}}$
has a symmetric monoidal structure $\otimes:\widehat{\textup{Cat}}_\infty^{\textup{L,st}}\times \widehat{\textup{Cat}}_\infty^{\textup{L,st}}\to \widehat{\textup{Cat}}_\infty^{\textup{L,st}}$
such that for $\mathcal{C},\mathcal{D}\in \widehat{\textup{Cat}}_\infty^{\textup{L,st}}$, $\mathcal{C}\otimes \mathcal{D}$ has the following universality:
there exists a functor $\mathcal{C}\times \mathcal{D}\to \mathcal{C}\otimes\mathcal{D}$,
which induces an equivalence
$\Fun^{\textup{L}}(\mathcal{C}\otimes\mathcal{D},\mathcal{E})\simeq
\Fun'(\mathcal{C}\times\mathcal{D},\mathcal{E})$
for every $\mathcal{E}\in \widehat{\textup{Cat}}_\infty^{\textup{L,st}}$,
where the right hand side indicates the full subcategory
of $\Fun(\mathcal{C}\times\mathcal{D},\mathcal{E})$
spanned by functors which preserves small colimits separately in each variable.
A unit is equivalent to $\SP$.
Let $\sCat$
denote the subcategory of $\operatorname{Cat}_\infty$ which consists
of small stable idempotent complete $\infty$-categories.
Morphisms in $\sCat$
are functors that preserve finite colimits, that is,
exact functors.
There is a symmetric monoidal structure on $\sCat$.
For $\mathcal{C}, \mathcal{D}\in\sCat$ the
tensor product $\mathcal{C}\otimes\mathcal{D}$
has the following universality:
there is a functor $\mathcal{C}\times \mathcal{D}\to \mathcal{C}\otimes \mathcal{D}$ which preserves finite colimits separately in each variable, such that
if $\mathcal{E}\in \sCat$ and
$\Fun_{fc}(\mathcal{C}\times \mathcal{D},\mathcal{E})$
denotes the full subcategory of
$\Fun(\mathcal{C}\times \mathcal{D},\mathcal{E})$
spanned by functors which preserve finite colimits separately in each variable,
then the composition induces a categorical equivalence
\[
\Fun^{\textup{ex}}(\mathcal{C}\otimes\mathcal{D},\mathcal{E})\to \Fun_{fc}(\mathcal{C}\times \mathcal{D},\mathcal{E})
\]
where $\Fun^{\textup{ex}}(\mathcal{C}\otimes\mathcal{D},\mathcal{E})$
is the full subcategory of
$\Fun(\mathcal{C}\otimes\mathcal{D},\mathcal{E})$
spanned by exact functors.
A unit is equivalent to $\PSP$.
An object (resp. a morphism)
in $\CAlg(\widehat{\textup{Cat}}_\infty^{\textup{L,st}})$
can be regarded as 
a symmetric monoidal stable presentable $\infty$-category whose tensor operation preserves small colimits separately in each variable
(resp. a symmetric monoidal functor which preserves small colimits).
Similarly,
an object (resp. a morphism)
in $\CAlg(\sCat)$
can be regarded as 
a symmetric monoidal small stable idempotent complete $\infty$-category whose tensor opearation preserves finite colimits separately in each variable
(resp. a symmetric monoidal functor which preserves finite colimits).
If $R$ is a commutative ring spectrum,
we refer to an object in $\CAlg(\wCat^{\textup{L,st}})_{\Mod_R^\otimes/}$
(resp. $\CAlg(\textup{Cat}^{\textup{st}})_{\PMod_R^\otimes/}$)
simply as an $R$-linear symmetric monoidal stable presentable 
$\infty$-category (resp. an $R$-linear symmetric monoidal
small stable idempotent complete $\infty$-category) .
We refer to morphisms in $\CAlg(\wCat^{\textup{L,st}})_{\Mod_R^\otimes/}$
(or $\CAlg(\textup{Cat}^{\textup{st}})_{\PMod_R^\otimes/}$)
as $R$-linear symmetric monoidal functors.

\section{Derived group schemes and the $\infty$-categories of representations}

In this Section
we
recall the definitions of $\infty$-categories of representations of
derived affine group schemes and the tannakization of symmetric monoidal $\infty$-categories.

\subsection{Derived affine group scheme $G$ and $\infty$-categories $\Rep_G$ and $\PRep_G$.}
We refer to \cite[Appendix]{Tan} for the basic definitions concerning
derived group schemes.
Let $R$ be a commutative ring spectrum.
Let $G$ be a derived affine group scheme over $R$.
This can be viewed as a group object
$\psi:\NNNN(\Delta)^{op}\to \Aff_R:=(\CAlg_R)^{op}$
(see \cite[Definition A.2]{Tan}).
In this note, we refer to an object in $\Aff_R$ as an affine
 (derived) scheme over $R$
and call $\Aff_R$ the $\infty$-category
of affine (derived) schemes over $R$.
From Grothendick's viewpoint of ``functor of points'',
a derived affine group scheme over $R$ is a functor
$(\Aff_R)^{op} \to \Grp(\SSS)$ such that the composite
$(\Aff_R)^{op} \to \SSS$ with the forgetful functor $\Grp(\SSS)\to \SSS$
is represented by an affine scheme,
where $\Grp(\SSS)$ is the $\infty$-category of group objects in $\SSS$.
We will recall the definition of the symmetric monoidal
$\infty$-category $\Rep^\otimes_{G}$.
Set $G=\Spec B$ so that $B$ is a commutative Hopf ring spectrum
over $R$ which is described by a cosimplicial object
$\phi:=\psi^{op}:\NNNN(\Delta) \to \CAlg_R$. We here abuse notation and
$B$ indicates also the the underlying object $\phi([1])$
in $\CAlg_R$.
Let 
\[
\Theta:\CAlg \longrightarrow \CAlg(\wCat^{\textup{L,st}})
\]
be a functor which carries $A\in \CAlg$ to the symmetric
monoidal $\infty$-category $\Mod_A$
and sends a map $A\to A'$ in $\CAlg$
to a colimit-preserving symmetric monoidal base change functor
$\Mod_A\to \Mod_A': M\mapsto M\otimes_AA'$ (see \cite[Appendix A.6]{Tan}).
This functor induces
\[
\Theta_R:\CAlg_R\simeq \CAlg_{R/} \longrightarrow \CAlg(\wCat^{\textup{L,st}})_{\Mod_R^\otimes/}.
\]
Consider the composition $\NNNN(\Delta)\stackrel{\phi}{\to} \CAlg_R \stackrel{\Theta_R}{\to} \CAlg(\wCat^{\textup{L,st}})_{\Mod_R^\otimes/}$.
We define $\Rep^\otimes_G$ to be a limit of this composition.
We call it the symmetric monoidal
$\infty$-category of representations of $G$.
The underlying $\infty$-category is stable and presentable.
Since the forgetful functor $\CAlg(\wCat^{\textup{L,st}})_{\Mod_R^\otimes/}\to \wCat$
is limit-preserving,
we see that the underlying $\infty$-category of $\Rep_G^\otimes$,
which we denote by $\Rep_G$, is a limit of the composition
$\NNNN(\Delta)\stackrel{\Theta_R \circ \phi}{\longrightarrow} \CAlg(\wCat^{\textup{L,st}})_{\Mod_R^\otimes/}\to \wCat$.
There is the natural symmetric monoidal functor
$\Rep_G^\otimes\to \Mod_R^\otimes$ and we let $\PRep^\otimes_G$ be
the inverse image
of the full subcategory $\PMod_R^\otimes$.
Alternatively, there is a natural categorical equivalence
$\PRep_G\simeq \lim_{[n]\in \Delta}\PMod_{\phi([n])}$ and $\PRep_G^\otimes$ is a
 symmetric monoidal full subcategory of $\Rep^\otimes_G$
spanned by dualizable objects. We call it the symmetric monoidal
$\infty$-category of
perfect representations of $G$.

\subsection{$\infty$-categories of modules over presheaves.}
Let $(\CAlg_R)^{op} \hookrightarrow \Fun(\CAlg_R,\widehat{\mathcal{S}})$
be Yoneda embedding, where $\widehat{\mathcal{S}}$ denotes
the $\infty$-category of (not necessarily small) spaces, i.e. Kan complexes.
We shall refer to objects in $\Fun(\CAlg_R,\widehat{\mathcal{S}})$
as presheaves on $\CAlg_R$ or simply functors.
By left Kan extension of $\Theta_R$,
we have a colimit-preserving functor
\[
\overline{\Theta}_R:\Fun(\CAlg_R,\widehat{\mathcal{S}}) \to (\CAlg(\wCat)_{\Mod_R^\otimes/})^{op}.
\]
Let $\NNNN(\Delta)^{op}\stackrel{\psi}{\to} (\CAlg_R)^{op} \hookrightarrow \Fun(\CAlg_R,\widehat{\mathcal{S}})$ be the composition
and let $\mathsf{B}G$ denote the colimit.
Remember $\overline{\Theta}_R(\mathsf{B}G)=\Mod^\otimes_{\mathsf{B}G}\simeq \Rep_G^\otimes$.

Let $X \in \Fun(\CAlg_R,\widehat{\mathcal{S}})$.
Let $\PMod_X^\otimes$ denote
the symmetric monoidal full subcategory of the
underlying symmetric monoidal $\infty$-category $\overline{\Theta}_R$
spanned by dualizable objects.
Suppose that $\PMod_X^\otimes$ is a small stable idempotent complete
symmetric monoidal $\infty$-category whose tensor operation
$\otimes:\PMod_X\times \PMod_X\to \PMod_X$ preserves finite
colimits separately in each variable.
We refer to $\PMod_X^\otimes$ as the symmetric monoidal $\infty$-category
of perfect complexes on $X$.
We here call presheaves enjoying
this condition admissible presheaves (functors).
For example, affine derived schemes
and $\mathsf{B}G$ with $G$ a derived affine
group scheme are admissible. Indeed,
$\mathsf{B}G$ is described as the colimit of a simplicial affine derived
schemes $a:\NNNN(\Delta)^{op}\to \Aff_R$
and $\textup{Cat}_\infty^{\textup{st}}\hookrightarrow
 \textup{Cat}_\infty$ preserves small limits.
It follows that
$\PMod_{\mathsf{B}G}= \PRep_G \simeq \lim_{[n]}\PMod_{a([n])}$ is stable and idempotent complete
where $\lim_{[n]\in \Delta}\PMod_{a([n])}$
the limit of the cosimplicial diagram of
$\infty$-categories. 
Let $\Fun(\CAlg_R,\widehat{\mathcal{S}})^{\textup{adm}}$
be the full subcategory of $\Fun(\CAlg_R,\widehat{\mathcal{S}})$
spanned by admissible presheaves.
Applying $\overline{\Theta}_R$ and taking full subcategories of $\overline{\Theta}_R(X)$ spanned by dualizable objects we have the functor
\[
\overline{\theta}_R:\Fun(\CAlg_R,\widehat{\mathcal{S}})^{\textup{adm}} \to \CAlg(\textup{Cat}^{\textup{st}}_\infty)^{op}
\]
which carries $X$ to $\PMod_X^\otimes$.
We remark that by \cite[3.3.3.2, 5.1.2.2]{HTT}
$P$ in $\lim_{\Spec A\to X}\PMod_A$ ($\Spec A\to X$
run over $(\Aff_R)_{/X}$)
is a finite colimit of a (finite) diagram $I\to \PMod_X$
if and only if for each $\Spec A\to X$ the image of $P$ in $\PMod_A$ is a finite colimit of the induced diagram.

\subsection{Automorphisms.}

Let us review the automorphism group of a symmetric monoidal
functor.
Let $\mathcal{C}^\otimes$ be a symmetric monoidal small
$\infty$-category.
Let $\omega:\mathcal{C}^\otimes \to \PMod_R^\otimes$
be a symmetric monoidal functor.
We write $\mathcal{C}$ for its underlying $\infty$-category.
Let $\theta_{\mathcal{C}^\otimes}:\CAlg(\uCat)\to \mathcal{S}$ be the functor corresponding to $\mathcal{C}^\otimes$
via the Yoneda embedding $\CAlg(\uCat)^{op}\subset \Fun(\CAlg(\uCat),\mathcal{S})$. We denote the restriction of $\bar{\theta}_R$ to
$\Aff_R$ by $\theta_R$.
Then
the composite
\[
\xi:\CAlg_R \stackrel{\theta_R^{op}}{\longrightarrow} \CAlg(\uCat) \stackrel{\theta_{\mathcal{C}^\otimes}}{\longrightarrow} \mathcal{S}
\]
carries $A$ to the space equivalent to
$\Map^\otimes(\mathcal{C}^\otimes,\PMod_A^\otimes)$.
We can extends $\xi$ to
$\xi_\ast:\CAlg_R\to \mathcal{S}_{\ast}$
by using the symmetric monoidal functor $\omega$.
Here $\mathcal{S}_\ast$ denotes the $\infty$-category of pointed
spaces, that is, $\mathcal{S}_{\Delta^0/}$.
To explain this, let $\mathcal{M}\to \CAlg_R$ be a left fibration
corresponding to $\xi$.
An extension of $\xi$ to $\xi_\ast$ amounts to giving
a section $\CAlg_R\to \mathcal{M}$ of the
left fibration
$\mathcal{M}\to \CAlg_R$.
According to \cite[3.3.3.4]{HTT}
a section corresponds to an object
in the $\infty$-category $\mathcal{L}$
which is the limit of the diagram of spaces (or $\infty$-categories) given by $\xi$;
$A\mapsto \Map^\otimes(\mathcal{C}^\otimes, \PMod_A^\otimes)$.
Thus if $\lim\PMod_A^\otimes$ denotes the limit of
$\bar{\theta}_R^{op}:\CAlg_R\to \CAlg(\uCat)$, then
$\mathcal{L}$ is equivalent to $\Map^\otimes(\mathcal{C}^\otimes,\lim\PMod_A^\otimes)$ as $\infty$-categories (or equivalently spaces).
The natural functor $\PMod_R^\otimes \stackrel{\sim}{\to} \lim\PMod_A^\otimes$
induces $p:\Map^\otimes(\mathcal{C}^\otimes,\PMod_R^\otimes) \to \Map^\otimes(\mathcal{C}^\otimes,\lim\PMod_A^\otimes)\simeq \mathcal{L}$.
The image $p(\omega)$ in $\mathcal{L}$ gives rise to
a section $\CAlg_R\to \mathcal{M}$.
Consequently, we have $\xi_\ast:\CAlg_R\to \mathcal{S}_\ast$
which extends $\xi$.
We define $\Aut(\omega)$
to be the composite
$\CAlg_R\stackrel{\xi_\ast}{\longrightarrow} \mathcal{S}_\ast\stackrel{\Omega_\ast}{\longrightarrow} \Grp(\mathcal{S})$,
where the second functor is the based loop functor.
We refer to $\Aut(\omega)$ as the automorphism group functor of $\omega:\mathcal{C}^\otimes\to \PMod_R^\otimes$.

\subsection{Tannakization.}
We recall one of main results of \cite[Section 4]{Tan}.

\begin{Theorem}
\label{premain}
Let $\mathcal{C}^\otimes$ be a symmetric monoidal small $\infty$-category.
Let $\omega:\mathcal{C}^\otimes\to \PMod_R^\otimes$ be a symmetric monoidal
functor.
There exists a derived affine
group scheme $G$ over $R$ which represents
the automorphism group functor $\Aut(\omega)$.
Moreover, there is
a symmetric monoidal functor
$u:\mathcal{C}^\otimes\to \PRep_G^\otimes$ which makes the outer triangle in
\[
\xymatrix{
     &  \PMod_G^\otimes \ar@{..>}[d] \ar[ddr]^{\textup{forget}}  &     \\
   & \PRep_H^\otimes \ar[rd]_{\textup{forget}} &  \\
\mathcal{C}^\otimes \ar[ru] \ar[ruu]^{u} \ar[rr]_\omega & & \PMod_R^\otimes
}
\]
commute in the $\infty$-category of symmetric monoidal $\infty$-categories
such that these possess the following universality:
for any inner triangle consisting of solid arrows in the above diagram
where $H$ is a derived affine group scheme,
there exists a unique (in an appropriate sense)
morphism $f:H\to G$ of derived affine group
schemes which induces $\PRep_G^\otimes\to \PRep_H^\otimes$
(indicated by the dotted arrow) filling the above diagram.
\end{Theorem}

We usually refer to $G$ as the tannakization of
$\omega:\mathcal{C}^\otimes\to \PMod_R^\otimes$.
(In this note, we do not use this Theorem in an essential way.)

\section{Automorphism of fiber functors}

Let $Y$ be a derived stack over $R$ (we fix our convention below)
and $\PMod_Y^\otimes$ the
$\infty$-category of perfect complexes on $Y$,
which we regard as an object in $\CAlg(\sCat)_{\PMod^\otimes_R/}$.
Let $\Spec R\to Y$ be a section of the structure morphism $Y \to \Spec R$.
There is the pullback functor $\PMod_Y^\otimes\to \PMod_R^\otimes$
in $\CAlg(\sCat)$.
In this Section,
we study the automorphisms of this functor.
Our goal is Theorem~\ref{perfectkan} and Corollary~\ref{main}.

\vspace{2mm}

We start with our setup of derived stacks.
A functor $Y:\CAlg_R\to \widehat{\SSS}$
is said to be a derived stack (over $R$)
if two condition hold:
\begin{enumerate}
\renewcommand{\labelenumi}{(\roman{enumi})}
\item there exists a groupoid object
$\NNNN(\Delta)^{op}\to \Aff_R$
(cf. \cite[Definition A.2]{Tan}) such that $Y$ is equivalent to
the colimit of
the composite $\NNNN(\Delta)^{op}\to\Aff_R \hookrightarrow \Fun(\CAlg_R,\widehat{\SSS})$,
\item $Y$ has affine diagonal, that is, for any two morphisms $\Spec A\to Y$ and $\Spec B\to Y$, the fiber product $\Spec A\times_Y\Spec B$ belongs to $\Aff_R\subset \Fun(\CAlg_R,\widehat{\SSS})$.
\end{enumerate}
In this note, despite $Y$ in the above definition
is usually called a pre-stack,
we will not equip $\CAlg_R$
with Grothendieck topology such as flat, \'etale topologies
since the sheafification $Y'$ of $Y$ by such topologies does
induce a categorical equivalence
$\Mod_{Y'}\to \Mod_Y$ by the flat descent theory.
In addition, such topologies are irrelevant for our argument below.
(Conversely, for our purpose one can replace $\Fun(\CAlg_R,\widehat{\SSS})$ in the above definition
by the full subcategory
of sheaves with respect to flat topology (see e.g. \cite{HAG2}, \cite[VII, 5.4]{DAGn}
for flat morphisms)).
At any rate, we remark that
our definition of derived stacks is not standard (compare \cite{HAG2}, \cite{DAGn}). We note that our derived stacks are admissible functors.

\begin{Example}
\label{quotient}
We present quotient stacks arising from
the action of a derived affine group scheme on an
affine scheme as examples of derived stacks.
Let $F:\NNNN(\Delta)^{op}\to \Aff_R$ be a groupoid object, which
we regard as a derived stack.
Let $G:\NNNN(\Delta)^{op}\to \Aff_R$ be a group object, that is,
a derived affine group scheme.
Let $F\to G$ be a morphism (i.e., natural transformation)
which induces a cartesian diagram
\[
\xymatrix{
F([n]) \ar[r] \ar[d] & F([m]) \ar[d] \\
G([n]) \ar[r] & G([m])
}
\]
in $\Aff_R$ for each $[m]\to [n]$.
If we write $X$ for $F([0])$, then we can think that the morphism
$F\to G$ with the
above property means an action of $G$ on $X$.
In this situation, we say that $G$ acts on $X$ and denote by $[X/G]$
the colimit of $\NNNN(\Delta)^{op}\stackrel{F}{\to} \Aff_R\hookrightarrow \Fun(\CAlg_R,\widehat{\SSS})$. We refer to $[X/G]$ as the quotient stack.
We can think of $\mathsf{B}G$ as the quotient stack $[\Spec R/G]$
where $G$ acts trivially on $\Spec R$.
\end{Example}

Let $\pi:\Spec R\to Y$ denote the fixed section
and $\pi^*:\Mod_Y^\otimes\to \Mod_R^\otimes$ the associated
symmetric monoidal functor which preserves small colimits.
Since $\Mod_Y$ and $\Mod_R$ are presentable, by adjoint functor
theorem (see \cite[5.5.2.9]{HTT}) there is a right adjoint
functor $\pi_*:\Mod_R \to \Mod_Y$.
Moreover, according to \cite[8.3.2.6]{HA}
the right adjoint functor is extended to a right adjoint functor
to relative to $\NNNN(\FIN)$ (see \cite[8.3.2.2]{HA})
\[
\xymatrix{
\Mod_R^\otimes \ar[rd] \ar[rr] &  &  \Mod_Y^\otimes \ar[ld] \\
 & \NNNN(\FIN). & 
}
\]
It yields a right adjoint functor
\[
\CAlg(\Mod_R^\otimes)\to \CAlg(\Mod_Y^\otimes)
\]
of the
functor
$\CAlg(\Mod_Y^\otimes)\to \CAlg(\Mod_R^\otimes)$ determined by $\pi^*$.

Let $\phi:\NNNN(\Delta) \to \CAlg_R$ be a cosimplicial diagram
such that the colimit of composition
$\NNNN(\Delta)^{op} \stackrel{\phi^{op}}{\to} \Aff_R\hookrightarrow \Fun(\CAlg_R,\widehat{\SSS})$
is equivalent to $Y$.
Recall from Section 3.1 the functor $\Theta_R:\CAlg_R\to \CAlg(\widehat{\textup{Cat}}_\infty^{\textup{L,st}})_{\Mod_R^\otimes/}$.
Note that
by definition $\Mod_Y^\otimes$ is a limit of the composition $\phi'':\NNNN(\Delta)\stackrel{\phi}{\to} \CAlg_R \stackrel{\Theta_R}{\to} \CAlg(\widehat{\textup{Cat}}_\infty^{\textup{L,st}})_{\Mod_R^\otimes/} \to \CAlg(\widehat{\textup{Cat}}_\infty^{\textup{L,st}})$ where the last functor is the forgetful functor.
Let $p:\mathcal{M}_{\phi'} \to \NNNN(\Delta)$ be the coCartesian fibration
corresponding to the composition
$\phi':\NNNN(\Delta)\stackrel{\phi''}{\to} \CAlg(\widehat{\textup{Cat}}_\infty^{\textup{L,st}}) \to \wCat$ where the last functor is the forgetful functor.
We denote by $\Fun_{\NNNN(\Delta)}'(\NNNN(\Delta),\mathcal{M}_{\phi'})$
the full subcategory of $\Fun_{\NNNN(\Delta)}(\NNNN(\Delta),\mathcal{M}_{\phi'})$
spanned by sections $\NNNN(\Delta)\to \mathcal{M}_{\phi'}$
which carries all edges of $\NNNN(\Delta)$
to $p$-coCartesian edges.
Then by \cite[3.3.3.2]{HTT} $\Mod_Y$ is equivalent to
$\Fun_{\NNNN(\Delta)}'(\NNNN(\Delta),\mathcal{M}_{\phi'})$
as $\infty$-categories.
Consider the base change of $\NNNN(\Delta)^{op}\stackrel{\phi^{op}}{\to}\Aff_R \hookrightarrow \Fun(\CAlg_R,\widehat{\mathcal{S}})$, where the second functor is Yoneda embedding,
by $\pi:\Spec R\to Y$.
Let $Y_n=\phi^{op}([n])\in \Aff_R$ for each $[n]\in\Delta$.
The $n$-th term of this base change $\tau:\NNNN(\Delta)^{op}\to\Fun(\CAlg_R,\widehat{\mathcal{S}})$ is equivalent to
$Y_n\times_Y\Spec R$
and in particular, it factors through $\Aff_R\subset \Fun(\CAlg_R,\widehat{\mathcal{S}})$. Taking the opposite categories
we have $\psi:\NNNN(\Delta)\to \CAlg_R$.
Note that $\Spec R$ is a colimit of $\tau$ since in the $\infty$-topos
$\Fun(\CAlg_R,\widehat{\mathcal{S}})$ colimits are universal
(see \cite[Chapter 6]{HTT}).
Thus the natural transformation $\psi^{op}\to \phi^{op}$ induces $\pi:\Spec R\to Y$,
and we can informally indicates our situation as follows:
\[
\xymatrix{
   \cdots  \ar@<5pt>[r] \ar@<-5pt>[r] \ar[r] &  Y_1\times_Y\Spec R \ar[d] \ar@<3pt>[r] \ar@<-3pt>[r]  \ar@<2.5pt>[l] \ar@<-2.5pt>[l] &   Y_0\times_Y\Spec R \ar[l] \ar[d] \ar[r] &  \Spec R \ar[d]^{\pi}   \\
    \cdots \ar@<5pt>[r] \ar@<-5pt>[r] \ar[r] &   Y_1   \ar@<2.5pt>[l] \ar@<-2.5pt>[l]   \ar@<3pt>[r] \ar@<-3pt>[r]      &   Y_0 \ar[l] \ar[r] & Y
}
\]
(here $\psi^{op},\phi^{op}:\NNNN(\Delta)^{op}\to \Aff_R$).
We define $\psi':\NNNN(\Delta) \to \wCat$ in the same way that we define $\phi'$,
and we let $q:\mathcal{M}_{\psi'}\to \NNNN(\Delta)$
the coCartesian fibration corresponding to $\psi'$.
The natural transformation $\phi\to \psi$
corresponds to a map between coCartesian fibrations
$\mathcal{M}_{\phi'}\to \mathcal{M}_{\psi'}$
over $\NNNN(\Delta)$, which carries coCartesian edges to coCartesian edges.
Again by \cite[8.3.2.6]{HA}
there is a right adjoint functor 
$\mathcal{M}_{\psi'}\to \mathcal{M}_{\phi'}$
of $\mathcal{M}_{\phi'}\to \mathcal{M}_{\psi'}$
relative to $\NNNN(\Delta)$.
Let us observe the following:

\begin{Lemma}
The map $\mathcal{M}_{\psi'}\to \mathcal{M}_{\phi'}$ of coCartesian
fibrations over $\NNNN(\Delta)$ carries $q$-coCartesian edges
to $p$-coCartesian edges.
\end{Lemma}

\Proof
It suffices to show that
if for any map $r:[m]\to [n]$ in $\Delta$
we describe the diagram induced by $\psi^{op}\to \phi^{op}$ as
\[
\xymatrix{
Y_n\times_Y\Spec R \ar[r]^{a} \ar[d]^{b} &  Y_m\times_Y\Spec R \ar[d]^{c} \\
Y_n \ar[r]^{d} & Y_m,
}
\]
then the natural base change morphism $d^*\circ c_*\to b_* \circ a^*$
is an equivalence.
It follows from \cite[Lemma 3.14]{BFN}.
\QED

Let 
\[
\alpha:\Fun_{\NNNN(\Delta)}'(\NNNN(\Delta),\mathcal{M}_{\phi'})\rightleftarrows \Fun_{\NNNN(\Delta)}'(\NNNN(\Delta),\mathcal{M}_{\psi'}):\beta
\]
be functors induced by the adjunction
$\mathcal{M}_{\phi'}\rightleftarrows \mathcal{M}_{\psi'}$, where $\Fun_{\NNNN(\Delta)}'(\NNNN(\Delta),\mathcal{M}_{\phi'})$
is the full subcategory of
$\Fun_{\NNNN(\Delta)}(\NNNN(\Delta),\mathcal{M}_{\phi'})$,
spanned by sections which carries all edges to coCartesian edges
and we define $\Fun_{\NNNN(\Delta)}'(\NNNN(\Delta),\mathcal{M}_{\phi'})$
in a similar way.
Note that by \cite[3.3.3.2]{HTT}
\[
\Fun_{\NNNN(\Delta)}'(\NNNN(\Delta),\mathcal{M}_{\phi'})\simeq \Mod_{Y}\ \ 
\textup{and}\ \ \Fun_{\NNNN(\Delta)}'(\NNNN(\Delta),\mathcal{M}_{\psi'})\simeq \Mod_{R},
\]
and $\Fun_{\NNNN(\Delta)}'(\NNNN(\Delta),\mathcal{M}_{\phi'})\to \Fun_{\NNNN(\Delta)}'(\NNNN(\Delta),\mathcal{M}_{\psi'})$
is equivalent to $\pi^*:\Mod_Y \to \Mod_R$ as functors.
Then observe that the pair $(\alpha,\beta)$ forms adjunction.
Namely, 
\begin{eqnarray*}
\Map_{\Fun_{\NNNN(\Delta)}'(\NNNN(\Delta),\mathcal{M}_{\psi'})}(\alpha(a),b)&\simeq& \lim_{[n]\in \Delta}\Map_{\psi'([n])}(\alpha(a_n),b_n) \\
&\to&  \lim_{[n]\in \Delta}\Map_{\phi'([n])}(\beta(\alpha(a_n)),\beta(b_n))  \\
&\stackrel{x}{\to}& \lim_{[n]\in \Delta}\Map_{\phi'([n])}(a_n,\beta(b_n)) \\
&\simeq& \Map_{\Fun_{\NNNN(\Delta)}'(\NNNN(\Delta),\mathcal{M}_{\phi'})}(a,\beta(b))
\end{eqnarray*}
is equivalence in $\mathcal{S}$,
where $a_n$ (resp. $b_n$) is the projection of $a$ (resp. $b$)
to $\phi'([n])$ (resp. $\psi'([n])$)
and $x$ is induced by the unit map of the adjunction
$\mathcal{M}_{\phi'}\rightleftarrows \mathcal{M}_{\psi'}$.
(The fiber of the adjunction
$\mathcal{M}_{\phi'}\rightleftarrows \mathcal{M}_{\psi'}$
over each object of $\NNNN(\Delta)$ forms adjunction.)
Notice that $\Fun_{\NNNN(\Delta)}(\NNNN(\Delta),\mathcal{M}_{\psi'})\to \Fun_{\NNNN(\Delta)}(\NNNN(\Delta),\mathcal{M}_{\phi'})$ is equivalent to $\pi_*:\Mod_R\to \Mod_Y$ as functors.
Consequently, we have

\begin{Lemma}
\label{bbb}
Let
\[
\xymatrix{
Y_n\times_Y\Spec R \ar[r]^{s_n} \ar[d]^{\pi_n} & \Spec R \ar[d]^\pi \\
Y_n \ar[r]^{t_n} & Y
}
\]
be the pullback diagram induced by $\psi^{op}([n])\to \phi^{op}([n])$.
Then the natural base change morphism
$(t_n)^*\circ \pi_*\to (\pi_n)_*\circ (s_n)^*$ is an equivalence
of functors from $\Mod_R$ to $\Mod_{Y_n}$.
\end{Lemma}

\begin{Corollary}
\label{aaa}
We abuse notation and 
we write
$(t_n)^*\circ \pi_*\to (\pi_n)_*\circ (s_n)^*$
for the natural base change morphism from $\CAlg(\Mod_R^\otimes)$ to $\CAlg(\Mod_{Y_n}^\otimes)$
which is determined by adjunctions $(\pi^*,\pi_*)$ and $((\pi_n)^*,(\pi_n)_*)$
relative to $\NNNN(\FIN)$.
Then $(t_n)^*\circ \pi_*\to (\pi_n)_*\circ (s_n)^*$ is an equivalence of
functors.
\end{Corollary}

Let $\mathbf{1}_R$ be a unit of $\Mod_R$ which we here regard
as an object in $\CAlg_R=\CAlg(\Mod_R)$.
Then
there is a lax symmetric monoidal functor $\Mod_R^\otimes \to \Mod_{\pi_*\mathbf{1}_R}^\otimes(\Mod_{Y}^\otimes)$ of symmetric monoidal $\infty$-categories
induced by $\pi_*$ by the construction of the
$\infty$-operad of module objects
\cite[3.3.3.8]{HA}.
For the notation $\Mod_{\pi_*\mathbf{1}_R}^\otimes(\Mod_{Y}^\otimes)$,
see Section 2.

\begin{Lemma}

\label{fs}

The functor
$\Mod_R^\otimes \to \Mod_{\pi_*\mathbf{1}_R}^\otimes(\Mod_{Y}^\otimes)$
is a symmetric monoidal equivalence.
\end{Lemma}

\Proof
We first obeserve that $\Mod_R^\otimes \to \Mod_{\pi_*\mathbf{1}_R}^\otimes(\Mod_{Y}^\otimes)$ is symmteric monoidal.
Since it is lax symmetric monoidal, combined with
Lemma~\ref{bbb} we are reduced to showing the following obvious
claim:
for a morphism $x:\Spec A\to \Spec B$ of affine derived schemes
and $M,N\in \Mod_A$, the natural map
$x_*(M)\otimes_{A}x_*(N)\to x_*(M\otimes_AN)$ is an equivalence
where $x_*:\Mod_A\to \Mod_A(\Mod_B^\otimes)$ is the natural pushforward
functor.

We now adopt notation similar to Lemma~\ref{bbb}.
Since the natural equivalence
$(t_n)^*\circ \pi_*\mathbf{1}_R \simeq (\pi_n)_*\circ (s_n)^*\mathbf{1}_R$
by the above result,
we have
\[
(\pi_n)_*:\Mod_{\psi([n])}=\Mod_{Y_n\times_Y\Spec R}\simeq \Mod_{(\pi_n)_*\circ(s_n)^*\mathbf{1}_R}(\Mod^\otimes_{\phi([n])})\simeq \Mod_{(t_n)^*\circ \pi_*\mathbf{1}_R}(\Mod^\otimes_{\phi([n])})
\]
for each $n$.
Then we identify $\Mod_R \to \Mod_{\pi_*\mathbf{1}_R}(\Mod_{Y}^\otimes)$
with the limit
\[
\lim_{[n]\in \Delta}\Mod_{\psi([n])} \simeq \lim_{[n]\in \Delta} \Mod_{Y_n\times_Y\Spec R}\simeq \lim_{[n]\in \Delta}\Mod_{(t_n)^*\circ \pi_*\mathbf{1}_R}(\Mod^\otimes_{\phi([n])})
\]
which is an equivalence in $\wCat$. It follows that
$\Mod_R^\otimes \to \Mod_{\pi_*\mathbf{1}_R}^\otimes(\Mod_{Y}^\otimes)$
is a symmetric monoidal equivalence.
\QED

Let $\Aut(\pi^*):\CAlg_R\to \Grp(\widehat{\mathcal{S}})$
be the automorphism group functor of $\pi^*$ (defined
as in the previous Section),
which carries $A\in \CAlg_R$ to the automorphisms of
composition $\Mod_Y^\otimes\to \Mod_R^\otimes\to \Mod_A^\otimes$
in $\CAlg(\wCat^{\textup{L,st}})$ where
the second functor is the base change by $R\to A$.

Let $\Delta_+$ be the category of finite (possibly empty)
linearly ordered sets and we write $[-1]$ for the empty set.
Let $\iota:\Delta^1\to \NNNN(\Delta_+)$
be a map which carries $\{0\}$ and $\{1\}$
to $[-1]$ and $[0]$ respectively.
It is a fully faithful functor.
Let $(\Delta^1)^{op}\to \Fun(\CAlg_R,\widehat{\SSS})$ be a map
corresponding to $\pi:\Spec R\to Y$.
Let $\rho:\NNNN(\Delta_+)^{op}\to \Fun(\CAlg_R,\widehat{\SSS})$
be a right Kan extension along $\iota^{op}:(\Delta^1)^{op}\to \NNNN(\Delta_+)^{op}$ which is called \v{C}ech nerve (cf. \cite[6.1.2.11]{HTT}).
By our assumption, for each $n\ge 0$, $\rho([n])$ belongs to $\Aff_R$
and the restriction of $\rho$ to $\NNNN(\Delta)^{op}$ is a derived affine
group scheme which we denote by $G_\pi$.
By the definition of $G_\pi$ and $\Mod^\otimes_{G_\pi}$,
we see that $\pi^*:\Mod_Y^\otimes\to \Mod_R^\otimes$ factors through
the forgetful functor
$\Rep_{G_\pi}^\otimes\to \Mod_R^\otimes$.
Note that
the derived group scheme $G_\pi:(\Aff_R)^{op}\to \Grp(\SSS)$ represents the automorphism group $\Aut(\pi):\CAlg_R\to \Grp(\mathcal{S})$
of $\pi:\Spec R\to Y$.
Here for any $A\in \CAlg_R$, $\Aut(\pi)(A)$ is the mapping space in
$\Map_{\Fun(\CAlg_R,\widehat{\mathcal{S}})}(\Spec A,Y)$
from $\Spec A\to \Spec R\stackrel{\pi}{\to} Y$ to itself, endowed
with the group structure (the construction
is similar to that of $\Aut(\omega)$ in the previous Section).
We have the natural morphism
$G_\pi\simeq \Aut(\pi)\to \Aut(\pi^*)$.

\begin{Proposition}
\label{naturalrep}
The natural morphism $G_\pi\to \Aut(\pi^*)$ is an equivalence,
that is, $\Aut(\pi^*)$ is represented by $G_\pi$.
\end{Proposition}

\Proof
For simplicity, let $G:=G_\pi$.
Let $G_1:\CAlg_R\to \widehat{\SSS}$ and (resp. $\Aut(\pi^*)_1$)
be the composite of $G:\CAlg_R\to \Grp(\widehat{\SSS})$ (resp. $\Aut(\pi^*)$)
and the forgetful functor $\Grp(\widehat{\SSS})\to \widehat{\SSS}$.
For each $A\in \CAlg_R$, it will suffice to show that
the induced map 
$G_1(A)\to \Aut(\pi^*)_1(A)$ is an equivalence in $\widehat{\SSS}$.

For $A\in \CAlg_R$, let $\pi_A:\Spec A\to \Spec R\to Y$
denote the composition.
Let $\mathbf{1}_A$ be the unit of $\Mod_A$ which we here think of as an object
of $\CAlg(\Mod_A^\otimes)$.
Applying \cite[6.3.5.18]{HA} together with Lemma~\ref{fs} and adjunction
we deduce
\begin{eqnarray*}
\Map_{\CAlg(\widehat{\textup{Cat}}_\infty^{\textup{L,st}})_{\Mod_Y^\otimes/}}(\Mod_A^\otimes,\Mod_A^\otimes) &\simeq& \Map_{\CAlg(\Mod_Y^\otimes)}((\pi_A)_*\mathbf{1}_A,(\pi_A)_*\mathbf{1}_A)   \\
&\simeq& \Map_{\CAlg(\Mod_A)}((\pi_A)^*(\pi_A)_*\mathbf{1}_A, \mathbf{1}_A).
\end{eqnarray*}
Unwinding the definitions
we have
\begin{eqnarray*}
\Map_{\CAlg(\Mod_A)}((\pi_A)^*(\pi_A)_*\mathbf{1}_A, \mathbf{1}_A) &\simeq& \Map_{(\Aff)_{/\Spec A}}(\Spec A,\Spec A\times_Y\Spec A) \\
&\simeq& \Map_{(\Aff)_{/\Spec A}}(\Spec A,G_1\times_R A\times_RA) \\
&\simeq& \Map_{\Aff_{/Y}}(\Spec A,\Spec A)
\end{eqnarray*}
where $G_1$ is $\Spec R\times_Y\Spec R \simeq \rho([1])$, and
$G_1\times_RA\times_RA\to \Spec A \in (\Aff)_{/\Spec A}$
is the second projection.
Note that through natural equivalences
a morphism $\Spec A\to \Spec A$ over $Y$, which we regard
as an object of $\Map_{\Aff_{/Y}}(\Spec A,\Spec A)$,
induces a symmetric monoidal functor $\Mod_A^\otimes \to \Mod_A^\otimes$
under $\Mod_{Y}^\otimes$ which we think of as an object of 
$\CAlg(\widehat{\textup{Cat}}_\infty^{\textup{L,st}})_{\Mod_Y^\otimes/}$.

Next using the natural equivalence
\[
\Map_{\CAlg(\widehat{\textup{Cat}}_\infty^{\textup{L,st}})_{\Mod_Y^\otimes/}}(\Mod_A^\otimes,\Mod_A^\otimes) \simeq \Map_{\Aff_{/Y}}(\Spec A,\Spec A)
\]
we consider the automorphisms of $\pi^*$.
To this end let
$T_A$ be the fiber product
\[
\Map_{\Aff_{/Y}}(\Spec A,\Spec A)\times_{\Map_{\Aff}(\Spec A,\Spec A)}\{\textup{Id}_{\Spec A}\}
\]
in $\widehat{\SSS}$ where the diagram is induced by the forgetful functor
$\Map_{\Aff_{/Y}}(\Spec A,\Spec A)\to\Map_{\Aff}(\Spec A,\Spec A)$.
Similarly, we define $S_A$ to be th fiber product
\[
\Map_{\CAlg(\widehat{\textup{Cat}}_\infty^{\textup{L,st}})_{\Mod_Y^\otimes/}}(\Mod_A^\otimes,\Mod_A^\otimes)\times_{\Map_{\CAlg(\widehat{\textup{Cat}}_\infty^{\textup{L,st}})}(\Mod_A^\otimes,\Mod_A^\otimes)}\{\textup{Id}\}
\]
in $\widehat{\SSS}$, which is equivalent to $T_A$.
There are natural equivalences
\begin{eqnarray*}
T_A &\simeq& \Map'_{(\Aff)_{/\Spec A}}(\Spec A,G_1\times_R A\times_RA) \\
&\simeq& \Map_{(\Aff)_{/\Spec A}}(\Spec A,G_1\times_R A) \\
&\simeq& \Map_{\Aff}(\Spec A,G_1)
\end{eqnarray*}
in $\widehat{\mathcal{S}}$
where $\Map'_{(\Aff)_{/\Spec A}}(\Spec A,G_1\times_R A\times_RA)$
is the fiber product 
\[
\Map_{(\Aff)_{/\Spec A}}(\Spec A,G_1\times_R A\times_RA)\times_{\Map_{\Aff}(\Spec A,\Spec A)}\{\textup{Id}_{\Spec A}\}
\]
in $\widehat{\SSS}$ where the diagram is induced by
the projection $\textup{pr}_3:G_1\times_RA\times_RA \to \Spec A$.
Thus we have an equivalence $\Map_{\Aff}(\Spec A,G_1)\simeq S_A$.
Hence we have the required equivalence
$G_1(A)\simeq \Aut(\pi^*)_1(A)$.
\QED

Let $\mathcal{C}^\otimes, \mathcal{D}^\otimes \in \CAlg(\wCat^{\textup{L,st}})$.
Suppose that
$\mathcal{C}$ is compactly generated, that is, the natural colimit-preserving functor
$\Ind(\mathcal{C}_\circ)\to \mathcal{C}$
is a categorical equivalence,
and $\otimes:\mathcal{C}\times\mathcal{C}\to \mathcal{C}$
induces
$\mathcal{C}_\circ\times\mathcal{C}_\circ\to \mathcal{C}_\circ$,
which makes $\mathcal{C}_\circ$ a symmetric monoidal $\infty$-category,
where $\mathcal{C}_\circ$ is the full subcategory of compact objects
in $\mathcal{C}$ and $\Ind(-)$ indicates the Ind-category (see \cite[5.3.5]{HTT}). Note that under this assumption, a unit object is compact.
Recall the following
result which follows from \cite[5.3.6.8]{HTT} and \cite[6.3.1.10]{HA}.

\begin{Proposition}
\label{symKan}
Let $\Map^{\otimes,\textup{L}}(\mathcal{C}^\otimes,\mathcal{D}^\otimes)$
be $\Map_{\CAlg(\wCat^{\textup{L,st}})}(\mathcal{C}^\otimes,\mathcal{D}^\otimes)$.
Let
$\Map^{\otimes,\textup{ex}}(\mathcal{C}_\circ^\otimes,\mathcal{D}^\otimes)$
be the full subcategory of $\Map_{\CAlg(\wCat)}(\mathcal{C}_\circ^\otimes,\mathcal{D}^\otimes)$ spanned by symmetirc monoidal functors which preserves
finite colimits.
The natural inclusion
$\mathcal{C}_\circ^\otimes \to \mathcal{C}^\otimes$ induces
an equivalence
\[
\Map^{\otimes,\textup{L}}(\mathcal{C}^\otimes,\mathcal{D}^\otimes)\to \Map^{\otimes,\textup{ex}}(\mathcal{C}_\circ^\otimes,\mathcal{D}^\otimes)
\]
in $\widehat{\mathcal{S}}$.
\end{Proposition}

Let us recall the definition of perfectness of stacks
introduced by Ben-Zvi, Francis, and Nadler
in their work on derived Morita theory \cite{BFN}
(this notion is also important to our previous paper \cite{FI}).
We say that a derived stack $Y$ is perfect
if the natural functor
$\Ind(\PMod_Y)\to \Mod_Y$ is a categorical
equivalence.
As a corollary of results of this Section, we have:

\begin{Theorem}
\label{perfectkan}
Let $Y$ be a perfect derived stack over $R$ and $\pi:\Spec R\to Y$ is a
section of the structure morphism $Y\to\Spec R$.
Let $\pi^*:\Mod_Y^\otimes \to \Mod_R^\otimes$ be the morphism
in $\CAlg(\wCat^{\textup{L,st}})$ induced by $\pi:\Spec R\to Y$,
and let $\pi_\circ^* :\PMod_Y^\otimes \to \PMod_R^\otimes$
denote its restriction which belongs to $\CAlg(\textup{Cat}_\infty^{\textup{st}})$.
Let $\Aut(\pi_\circ^*):\CAlg_R\to \Grp(\mathcal{S})$ be the automorphism
functor of $\pi_\circ^*$.
Then the restriction
induces an equivalence of functors $\Aut(\pi^*)\to \Aut(\pi^*_\circ)$.
In particular, the tannakization of
$\pi^*_\circ:\PMod_Y^\otimes\to \PMod_R^\otimes$ is equivalent to $G_\pi$.
(see the setup
before Proposition~\ref{naturalrep} for the notation $G_\pi$.)
\end{Theorem}

\Proof
Combine Proposition~\ref{naturalrep} and~\ref{symKan}.
\QED

\begin{Corollary}
\label{main}
Let $Y$ be a derived stack over $R$ equipped with $\pi:\Spec R\to Y$ as in
Theorem~\ref{perfectkan}.
Suppose either one of followings:
\begin{enumerate}
\renewcommand{\labelenumi}{(\roman{enumi})}
\item a derived stack $Y$ over $R$ belongs to $\Aff_R$,

\item let $G$ be an affine group scheme of finite type over a field
$k$ of characteristic zero, which we regard as a derived affine group
scheme over $R=Hk$. Suppose that $G$ acts on $X\in \Aff_R$
and let $Y=[X/G]$ be the quotient stack (see Example~\ref{quotient}).
\end{enumerate}
Then the tannakization of
$\pi_\circ^*:\PMod_Y^\otimes \to \PMod_R^\otimes$ is equivalent to $G_\pi$.
\end{Corollary}

\Proof
According to Proposition~\ref{naturalrep}
and Theorem~\ref{perfectkan},
it will suffices to show that $Y$ is perfect, that is,
the
natural functor $\Ind(\PMod_Y)\to \Mod_Y$
is a categorical equivalence.
Then our claim follows from \cite[3.19, 3.22]{BFN}.
\QED

\section{Bar constructions}

This Section contains no new result.
In this Section, we review the relation between
bar constructions and the case (i) of Corollary~\ref{main}.
Let $A \in \CAlg_R$ and let $s:R\to A$
be the natural morphism in $\CAlg_R$
(note $R$ is an initial object in $\CAlg_R$).
Suppose that $t:A \to R$ is a cosection
of $s$, that is, $t \circ s$ is equivalent to the identity of $R$.
Recall that $\Delta_+$ is the category of finite (possibly empty)
linearly ordered sets and we write $[-1]$ for the empty set.
Let $\iota:\Delta^1\to \NNNN(\Delta_+)$
be a map which carries $\{0\}$ and $\{1\}$
to $[-1]$ and $[0]$ respectively.
It is a fully faithful functor.
Let $f:\Delta^1\to \CAlg_R$ be the map corresponding to
$A \to R$.
Since $\CAlg_R$ admits small colimits,
there is a left Kan extension
\[
g:\NNNN(\Delta_+)\to \CAlg_R
\]
of $f$ along $\iota$.
We refer to $g^{op}:\NNNN(\Delta_+)^{op} \to \Aff_R$
as the \v{C}ech nerve of $f^{op}:(\Delta^1)^{op}\to \Aff_R$.
This construction is called the bar construction for $t:A\to R$.
The underlying simplicial object $\NNNN(\Delta)^{op} \to \NNNN(\Delta_+)^{op}\to \Aff_R$ is a group object (see \cite[Appendix]{Tan} or \cite[7.2.2.1]{HTT}
for the definition of group objects).
Let $G$ be a derived affine group scheme
corresponding to the simplicial object.

Let $t_\circ^*:\PMod_A^\otimes \to \PMod_R^\otimes$
be the morphism in $\CAlg(\sCat)_{\PMod_R^\otimes/}$.
The case (i) of Corollary~\ref{main} says:

\begin{Theorem}
\label{bar}
$\Aut(t_\circ^*)$ is represented by $G$.
\end{Theorem}

\begin{Remark}
For the readers who are familiar with commutative differential graded 
algebras (dg-algebras for short),
we relate the bar construction of commutative dg-algberas with $G$.
Let $k$ be a field of characteristic zero.
Let $\mathsf{dga}_k$ be the category of commutative dg-algebras
over $k$ (cf. \cite{Hi}). A morphism $P^\bullet\to Q^\bullet$ in $\dga_k$
is a weak equivalence (resp. fibration) if it induces
a bijection $H^n(P^\bullet)\to H^n(Q^\bullet)$ for each $n\in \ZZ$
(resp. $P^n\to Q^n$ is a surjective morphism
of $k$-vector spaces for each $n\in \ZZ$).
There is a model category structure on $\dga_k$
whose weak equivalences and fibrations are defined 
in this way (see \cite[2.2.1]{Hi}).
Let $\NNNN(\dga_k^c)_\infty$ be the $\infty$-category
obtained from the full subcategory $\dga_k^c$ spanned by
cofibrant objects by inverting weak equivalences (see \cite[1.3.4.15]{HA}).
According to \cite[8.1.4.11]{HA}, there is
a categorical equivalence
$\NNNN(\dga_k^c)_\infty\simeq \CAlg_{Hk}$.
Let $R=Hk$ and let $t:A\to k$ be an augmentation in $\dga_k$.
We abuse notation and we denote by $t:A\to R$ the induced morphism
in $\CAlg_R$.
The underlying derived scheme of $G$ is the fiber
product $\Spec R\times_{\Spec A}\Spec R$ in $\Aff_R$.
By this equivalence, the pushout $R\otimes_AR$ in $\CAlg_R$
corresponds to a homotopy pushout
$k\otimes^{\mathbb{L}}_Ak$ in the model category $\dga_k$,
which is weak equivalent to a homotopy pushout
$A\otimes^{\mathbb{L}}_{A\otimes_kA}k$
of
\[
\xymatrix{
A\otimes_kA\ar[r]^{t\otimes t} \ar[d]^{m} & k \\
A
}
\]
where $m$ is the multiplication.
We will review the construction of the concrete model of a homotopy pushout
$A\otimes^{\mathbb{L}}_{A\otimes_kA}k$ in $\dga_k$, which is known
as the bar construction of a commutative dg-algebra
(see for example \cite{Ol}, \cite{Te}).
Consider the adjoint pair
\[
T:\dga_{k,A/}\rightleftarrows \dga_{k,A\otimes_kA/}:U
\]
where $U$
is the forgetful functor induced by $A\to A\otimes_kA,\ x\mapsto x\otimes 1$,
and
$T$ is given by formula $M\mapsto M\otimes_A(A\otimes_kA)$.
Let $\alpha:\textup{Id}\to UT$ and $\beta:TU\to \textup{Id}$
be the unit map and counit map respectively.
To an object $C\in \dga_{k,A\otimes_kA/}$ one associates
a simplicial diagram $(T,U)_\bullet(C)$ in $\dga_{k,A/}$ as follows:
Define
\[
(T,U)_n(C)=(TU)^{\circ(n+1)}(C)=(TU)\circ\cdots \circ (TU)(C)
\]
where the right hand side is the $(n+1)$-fold composition.
For $0\le i\le n+1$,
\begin{eqnarray*}
d_i:(T,U)_{n+1}(C)&=&(TU)^{\circ i}\circ (TU)\circ (TU)^{\circ (n+1-i)}(C) \\
&\to& (TU)^{\circ i}\circ \textup{Id}\circ (TU)^{\circ (n+1-i)}(C)=(T,U)_{n}(C)
\end{eqnarray*}
is induced by $\beta$ in the middle term.
For $0\le i\le n$,
\begin{eqnarray*}
s_i:(T,U)_{n}(C)&=&(TU)^{\circ i}\circ T\circ \textup{Id} \circ U \circ (TU)^{\circ (n-i)}(C) \\
&\to& (TU)^{\circ i}\circ T\circ (UT)\circ U \circ (TU)^{\circ (n-i)}(C)=(T,U)_{n+1}(C)
\end{eqnarray*}
is induced by $\alpha:\textup{Id} \to (UT)$ in the middle term.
Let us consider $A$ to be an object in $\dga_{k,A\otimes_kA/}$ via $m:A\otimes_kA\to A$. Then by the above construction we obtain the simplicial object
$(T,U)_\bullet(A)\otimes_{A\otimes_kA}k$ in $\dga_k$.
The totalization
$\textup{tot}((T,U)_\bullet(A)\otimes_{A\otimes_kA}k)\in \dga_k$, which we call
the bar complex, represents the homotopy pushout
$A\otimes^{\mathbb{L}}_{A\otimes_kA}k$.
\end{Remark}


\section{Mixed Tate motives}

In this Section, as an application of the results we have proved; in particular Theorem~\ref{perfectkan} and Corollary~\ref{main}, we will describe the tannakization of the
stable $\infty$-category of mixed Tate motives equipped with
the realization functor
as the $\mathbb{G}_m$-equivariant bar construction
of a commutative dg-algebra.
The main goal of this Section is Theorem~\ref{tatemain}.
We emphasize that this section works without assuming Beilinson-Soul\'e vanishing conjecture. In what follows we often use model categories.
Our references for them are \cite{Ho1} and \cite[Appendix]{HTT}.

\subsection{Review of $\infty$-category of mixed motives}
Let $\mathbf{K}$ be a field of characteristic zero.
Let $\mathcal{A}$ be
the abelian category of $\mathbf{K}$-vector spaces.
We equip the category of complexes of $\mathbf{K}$-vector spaces, denoted by
$\textup{Comp}(\mathcal{A})$, with the projective model structure,
in which weak equivalences are quasi-isomorphisms,
and fibrations are degreewise surjective maps (cf. e.g. \cite[Section 2.3]{Ho1}, \cite[Appendix]{HTT}, \cite{CD1}).

Let $k$ be a perfect field.
Let $\textup{DM}^{eff}(k)$
be the category of complexes of $\mathcal{A}$-valued Nisnevich sheaves
with transfers (the indroductory references of this notion include
\cite{MVW} and \cite{Deg}). For a smooth scheme $X$ separated of finite type
over $k$, we denote by $L(X)$ the 
$\mathcal{A}$-valued Nisnevich sheaves
with transfers which is represented by $X$ (cf. \cite[page.15]{MVW}).
We equip $\textup{DM}^{eff}(k)$ with the symmetric monoidal model structure
in \cite[Example 4.12]{CD1}.
The triangulated subcategory of the homotopy category of this model
category $\textup{DM}^{eff}(k)$, spanned by right bounded complexes,
is equivalent to the triangulated category $\mathbf{DM}^{\textup{eff},-}_{Nis}(k,\mathbf{K})$ constructed in
\cite[Lecture 14]{MVW}.

The pointed algebraic torus
$\Spec (k)\to \mathbb{G}_m$ over $k$ induces
a split monomorphism $L(\Spec (k))\to L(\mathbb{G}_m)$
in $\textup{DM}^{eff}(k)$.
Then
we define $\mathbf{K}(1)$ to be
\[
\textup{Coker}(L(\Spec (k))\to L(\mathbb{G}_m))[-1].
\]
Let $\operatorname{DM}(k)$ be
the category of symmetric $\mathbf{K}(1)$-spectra in $(\textup{DM}^{eff}(k))^\mathfrak{S}$ (cf. \cite[Section 7]{CD1})
which is endowed with the symmetric monoidal model structure
in \cite[Example 7.15]{CD1}
(see loc. cit. for details).
Then we have a sequence of left Quillen symmetric monoidal functors
\[
\textup{Comp}(\mathcal{A})\longrightarrow \textup{DM}^{eff}(k) \stackrel{\Sigma^\infty}{\longrightarrow} \textup{DM}(k),
\]
where the first functor sends the unit to $L(\Spec (k))$, and the second
functor is the infinite suspension functor.

Recall the localization method in \cite[1.3.4.1, 1.3.1.15, 4.1.3.4]{HA} (see also \cite{DK}, \cite[Section 5]{Tan}); it associates to any (symmetric monoidal)
model category $\mathbb{M}$ a (symmetric monoidal) $\infty$-category
$\NNNN(\mathbb{M}^c)_\infty$. Here $\mathbb{M}^c$ is the full subcategory
spanned by cofibrant objects (this restriction is due to the technical
reason for the construction of symmetric monoidal $\infty$-categories).
We shall refer to the associated (symmetric monoidal)
$\infty$-category as the
(symmetric monoidal) $\infty$-category obtained from the model
category $\mathbb{M}$ by inverting weak equivalences.
Applying this localization,
we obtain
a symmetric monoidal functors of symmetric monoidal $\infty$-categories
\[
\Mod^\otimes_{H\mathbf{K}}\simeq \NNNN(\textup{Comp}(\mathcal{A})^c)_\infty\to \NNNN(\textup{DM}^{eff}(k)^c)_\infty \to \NNNN(\textup{DM}(k)^c)_\infty.
\]
where the first equivalence follows from \cite[8.1.2.13]{HA}.
Here $H\mathbf{K}$ denotes the Eilenberg-MacLane spectrum.
We shall write $\mathsf{DM}$ and $\mathsf{DM}^{eff}$ for $\NNNN(\textup{DM}(k)^c)_\infty$ and $\NNNN(\textup{DM}^{eff}(k)^c)_\infty$ respectively.
When we indicate that $\mathsf{DM}$ is the symmetric monoidal $\infty$-category, we denote it by $\mathsf{DM}^\otimes$.
The functor
$\Mod_{H\mathbf{K}}^\otimes\to \mathsf{DM}^\otimes$
is considered to be an $H\mathbf{K}$-linear structure.
For a proof of Theorem~\ref{tatemain},
the $H\mathbf{K}$-structure is not needed.
But $H\mathbf{K}$-linear structures are useful in other situations,
thus we will take into accout such structures in some Lemmata
and Propositions.
In \cite[Section 5]{Tan} we have constructed
another symmetric monoidal stable presentable $\infty$-category
$\mathsf{Sp}^\otimes_{\textup{Tate}}(\mathbf{HK})$
by using the recipe in \cite{CD2} and \cite{RO}.
We do not review the construction; but there is an equivalence
$\mathsf{DM}^\otimes \simeq \mathsf{Sp}^\otimes_{\textup{Tate}}(\mathbf{HK})$
(cf. \cite[Remark 6.6]{HA}).

It should be emphasized that there are several (quite different but equivalent)
constructions
of the category of mixed motives as differential graded categories
and model categories.
One can obtain $\infty$-categories from differential graded categories
and model categories. In our work,
it is important to treat ``the category of mixed motives''
as a {\it symmetric monoidal} $\infty$-category,
and therefore we choose the symmeric monoidal model category $\textup{DM}(k)$
constructed by Cisinski-D\'eglise.

\subsection{$\infty$-category of mixed Tate motives}
Let us recall the stable $\infty$-category of mixed Tate motives.
We also denote by $\mathbf{K}(1)$ its image of $\KKK(1)\in \textup{DM}^{eff}(k)$ in $\textup{DM}(k)$.
It is a cofibrant object and $\mathbf{K}(1)$ can be regard
as an object in the $\infty$-category $\mathsf{DM}$.
There exists the dual object of $\KKK(1)$ in $\DM$, which we will denote by $\KKK(-1)$.
Let $\mathsf{DTM}$ be the presentable stable subcategory generated by $\mathbf{K}(1)^{\otimes n}=\mathbf{K}(n)$
for $n\in \mathbb{Z}$,
where $\mathbf{K}(1)^{\otimes n}$ is the $n$-fold tensor product in
$\mathsf{DM}^\otimes$.
Namely, $\mathsf{DTM}$ is the smallest stable
subcategory in $\mathsf{DM}$, which admits coproducts (thus all small colimits)
and consists of $\mathbf{K}(n)$ for all $n\in \mathbb{Z}$.
The tensor product functor
$\otimes:\mathsf{DM}\times \mathsf{DM}\to \mathsf{DM}$
preserves small colimits and translations (suspensions and loops)
separately in each variable,
and thus the symmetric monoidal structure of $\mathsf{DM}$
induces a symmetric monoidal structure on $\mathsf{DTM}$.
We denote by $\mathsf{DTM}^\otimes$ the resulting
symmetric monoidal stable presentable $\infty$-category.
Note that the inclusion
$\mathsf{DTM}\hookrightarrow \mathsf{DM}$
preserves small colimits.
Let $\mathsf{DTM}_{gm}$ be the smallest stable subcategory
consisting of $\mathbf{K}(n)$ for $n\in \mathbb{Z}$.
Since $\mathbf{K}(n)$ is compact in $\mathsf{DM}$
for every $n\in \mathbb{Z}$,
every object in $\mathsf{DTM}_{gm}$
is compact in $\mathsf{DM}$.
Let $\Ind(\mathsf{DTM}_{gm})\to \mathsf{DTM}$ be
a (colimit-preserving) left Kan extension of $\mathsf{DTM}_{gm}\to \mathsf{DTM}$,
which is fully faithful by \cite[5.3.5.11]{HTT}.
Hence it identifies $\Ind(\mathsf{DTM}_{gm})$
with
$\mathsf{DTM}$.
The symmetric monoidal functor $\Mod_{H\mathbf{K}}^\otimes\to \mathsf{DM}^\otimes$ factors through $\mathsf{DTM}^\otimes \subset \mathsf{DM}^\otimes$
since $\mathsf{DTM}^\otimes \hookrightarrow \mathsf{DM}^\otimes$
preserves small colimits, and $\mathsf{DTM}$ contains
the unit of $\mathsf{DM}$.
The factorization
$\Mod_{H\mathbf{K}}^\otimes\to \mathsf{DTM}^\otimes \hookrightarrow \mathsf{DM}^\otimes$ is regarded as a map
in $\CAlg(\wCat^{\textup{L,st}})_{\Mod_{H\mathsf{K}}^\otimes/}$
which we also denote by $\mathsf{DTM}^\otimes \hookrightarrow \mathsf{DM}^\otimes$.

\begin{Lemma}
\label{tatecompactdual}
Let $\DTM_{\vee}$ be the full subcategory of $\DTM^{\otimes}$ spanned by
dualizable objects.
Let $\DTM_\circ$ be the full subcategory of $\DTM$ spanned by compact objects.
Then $\DTM_\circ=\DTM_\vee$.
\end{Lemma}

\Proof
Observe that every object in
$\DTM_{\vee}$ is compact
in $\DTM$.
To this end, it is enough to show that the unit object
of $\DTM^\otimes$ is compact (cf. \cite[2.5.1]{CD2}).
This is implied by
\cite[Theorem 2.7.10]{CD2}.
For any $n\in \ZZ$, $\KKK(n)$ belongs to $\DTM_\vee$.
Therefore $\DTM_{gm} \subset \DTM_\vee\subset \DTM_\circ$.
Notice that $\DTM_{gm}\subset \DTM_\circ$ can be
viewed as an idempotent completion (see e.g. \cite[Lemma 2.14]{BGT}).
Moreover $\DTM$ is idempotent complete by \cite[4.4.5.16]{HTT}.
It will suffice to prove that the inclusion
$\DTM_\vee\subset \DTM$ is closed under retracts.
It easily follows from the definition of dualizable objects.
\QED

\vspace{5mm}

Let $\prod_{S}\textup{DM}$ be a product 
of the category $\textup{DM}$, indexed by a small set $S$.
There is a combinatorial model structure on $\prod_{S}\textup{DM}$, called projective model
structure (cf. \cite[A. 2.8.2]{HTT}), in which weak equivalences (resp. fibrations) are termwise weak equivalences (resp. termwise fibrations) in
$\textup{DM}$. Notice that cofibrations in $\prod_{S}\textup{DM}$
are termwise cofibrations.
When $S=\NN$, $\prod_{\NN}\textup{DM}$ has  a symmetric monoidal structure
defined as follows:
Let $(M_i)_{i\in \NN}$ and $(N_j)_{j\in \NN}$ be two objects
in $\prod_{\NN}\textup{DM}$.
Then $(M_i)_{i\in \NN}\otimes(N_j)_{j\in \NN}$
is defined to be $(\oplus_{i+j=k}M_i\otimes N_j)_{k\in \NN}$.

\begin{Lemma}
\label{symmodel}
With the above symmetric monoidal structure, $\prod_{\NN}\textup{DM}$
is a symmetric monoidal model category in the sense of \cite[A 3.1.2]{HTT}.
\end{Lemma}

\Proof
We must prove that cofibrations
$\alpha:(M_i)=(M_i)_{i\in \NN}\to (M_i)=(M_i')_{i\in \NN}$
and 
$\beta:(N_i)=(N_i)_{i\in \NN}\to (N_i)=(N_i')_{i\in \NN}$
induce a cofibration
\[
\alpha\wedge \beta:(M_i)\otimes(N'_i)\coprod_{(M_i)\otimes (N_i)}(M'_i)\otimes (N_i)\to (M_i')\otimes (N_i'),
\]
and moreover if either
$\alpha$ or $\beta$ is a trivial cofibration,
then $\alpha\wedge \beta$ is also a trivial cofibration.
Unwinding the definition,
we are reduced to showing that
\[
\bigoplus_{i+j=k} \bigl( M_i\otimes N'_j\coprod_{M_i\otimes N_j}M_i'\otimes N_j\bigl) \to \bigoplus_{i+j=k} M_i'\otimes N_j'
\]
is a cofibration in $\textup{DM}$, and moreover it is a trivial cofibration
if either $\alpha$ or $\beta$ is a trivial cofibration.
This is implied by the left lifting property of (trivial) cofibrations
and the fact that $\textup{DM}$ is a symmetric monoidal model category.
\QED

Consider the symmetric monoidal
functor $\xi:\prod_{\NN}\textup{DM}\to \textup{DM}$,
which carries $(M_i)$ to $\oplus_i M_i\otimes \mathbf{K}(-i)$.
Here $\mathbf{K}(-1)$ is a cofibrant ``model''
of the dual of $\mathbf{K}(1)$, and $\mathbf{K}(-i)$ is
$i$-fold tensor product of $\mathbf{K}(-1)$
in the symmetric monoidal category $\textup{DM}$.
Since $\mathbf{K}(-i)$ is cofibrant, we see that $\xi$ is a
left Quillen adjoint functor.
By the localization, we obtain a symmetric monoidal left adjoint functor
\[
f:=\NNNN(\xi):\mathsf{DM}^\otimes_{\NN} :=\NNNN((\prod_{\NN}\textup{DM})^c)_\infty\to \NNNN(\textup{DM}^c)_\infty=\mathsf{DM}^\otimes.
\]
By the relative version of adjoint functor theorem \cite[8.3.2.6]{HA} (see also
\cite[VIII 3.2.1]{DAGn}),
$f$ has a lax symmetric monoidal right adjoint functor
which we denote by $g:\mathsf{DM}^\otimes\to \mathsf{DM}^\otimes_\NN$.
It yields $g:\CAlg(\mathsf{DM}^\otimes)\to \CAlg(\mathsf{DM}^\otimes_\NN)$.
We set $A:=g(\mathsf{1}_{\mathsf{DM}})$ in $\CAlg(\mathsf{DM}^\otimes_\NN)$,
where $\mathsf{1}_{\mathsf{DM}}$ is a unit in $\mathsf{DM}^\otimes$.
The adjoint pair
\[
f:\DM_\NN\rightleftarrows \DM:g
\]
induces
the adjoint pair
\[
f:\textup{h}(\DM_\NN)\rightleftarrows \textup{h}(\DM):g
\]
of homotopy categories.
Let $\mathsf{Hom}(N,-)$ denote the internal Hom object
given by the right adjoint of $(-)\otimes N:\DM\to \DM$.
Then $g$ is given by $M \mapsto (\mathsf{Hom}(\mathbf{K}(-i),M))_{i\in \NN}$.
Thus the underlying object $A$ in $\textup{h}(\DM)$
is $(\mathbf{K}(i))_{i\in \NN}$, that is, the $i$-th term is
$\mathbf{K}(i)$.
Moreover, by the straightforward calculation of adjunction maps,
we see that the commutative algebra structure of $A$ in the symmetric monoidal homotopy category $\textup{h}(\DM)$
is given by
\[
(\mathbf{K}(i))_{i\in \NN}\otimes(\mathbf{K}(j))_{j\in \NN}=(\oplus_{i+j=k}\mathbf{K}(i)\otimes\mathbf{K}(j))_{k\in \NN}\to (\mathbf{K}(k))_{k\in \NN}
\]
where the second map is induced by the identity maps $\mathbf{K}(i)\otimes \mathbf{K}(j)\simeq \mathbf{K}(k)\to \mathbf{K}(k)$.

Now recall from \cite{Sp1}
the notion of ``periodic'' commutative ring object
(in loc. cit. the notion of ``periodizable'' is introduced, and we
use this notion in a slightly modified form).
Let $\prod_{\ZZ}\textup{DM}$ be the product of $\textup{DM}$
indexed by $\ZZ$, which is a combinatorial model category
defined as above. By the tensor
product $(M_i)_{i\in \ZZ}\otimes(N_j)_{j\in \ZZ}=(\oplus_{i+j=k} M_i\otimes N_j)_{k\in \ZZ}$, $\prod_{\ZZ}\textup{DM}$ is
a symmetric monoidal model category in the same way that
$\prod_{\NN}\textup{DM}$ is so.
Let $\DM_\ZZ^\otimes$ be the symmetric monoidal $\infty$-category
obtained from $(\prod_{\ZZ}\textup{DM})^c$ by inverting weak equivalences.
A commutative algebra object $X$ in $\DM^\otimes_\ZZ$
is said to be periodic if
the underlying object is of the form
\[
(\ldots, \mathbf{K}(-1),\mathbf{K}(0), \mathbf{K}(1),\ldots),
\]
that is, $\mathbf{K}(i)$ sits in the $i$-th degree, and
the commutative algebra structure of $X$
in $\textup{h}(\DM_\ZZ^\otimes)$
induced by
that in $\DM_\ZZ^\otimes$
is 
determined by the identity maps
$\mathbf{K}(i)\otimes\mathbf{K}(j)\to \mathbf{K}(i+j)$.

A periodic commutative algebra object
actually exists. To construct it, we let $i:\DM_\NN^\otimes\to \DM_{\ZZ}^\otimes$
be the symmetric monoidal functor
informally given by
$(M_i)_{i\in \NN}\mapsto (\ldots,0,0,M_0,M_1,\ldots)$.
Namely, it is
determined by inserting 0 in each negative degree.
Then $P_+:=i(A)$ belongs to $\CAlg(\DM_\ZZ^\otimes)$.
According to \cite[Proposition 4.2]{Sp1} and its proof,
we have:

\begin{Proposition}[\cite{Sp1}]
There exists a morphism $P_+\to P$ in $\CAlg(\DM_\ZZ^\otimes)$
such that $P$ is periodic.
\end{Proposition}

\begin{Remark}
Let $\mathbf{K}(1)_1$ be the object of the form
$(\ldots,0,\mathbf{K}(1),0,\ldots)$ where $\mathbf{K}(1)$
sits in the 1-st degree. Let $\Sym^*_{P_+}:\Mod_{P_+}(\DM_\ZZ^\otimes)\to \CAlg(\Mod_{P_+}^\otimes(\DM_\ZZ^\otimes))$ be the
left adjoint of the forgetful functor.
Let
\[
\CAlg(\Mod_{P_+}^\otimes(\DM_\ZZ^\otimes))\rightleftarrows \CAlg(\Mod_{P_+}^\otimes(\DM_\ZZ^\otimes))[\Sym^*_{P_+}(\kappa)^{-1}]
\]
be the localization adjoint pair (cf. \cite[5.2.7.2, 5.5.4]{HTT})
which inverts
$\Sym^*_{P_+}(\kappa)$, where
$\kappa:\mathbf{K}(1)_1\otimes P_+\to P_+$ in $\Mod_{P_+}(\DM_\ZZ^\otimes)$
induced by the natural embedding $\mathbf{K}(1)_1 \to P_+$
in the 1-st degree.
The morphism $P_+\to P$ is obtained as the unit map
of this adjoint pair.
\end{Remark}

Let $\prod_{\ZZ}\textup{Comp}(\mathcal{A})$ be the product of the category
$\textup{Comp}(\mathcal{A})$, that is endowed with the projective
model structure.
As in Lemma~\ref{symmodel}, we see that
$\prod_{\ZZ}\textup{Comp}(\mathcal{A})$ is a symmetric monoidal
model category, whose tensor product is
given by $(A_i)_{i\in \ZZ}\otimes(B_j)_{j\in \ZZ}=(\oplus_{i+j=k}A_i\otimes B_j)_{k\in \ZZ}$.
Then the natural left Quillen adjoint symmetric monoidal functor
$\textup{Comp}(\mathcal{A})\to \textup{DM}$ naturally extends to
a left Quillen adjoint symmetric monoidal functor
$l:\prod_{\ZZ}\textup{Comp}(\mathcal{A})\to \prod_{\ZZ}\textup{DM}$.
It gives rise to the
symmetric monoidal left adjoint functor of $\infty$-categories
\[
l:\Mod_{H\mathbf{K},\ZZ}^\otimes:=\NNNN(\prod_{\ZZ}\textup{Comp}(\mathcal{A})^c)_\infty^\otimes\to \DM_\ZZ^\otimes.
\]
According to the relative version of adjoint functor theorem \cite[8.3.2.6]{HA} (see also
\cite[VIII 3.2.1]{DAGn}),
$l$ has a lax symmetric monoidal right adjoint functor $r$.
Let $Q:=r(P)\in \CAlg(\Mod_{H\mathbf{K},\ZZ}^\otimes)$.
Let $\textup{DM}\to \prod_{\ZZ}\textup{DM}_\ZZ$ be the
left Quillen symmetric monoidal functor which carries $M$ to $(M_i)$ where $M_0=M$ and
$M_i=0$ if $i\neq 0$.
Thus we have a symmetric monoidal functor $\DM\to \DM_\ZZ$,
and again by the relative version of adjoint functor theorem
we obtain a lax symmetric monoidal functor $s:\DM_\ZZ\to \DM$
as the right adjoint.
Therefore there exists a diagram of symmetric monoidal $\infty$-categories:
\[
\xymatrix{
  & \Mod_{l(Q)}(\DM_{\ZZ}^\otimes) \ar[d]^{u} & \\
\Mod_{Q}(\Mod^\otimes_{H\mathbf{K},\ZZ}) \ar[ru]^{\tilde{l}} \ar[r]^{u\circ \tilde{l}} \ar@<5pt>[d]^{b} & \Mod_{P}(\DM_\ZZ^\otimes) \ar[d]^{t} \ar[rd]^{s\circ t}& \\
\Mod_{H\mathbf{K},\ZZ} \ar@<5pt>[r]^{l} \ar@<5pt>[u]^{a} &  \DM_{\ZZ} \ar@<5pt>[l]^{r} \ar[r]^{s} & \DM
}
\]
such that
\begin{itemize}
\item $\tilde{l}$ is a symmetric monoidal functor induced by $l$,

\item $u$ is the symmetric monoidal base change functor induced by the counit map $l(Q)=l(r(P))\to P$,

\item $t$ is the forgetful monoidal functor which is a lax symmetric monoidal functor,

\item $a$ is the base change functor, and $b$ is the forgetful functor.

\end{itemize}
Let $z:=s\circ t\circ u\circ \tilde{l}$.
We recall the theorem by Spitzweck \cite[Theorem 4.3]{Sp1}
(see also its proof):

\begin{Theorem}[\cite{Sp1}]
\label{Spitzweckrep}
The composite $z:\Mod_{Q}(\Mod^\otimes_{H\mathbf{K},\ZZ})\to \DM$ gives
an equivalence 
\[
\Mod_{Q}(\Mod^\otimes_{H\mathbf{K},\ZZ})\simeq \DTM
\]
as symmetric monoidal $\infty$-categories.
\end{Theorem}

\begin{Remark}
This result is extended to a more general situation by a different method
\cite{PM}.
\end{Remark}

Furthermore, we can see that $z$ gives an equivalence of them
as {\it $H\mathbf{K}$-linear} symmetric monoidal $\infty$-categories.
To see this, it is enough to show that
$z$ is promoted to a
$H\mathbf{K}$-linear symmetric monoidal functor.
To treat problems of this type,
the following Lemma is useful.

\begin{Lemma}
\label{linearity}
Let $\mathcal{C}^\otimes$ be in $\CAlg(\wCat^{\textup{L,st}})$.
We denote by $\mathcal{C}$ the underlying $\infty$-category.
Suppose that a unit $\mathbf{1}$ of $\mathcal{C}^\otimes$ is
compact in $\mathcal{C}$.
Let $\mathcal{C}_{\mathbf{1}}\subset \mathcal{C}$ be the smallest stable subcategory
which admits small colimits and contains $\mathbf{1}$.
The $\infty$-category $\mathcal{C}_{\mathbf{1}}$
admits a symmetric monoidal structure induced by that of $\mathcal{C}^\otimes$.
Then there exist $A$ in $\CAlg$
and an equivalence $\Mod_A^\otimes\simeq \mathcal{C}^\otimes$
of symmetric monoidal $\infty$-categories.
Moreover, if $R$ is a commutative ring spectrum
and $p:\Mod_R^\otimes\to \mathcal{C}^\otimes$ is
a symmetric monoidal colimit-preserving functor, then
$p$ factors through $\mathcal{C}_{\mathbf{1}}^\otimes\subset \mathcal{C}^\otimes$ and there exists a morphism $R\to A$ in $\CAlg$, up to the
contractible space of choice, which induces
$\Mod_{R}^\otimes\to\mathcal{C}_{\mathbf{1}}^\otimes\simeq \Mod_A^\otimes$
(as the base change).
\end{Lemma}

\Proof
The first assertion follows from
\cite[8.1.2.7]{HA}; the characterization of symmetric monoidal
stable $\infty$-categories of module spectra.
Since $p$ preserves small colimits,
$p$ factors through $\mathcal{C}_{\mathbf{1}}^\otimes\subset \mathcal{C}^\otimes$.
The last assertion is implied by \cite[6.3.5.18]{HA}.
\QED

\begin{Remark}
Under the assumption of Lemma~\ref{linearity},
$A$ is considered to be the ``endomorphism algebra'' of the unit, and
we can say that
giving a $R$-linear structure, that is, a symmetric monoidal colimit-preserving functor $\Mod_R^\otimes \to \mathcal{C}^\otimes$
is equivalent to giving a morphism $R\to A$ in $\CAlg$.
\end{Remark}

Return to the case of $H\mathbf{K}$-linear symmetric monoidal $\infty$-category
$\DTM^\otimes$. The endomorphism algebra of the unit of $\DTM^\otimes$
is $H\mathbf{K}$ (i.e. $\mathbf{K}$),
and its $H\mathbf{K}$-linear structure is determined by
the identity $H\mathbf{K}\to H\mathbf{K}$.
Thus, to promote $z$ to a $H\mathbf{K}$-linear symmetric monoidal functor,
it is enough to show that $f\circ a\circ q:\Mod_{H\mathbf{K}}^\otimes\to \DTM^\otimes$
induces the identity morphism of endomorphism algebras of units
$H\mathbf{K}\to H\mathbf{K}$,
where $q$ is the inclusion $\Mod_{H\mathbf{K}}^\otimes\to \Mod_{H\mathbf{K},\ZZ}^\otimes$ into the degree zero part.
This claim is clear from our construction.

\subsection{Realization functor and augmentation}
Let $E$ be a mixed Weil theory with $\mathbf{K}$-coefficients
(cf. \cite[Definition 2.1]{CD2}). A
 mixed Weil theory is a presheaf of commutative dg
$\mathbf{K}$-algebras on the category of smooth affine schemes over $k$,
which satisfies Nisnevich descent property,
$\mathbb{A}^1$-homotopy, K\"unneth formula and axioms of dimensions, etc
(for the precise definition see \cite[2.1.4]{CD2}).
For example, algebraic de Rham cohomology determines a mixed Weil theory
with $\mathbf{K}=k$;
to any smooth affine scheme $X$ we associates
a commutative dg $k$-algebra $\Gamma(X,\Omega_{X/k}^\ast)$
where $\Omega_{X/k}^\ast$ is the algebraic de Rham complex arising from
the exterior $\mathcal{O}_X$-algebra generated by $\Omega_{X/k}^1$.
Another example is $l$-adic \'etale cohomology with
$\mathbf{K}=\QQ_l$ (see \cite[Section 3]{CD2}).
To a mixed Weil theory $E$ we can associate
\[
\mathsf{R}_E:\DM^{\otimes}\to \Mod_{H\mathbf{K}}^\otimes
\]
a morphism in
 $\CAlg(\wCat^{\textup{L,st}})_{\Mod_{H\mathbf{K}}^\otimes/}$ which we call the homological realization functor with respect to $E$ (see \cite[Section 5.1, 5.2]{Tan}, \cite[2.6]{CD2}).
From now on we usually omit the subscript $E$.
Then according to
\cite[2.7.14]{CD2}
when $E$ is the mixed Weil theory
associated to algebraic de Rham cohomology, for any smooth affine scheme $X$
the image $\mathsf{R}(h(X))$ in $\Mod_{H\mathbf{K}}$
is equivalent to the dual complex of
the derived global section
$\mathbf{R}\Gamma(X,\Omega_{X/k}^\ast)$
where by \cite[5.10]{Tan}
we identify $\Mod_{H\mathbf{K}}$ with the $\infty$-category
of unbounded complexes of $\mathbf{K}$-vector spaces.
We denote by $\mathsf{R}_T$ the composition
\[
\DTM^\otimes\hookrightarrow \DM^\otimes\to \Mod_{H\mathbf{K}}^\otimes
\]
which we call the homological realization of Tate motives
(with respect to $E$).
By the restrictions,
it gives rise to the morphism $\DTM_{\vee}^\otimes \to \PMod_{H\mathbf{K}}^\otimes$ in $\CAlg(\sCat)_{\PMod^\otimes_{H\mathbf{K}}/}$ which we denote also by $\mathsf{R}_T$.

Combined with Theorem~\ref{Spitzweckrep}
we have the sequence of symmetric monoidal colimit-preserving functors
\[
\Mod^\otimes_{H\mathbf{K},\ZZ}\stackrel{a}{\longrightarrow} \Mod^\otimes_{Q}(\Mod_{H\mathbf{K},\ZZ})\simeq \DTM^\otimes \stackrel{\mathsf{R}_T}{\longrightarrow} \Mod^\otimes_{H\mathbf{K}}.
\]
By the relative version of adjoint functor theorem,
the composition admits a lax symmetric monoidal right adjoint functor
$\xi$.
In particular, if we set $R=\xi(\mathsf{1}_{H\mathbf{K}})$
with $\mathsf{1}_{H\mathbf{K}}$ the unit of $\Mod_{H\mathbf{K}}^\otimes$,
then $R$ belongs to $\CAlg(\Mod_{H\mathbf{K},\ZZ}^\otimes)$.
By the functoriality and the construction of $Q$,
we have the natural morphism $Q\to R$
in $\CAlg(\Mod_{H\mathbf{K},\ZZ}^\otimes)$.
There is a commutative diagram (up to homotopy) of symmetric monoidal $\infty$-categories
\[
\xymatrix{
\Mod^\otimes_Q(\Mod_{H\mathbf{K},\ZZ}^\otimes) \ar[r]^{\sim}_{z} \ar[d] & \DTM^\otimes \ar[r]_{\mathsf{R}_T} \ar[d] & \Mod^\otimes_{H\mathbf{K}} \\
\Mod^\otimes_{R}(\Mod_{H\mathbf{K},\ZZ}^\otimes) \ar[r]_{\tilde{z}} & \Mod^\otimes_{f(R)}(\DTM^\otimes) \ar[r]_{\tilde{\mathsf{R}}_T} & \Mod^\otimes_{\mathsf{R}_T(z(R))}(\Mod_{H\mathbf{K}}^\otimes) \ar[u]
}
\]
where $\tilde{z}$ and $\tilde{\mathsf{R}}_T$ are induced by
$z$ and $\mathsf{R}_T$ respectively,
the left and central vertical arrows are base change functors,
and the right vertical arrow is the counit map $\mathsf{R}_T(z(R))\to H\mathbf{K}$
in $\CAlg(\Mod_{H\mathbf{K}}^\otimes)$.
Note that all functors in the diagram are
$H\mathbf{K}$-linear symmetric monoidal functors.
The commutativity of the right square follows from
the observation that the counit map $\mathsf{R}_T(z(R))\to H\mathbf{K}$
is an augmentation of the structure map $H\mathbf{K}\to \mathsf{R}_T(z(R))$.

\begin{Lemma}
\label{ppp}
The composite $h:\mathcal{C}^\otimes:=\Mod_{R}^\otimes(\Mod_{H\mathbf{K},\ZZ}^\otimes)\to \mathcal{D}^\otimes:=\Mod^\otimes_{H\mathbf{K}}$ in the above diagram gives an equivalence of $H\mathbf{K}$-linear symmetric monoidal $\infty$-categories.
\end{Lemma}

\Proof
It will suffice to show that the underlying functor is a categorical
equivalence.

The symmetric monoidal functor $h$ is $H\mathbf{K}$-linear.
Thus $h$ is essentially surjective.

Next we will show that $h$ is fully faithful.
Let $\mathbf{K}_n:=(\ldots0,\mathbf{K},0\ldots)$ be the
object in $\Mod_{H\mathbf{K},\ZZ}$ such that $\mathbf{K}$ sits in the
$n$-th degree. Let $R(n)$ be the image of
$\mathbf{K}_n$ by the base change functor
$\Mod_{H\mathbf{K},\ZZ}\to \Mod_R(\Mod_{H\mathbf{K},\ZZ}^\otimes)$.
(For any $n\in \ZZ$, $h(R(n))\simeq H\mathbf{K}$.)
It is enough to prove that
\[
\Map_{\mathcal{C}}(R(i),R(j))\to \Map_{\mathcal{D}}(h(R(i)),h(R(j)))
\]
is an equivalence in $\mathcal{S}$.
Indeed, $\mathcal{C}$ is generated by the sets $\{R(i)\}_{i\in \ZZ}$
under finite (co)limits, translations, and filtered colimits.
Since $R(i)$ and $h(R(i))$ is compact for each $i\in\ZZ$
and $h$ is colimit-preserving, we are reduced to showing
that the above map is an equivalence in $\mathcal{S}$.
(Assuming it to hold, note first that
$\Map_{\mathcal{C}}(R(i),N)\to\Map_{\mathcal{D}}(h(R(i)),h(N))$
is an equivalence in $\SSS$
for $N$ being in the smallest 
stable subcategory $\mathcal{C'}$ generated by $\{R(i)\}_{i\in \ZZ}$.
Then since $R(i)$ and $h(R(i))$ are compact, $\Ind(\mathcal{C}')\simeq \mathcal{C}$, and $h$ preserves small colimits, thus
for any $N\in \mathcal{C}$, $\Map_{\mathcal{C}}(R(i),N)\to\Map_{\mathcal{D}}(h(R(i)),h(N))$ is an equivalence. Since $\mathcal{C}$ is generated by $\{R(i)\}_{i\in \ZZ}$ under finite colimits, translations and filtered colimits,
we conclude that for any $M, N\in \mathcal{C}$, $\Map_{\mathcal{C}}(M,N)\to\Map_{\mathcal{D}}(h(M),h(N))$ is an equivalence.)
Note that $\Map_{\mathcal{C}}(R(i),R(j))\simeq \Map_{\mathcal{C}}(R(i-j),R)$,
and therefore we may and will assume that $j=0$.
Then by using adjunctions
we can identify
$\Map_{\mathcal{C}}(R(i),R)\to \Map_{\mathcal{D}}(h(R(i)),h(R))$
with the composition
\begin{eqnarray*}
\Map_{\mathcal{C}}(R(i),R) &\stackrel{\sim}{\to}& \Map_{\Mod_Q(\Mod_{H\mathbf{K},\ZZ}^\otimes)}(Q(i),R)  \\
&\stackrel{\sim}{\to}& \Map_{\Mod_{H\mathbf{K}}}(\mathsf{R}_T(z(Q(i))),H\mathbf{K}) \\
&\stackrel{\sim}{\to} & \Map_{\Mod_{H\mathbf{K}}}(H\mathbf{K},H\mathbf{K}).
\end{eqnarray*}
This proves our Lemma.
\QED

\begin{Proposition}
\label{kimon?}
There exists a $H\mathbf{K}$-linear symmetric monoidal
equivalence
\[
\Mod_{H\mathbf{K},\ZZ}^\otimes \to \Mod_{\mathsf{B}\mathbb{G}_m}^\otimes.
\]
\end{Proposition}

\Proof
We will construct a symmetric monoidal functor
$\Mod_{H\mathbf{K},\ZZ}^\otimes \to \Mod_{\mathsf{B}\mathbb{G}_m}^\otimes$,
which preserves colimits.

For this purpose, we will construct $\Mod_{\mathsf{B}\mathbb{G}_m}^\otimes$
in a somewhat explicit way.
Regard the group scheme $\mathbb{G}_m$ over $\mathbf{K}$
as the simplicial scheme, denoted by $G_\bullet$ such that
$G_i$ is the $i$-fold product $\mathbb{G}_m^{\times i}$.
This corresponds to the cosimplicial $\mathbf{K}$-algebra
$\Gamma(G)^{\bullet}$ such that $\Gamma(G)^{i}\simeq \mathbf{K}[t_1^{\pm},\ldots,t^{\pm}_i]$.
The cosimplicial $\mathbf{K}$-algebra $\Gamma(G)^\bullet$
naturally induces the cosimplicial diagram $\rho:\NNNN(\Delta)\to \wCat$
such that $\rho([i])=\NNNN(\textup{Comp}(\Gamma(G)^i)^c)$.
Here $\textup{Comp}(\Gamma(G)^i)$ denotes the category of chain complexes of
$\Gamma(G)^i$-modules which is endowed with the projective model structure,
and $\textup{Comp}(\Gamma(G)^i)^c$ is its full subcategoy of cofibrant objects.
Each category $\textup{Comp}(\Gamma(G)^i)^c$ has the (natural)
symmetric monoidal structure, and thus
$\rho$ is promoted to $\rho:\NNNN(\Delta)\to \CAlg(\wCat)$,
where $\CAlg(\wCat)$ is the $\infty$-category of
symmetric monoidal $\infty$-categories (i.e., commutative algebra
objects in the Cartesian symmetric monoidal $\infty$-category $\wCat$).
The symmetric monoidal category $\textup{Comp}(\Gamma(G)^i)^c$
admits the subset of edges of weak equivalences.
Inverting weak equivalences in
$\textup{Comp}(\Gamma(G)^i)^c$,
we have $\rho':\NNNN(\Delta)\to \CAlg(\wCat)$
and the natural transformation $\rho\to \rho'$.
such that $\rho'([i])$ is a symmetric monoidal
$\infty$-category obtained from $\textup{Comp}(\Gamma(G)^i)^c$
by inverting weak equivalences.

By the explicit unstraightening functor \cite[3.2.5.2]{HTT},
the maps $\rho,\rho':\NNNN(\Delta)\rightrightarrows \CAlg(\wCat)$
gives rise to
coCartesian fibrations
$\mathcal{C}_{pre}^\otimes\to \NNNN(\FIN)\times \NNNN(\Delta)$
and $\mathcal{C}^\otimes\to \NNNN(\FIN)\times \NNNN(\Delta)$.
The natural transformation $\rho\to \rho'$
induces a map of coCartesian fibrations
\[
\xymatrix{
\mathcal{C}_{pre}^\otimes \ar[rr]^{\sigma} \ar[rd] &  &\mathcal{C}^\otimes \ar[dl] \\
 & \NNNN(\FIN)\times \NNNN(\Delta) & 
}
\]
which preserves coCartesian edges.
Note for each $[i]\in \Delta$, the fiber $\rho^{-1}([i])\to \NNNN(\FIN)\times\{[i]\}\cong \NNNN(\FIN)$ is the symmetric monoidal
$\infty$-category associated to the diagram of $\textup{Comp}(\Gamma(G)^i)^c$'s.
The fiber $(\rho')^{-1}([i])\to \NNNN(\FIN)$ is the symmetric monoidal
$\infty$-category obtained from $\textup{Comp}(\Gamma(G)^i)^c$
by inverting weak equivalences.

Next we define a map of simplicial sets $\overline{\SEC}(\mathcal{C}_{pre}^\otimes)\to \NNNN(\FIN)$ as follows. For any $a:T\to \NNNN(\FIN)$,
giving a map $T\to \overline{\SEC}(\mathcal{C}_{pre}^\otimes)$ over $\NNNN(\FIN)$
amounts to giving $\phi:T\times \NNNN(\Delta)\to \mathcal{C}_{pre}^\otimes$
which commutes with $a\times\textup{Id}:T\times\NNNN(\Delta)\to \NNNN(\FIN)\times\NNNN(\Delta)$ and $\mathcal{C}_{pre}^\otimes\to \NNNN(\FIN)\times \NNNN(\Delta)$.
Let $\SEC(\mathcal{C}_{pre}^\otimes)$ be the largest subcomplex of
$\overline{\SEC}(\mathcal{C}_{pre}^\otimes)$, which consists of
the following vertexes:
a vertex $v\in \overline{\SEC}(\mathcal{C}_{pre}^\otimes)$ lying over
$\langle i \rangle$ belongs to $\SEC(\mathcal{C}_{pre}^\otimes)$
exactly when $v:\{\langle i\rangle\}\times \NNNN(\Delta)\to \mathcal{C}_{pre}^\otimes$ carries all edges in $\{\langle i\rangle\}\times \NNNN(\Delta)$
to coCartesian edges in $\mathcal{C}_{pre}^\otimes$.
We define $\overline{\SEC}(\mathcal{C}^\otimes)\to \NNNN(\FIN)$
and $\SEC(\mathcal{C}^\otimes)\to \NNNN(\FIN)$
in a similar way. According to \cite[3.1.2.1 (1)]{HTT},
 we see that $\overline{\SEC}(\mathcal{C}_{pre}^\otimes)\to \NNNN(\FIN)$
and $\overline{\SEC}(\mathcal{C}^\otimes)\to \NNNN(\FIN)$
are coCartesian fibrations (notice that
$\overline{\SEC}(\mathcal{C}_{pre}^\otimes)=\NNNN(\FIN)\times_{\Fun(\NNNN(\Delta),\NNNN(\FIN)\times \NNNN(\Delta))}\Fun(\NNNN(\Delta),\mathcal{C}_{pre}^\otimes)$) where $\NNNN(\FIN)\to \Fun(\NNNN(\Delta),\NNNN(\FIN)\times \NNNN(\Delta))$ is induced by the identity
$\NNNN(\FIN)\times\NNNN(\Delta)\to \NNNN(\FIN)\times \NNNN(\Delta)$).
Moreover, by \cite[3.1.2.1 (2)]{HTT}
we deduce that
$\SEC(\mathcal{C}_{pre}^\otimes)\to \NNNN(\FIN)$
and $\SEC(\mathcal{C}^\otimes)\to \NNNN(\FIN)$
are coCartesian fibrations.
By construction, furthermore $\SEC(\mathcal{C}_{pre}^\otimes)\to \NNNN(\FIN)$
is a symmetric monoidal $\infty$-category.
Since the procedure of inverting weak equivalences
commutes with finite products \cite[4.1.3.2]{HA},
we see that $\SEC(\mathcal{C}^\otimes)\to \NNNN(\FIN)$
is also a symmetric monoidal $\infty$-category.
We will abuse notation and denote by $\SEC(\mathcal{C}_{pre}^\otimes)$
and $\SEC(\mathcal{C}^\otimes)$ the underlying $\infty$-categories.
Note that $\sigma$ (which preserves coCartesian edges)
induces a symmetric monoidal
 functor $\SEC(\mathcal{C}_{pre}^\otimes)\to \SEC(\mathcal{C}^\otimes)$.

Observe that the symmetric monoidal $\infty$-category
$\SEC(\mathcal{C}^\otimes)\to \NNNN(\FIN)$
is equivalent to the symmetric monoidal
$\infty$-category $\Mod_{\mathsf{B}\mathbb{G}_m}^\otimes$.
By \cite[3.3.3.2]{HTT} and \cite[3.2.2.4]{HA}, the symmetric monoidal
$\infty$-category $\SEC(\mathcal{C}^\otimes)$
is a limit of the diagram $\rho':\NNNN(\Delta)\to \CAlg(\wCat)$.
Note that by \cite[8.1.2.13]{HA} $\rho'([i])$ is equivalent to
$\Mod_{\Gamma(G)^i}^\otimes$.
Beside, the functor $\Theta:\CAlg \to \CAlg(\wCat^{\textup{L,st}})$
which carries $A$ to $\Mod_A^\otimes$ (see Section 3.1)
is fully faithful \cite[6.3.5.18]{HA}.
For a symmetric monoidal functor $\phi:\Mod_A^\otimes\to \Mod_B^\otimes$
in $\CAlg(\wCat^{\textup{L,st}})$, one can recover $f:A\to B$ with
$\Theta(f)\simeq \phi$
as the induced morphism from the endomorphism spectrum of a unit of $\Mod_A^\otimes$
to that of the unit in $\Mod_B^\otimes$.
Therefore from the construction of $\rho'$ (and $\rho$)
and the definition of $\Mod_{\mathsf{B}\mathbb{G}_m}^\otimes$,
we conclude that 
$\SEC(\mathcal{C}^\otimes)\to \NNNN(\FIN)$
is equivalent to $\Mod_{\mathsf{B}\mathbb{G}_m}^\otimes$.

Therefore, to construct $\Mod_{H\mathbf{K},\ZZ}^\otimes\to \Mod_{\mathsf{B}\mathbb{G}_m}^\otimes$ it suffices to construct a symmetric monoidal functor
from $\prod_{\ZZ}\textup{Comp}(\mathcal{A})^c$ to $\SEC(\mathcal{C}_{pre}^\otimes)$
which carries weak equivalences in $\prod_{\ZZ}\textup{Comp}(\mathcal{A})^c$
to edges in $\SEC(\mathcal{C}_{pre}^\otimes)$
whose images in $\SEC(\mathcal{C}^\otimes)$
are equivalences (note the universality of $\Mod_{H\mathbf{K},\ZZ}^\otimes$
\cite[4.1.3.4]{HA}).
Let $\mathbf{K}_n$ in $\prod_{\ZZ}\textup{Comp}(\mathcal{A})^c$
be the $\mathbf{K}$ which sits in the $n$-th degree with respect  to $\prod_{\ZZ}$. To $\mathbf{K}_n$ we attach the weight $n$ representation of
$\mathbb{G}_m$ on $\mathbf{K}$.
The weight $n$ representation gives rise to
an object of $\SEC(\mathcal{C}_{pre}^\otimes)$
in the obvious way, which we denote by $\mathbf{K}_n'$.
For $(M_i)_{i\in\ZZ}\in \prod_{\ZZ}\textup{Comp}(\mathcal{A})^c$,
we attach $\oplus_{i\in \ZZ}M_i\otimes \mathbf{K}_i'$.
Here we consider $M_i$ to be an object in $\SEC(\mathcal{C}_{pre}^\otimes)$,
that is
the complex endowed with the trivial action of $\mathbb{G}_m$.
This naturally induces a symmetric monoidal functor having the desired
property.
To prove that the induced functor $\Mod_{H\mathbf{K},\ZZ}\to \Mod_{\mathsf{B}\mathbb{G}_m}$ preserves small colimits, it is enough to show that
the composite $\Mod_{H\mathbf{K},\ZZ}\to \Mod_{\mathsf{B}\mathbb{G}_m}\to \Mod_{H\mathbf{K}}$, where the second functor is forgetful, preserves small colimits
since the forgetful functor is conservative and preserves small colimits
(an exact functor $p:\mathcal{K}\to \mathcal{L}$ between stable $\infty$-categories is said
to be conservative if for any $K\in \mathcal{K}$,
$p(K)\simeq 0$ implies that $K\simeq 0$).
The composite carries $(M_i)_{i\in \ZZ}$ to $\oplus_{i\in \ZZ}M_i$
and thus we conclude that the composite preserves small colimits.
To prove that $\Mod_{H\mathbf{K},\ZZ}^\otimes\to \Mod_{\mathsf{B}\mathbb{G}_m}^\otimes$
is promoted to a $H\mathbf{K}$-linear symmetric monoidal functor,
according to Lemma~\ref{linearity} (see also the discussion at the end of 6.3),
it suffices to observe that
$\Mod_{H\mathbf{K},\ZZ}^\otimes\to \Mod_{\mathsf{B}\mathbb{G}_m}^\otimes$ induces the identity morphism $H\mathbf{K}\to H\mathbf{K}$
of endomorphism algebras of units.
To see this, we are reduced to showing that
the composite
$\Mod_{H\mathbf{K},\ZZ}^\otimes\to \Mod_{\mathsf{B}\mathbb{G}_m}^\otimes\to \Mod_{H\mathbf{K}}^\otimes$, where the second functor is the forgetful functor,
induces the identity morphism
$H\mathbf{K}\to H\mathbf{K}$
of endomorphism algebras of units.
This is clear.

We have constructed
a symmetric monoidal colimit-preserving functor $\Mod^{\otimes}_{H\mathbf{K},\ZZ}\to \Mod^\otimes_{\mathsf{B}\mathbb{G}_m}$ with the (lax symmetric monoidal) right adjoint functor (the existence is assured by the relative version of adjoint functor theorem). To see that $\Mod^{\otimes}_{H\mathbf{K},\ZZ}\to \Mod^\otimes_{\mathsf{B}\mathbb{G}_m}$ is an equivalence of symmetric monoidal $\infty$-categories,
it is enough to show that it induces a categorical equivalence
$\Mod_{H\mathbf{K},\ZZ}\to \Mod_{\mathsf{B}\mathbb{G}_m}$ of underlying $\infty$-categories.
Moreover, by \cite[5.8]{Tan}, it suffices to check that
it induces an equivalence $\textup{h}(\Mod_{H\mathbf{K},\ZZ})\to \textup{h}(\Mod_{\mathsf{B}\mathbb{G}_m})$ of their homotopy categories.
The desired equivalence now follows from \cite[Section 8, Theorem 8.5]{Spi}
(see also the strictification theorem \cite[18.7]{HS}).
\QED


Let $A$ be an object in $\CAlg(\Mod_{\mathsf{B}\mathbb{G}_m}^\otimes)$.
Let $\overline{A}$ denote the image of $A$ in
$\CAlg(\Mod_{H\mathbf{K}}^\otimes)$ (via the pullback of $\Spec H\mathbf{K}\to \mathsf{B}\mathbb{G}_m$).
With the notation in the proof of Proposition~\ref{kimon?}, there is the natural augmented simplicial diagram $G_\bullet\to \mathsf{B}\mathbb{G}_m$.
This induces the natural functor $\CAlg(\Mod_{\mathsf{B}\mathbb{G}_m})\to \lim_{[i]\in\Delta} \CAlg(\Mod_{H\Gamma(G)^i})$.
We write $A^\bullet$ for the image of $A$ in $\lim_{[i]\in \Delta} \CAlg(\Mod_{H\Gamma(G)^i})$.
It gives rise to the quotient stack $[\Spec \overline{A}/\mathbb{G}_m]$ (see Example~\ref{quotient}).
The construction of the quotient derived stack is as follows:
The cosimplicial diagram $\{\Gamma(G)^{i}\}_{[i]\in \Delta}$ of
ordinary commutative $\mathbf{K}$-algebras
has a natural map from the constant simplicial diagram $\{\mathbf{K}\}$.
Both cosimplicial diagrams
naturally induce $c_{\mathbf{K}}, c_G:\NNNN(\Delta)\to \widehat{\textup{Cat}}_{\infty}$
such that $c_{\mathbf{K}}$ is the constant diagram of $\CAlg_{H\mathbf{K}}$,
and $c_{G}([i])=\CAlg_{H\Gamma(G)^{i}}$ and $[i]\to [j]$ maps to
$\CAlg_{H\Gamma(G)^{i}}\to \CAlg_{H\Gamma(G)^{j}}; R\mapsto H\Gamma(G)^{j}\otimes_{H\Gamma(G)^{i}}R$.
By \cite[3.2.0.1, 4.2.4.4]{HTT} the cosimplicial diagrams $c_{\mathbf{K}}$ and $c_{G}$
give rise to coCartesian fibrations $\textup{pr}_2:\CAlg_{H\mathbf{K}}\times\NNNN(\Delta)\to \NNNN(\Delta)$ and $\overline{\CAlg}_{G}\to \NNNN(\Delta)$
respectively.
The morphism $c_{\mathbf{K}}\to c_{G}$ induced by $\{\mathbf{K}\}\to \{\Gamma(G)^{i}\}_{[i]\in \Delta}$ gives rise to a morphism of coCartesian fibrations
$\alpha:\CAlg_{H\mathbf{K}}\times \NNNN(\Delta)\to \overline{\CAlg}_{G}$
over $\NNNN(\Delta)$ that preserves coCartesian edges. By \cite[8.3.2.7]{HA}, there is a right adjoint $\beta$
of
$\alpha$ relative to $\NNNN(\Delta)$.
Let $s:\NNNN(\Delta)\to \overline{\CAlg}_{G}$ be a section that
corresponding to $A^\bullet$ (cf. \cite[3.3.3.2]{HTT}).
Then the composite
\[
\xi:\NNNN(\Delta)\stackrel{s}{\to} \overline{\CAlg}_{G}\stackrel{\beta}{\to} \CAlg_{H\mathbf{K}}\times \NNNN(\Delta)\stackrel{\textup{pr}_1}{\to} \CAlg_{H\mathbf{K}}
\]
gives rise to a simplicial diagram $\xi^{op}:\NNNN(\Delta)^{op}\to \Aff_{H\mathbf{K}}$.
We define $[\Spec \overline{A}/\mathbb{G}_m]$
to be a colimit (geometric realization)
of the simplicial diagram $\xi^{op}$ in $\Fun(\CAlg_{H\mathbf{K}},\mathcal{S})$.
If $s_1:\NNNN(\Delta)\to \overline{\CAlg}_{G}$ is the section corresponds to
the initial object of $\lim_{[i]\in \Delta} \CAlg(\Mod_{H\Gamma(G)^i})$,
then the composite $\textup{pr}_1\circ \beta \circ s_1:\NNNN(\Delta)\to \CAlg_{H\mathbf{K}}$ is equivalent to the cosimplicial diagram
$\{H\Gamma(G)^i\}$. Thus we have a natural morphism $\pi:[\Spec \overline{A}/\mathbb{G}_m]\to B\mathbb{G}_m$. It is easy to see that this morphism
makes it a quotient stack.

\begin{Proposition}
\label{yoyu}
There exists a natural equivalence
\[
\Mod_A^\otimes(\Mod_{\mathsf{B}\mathbb{G}_m}^\otimes)\simeq \Mod_{[\Spec \overline{A}/\mathbb{G}_m]}^\otimes.
\]
\end{Proposition}

\Proof
We first construct
a symmetric monoidal colimit-preserving functor
\[
\Mod_A^\otimes(\Mod_{\mathsf{B}\mathbb{G}_m}^\otimes)\longrightarrow \Mod^\otimes_{[\Spec \overline{A}/\mathbb{G}_m]}.
\]
Let $\pi^*:\Mod_{\mathsf{B}\mathbb{G}_m}^\otimes\to  \Mod^\otimes_{[\Spec \overline{A}/\mathbb{G}_m]}$ be the symmetric monoidal functor
induced by the natural morphism
$\pi:[\Spec \overline{A}/\mathbb{G}_m]\to \mathsf{B}\mathbb{G}_m$.
By the relative version of adjoint functor theorem,
there is a lax symmetric monoidal right adjoint fucntor
$\pi_*:\Mod_{[\Spec \overline{A}/\mathbb{G}_m]} \to \Mod_{\mathsf{B}\mathbb{G}_m}$.
If $\mathsf{1}_{[\Spec \overline{A}/\mathbb{G}_m]}$ is a unit of
$\Mod^\otimes_{[\Spec \overline{A}/\mathbb{G}_m]}$, by the definition of
$[\Spec \overline{A}/\mathbb{G}_m]$ and the base-change formula,
$\pi_*(\mathsf{1}_{[\Spec \overline{A}/\mathbb{G}_m]})$ is equivalent to $A$ in $\CAlg(\Mod_{\mathsf{B}\mathbb{G}_m}^\otimes)$.
Thus
we have the composition of symmetric monoidal colimit-preserving functors
\[
h:\Mod_{A}^\otimes(\Mod_{\mathsf{B}\mathbb{G}_m}^\otimes)\to \Mod_{\pi^*(A)}(\Mod^\otimes_{[\Spec \overline{A}/\mathbb{G}_m]})\to \Mod^\otimes_{[\Spec \overline{A}/\mathbb{G}_m]}
\]
where the second functor is induced by the counit map
$\pi^*(A)\simeq \pi^*(\pi_*(\mathsf{1}_{[\Spec \overline{A}/\mathbb{G}_m]}))\to \mathsf{1}_{[\Spec \overline{A}/\mathbb{G}_m]}$.
Note that the composite is naturally a $H\mathbf{K}$-linear symmetric monoidal functor.

Next we will show that $h$ gives an equivalence of symmetric monoidal
$\infty$-categories.
It will suffice to prove that
the underlying functor of $\infty$-categories
is a categorical equivalence.
We first show that $h$ is fully faithful.
Let $\mathsf{1}_{\mathsf{B}\mathbb{G}_m}(i)\in \Mod_{\mathsf{B}\mathbb{G}_m}^\otimes$
be the object corresponding to
$\mathbf{K}_n$ in the proof of Lemma~\ref{ppp}.
Let $A(i)$ be the image of $\mathsf{1}_{\mathsf{B}\mathbb{G}_m}(i)$
under the natural functor
$\Mod_{\mathsf{B}\mathbb{G}_m} \to \Mod_A(\Mod_{\mathsf{B}\mathbb{G}_m}^\otimes)$.
Unwinding the definition of $h$ and using adjunctions,
we see that
\[
\Map_{\Mod_{A}^\otimes(\Mod_{\mathsf{B}\mathbb{G}_m}^\otimes)}(A(i),A(j))\to \Map_{\Mod^\otimes_{[\Spec \overline{A}/\mathbb{G}_m]}}(h(A(i)),h(A(j)))
\]
can be identified with
\begin{eqnarray*}
\Map_{\Mod_{A}(\Mod_{\mathsf{B}\mathbb{G}_m}^\otimes)}(A(i),A(j)) &\simeq&  \Map_{\Mod_{A}(\Mod_{\mathsf{B}\mathbb{G}_m}^\otimes)}(A(i-j),A) \\
&\simeq& \Map_{\Mod_{\mathsf{B}\mathbb{G}_m}}(\mathsf{1}_{\mathsf{B}\mathbb{G}_m}(i-j),A) \\
&\simeq& \Map_{\Mod_{[\Spec \overline{A}/\mathbb{G}_m]}}(\pi^*(\mathsf{1}_{\mathsf{B}\mathbb{G}_m}(i-j)),\mathsf{1}_{[\Spec \overline{A}/\mathbb{G}_m]}) \\
&\simeq& \Map_{\Mod_{[\Spec \overline{A}/\mathbb{G}_m]}}(\mathsf{1}_{[\Spec \overline{A}/\mathbb{G}_m]}(i),\mathsf{1}_{[\Spec \overline{A}/\mathbb{G}_m]}(j)).
\end{eqnarray*}
Note that $A(i)$ and $h(A(i))$ are compact for each $i$,
and $h$ preserves small colimits.
The stable presentable $\infty$-category
$\Mod_{A}(\Mod_{\mathsf{B}\mathbb{G}_m}^\otimes)$
is generated by $\{A(i)\}_{i\in \ZZ}$, that is, $\Mod_{A}^\otimes(\Mod_{\mathsf{B}\mathbb{G}_m}^\otimes)$ is
the smallest stable subcategory which contains
the set $\{A(i)\}_{i\in \ZZ}$ of objects and admits filtered colimits.
Therefore for any $N\in \Mod_{A}(\Mod_{\mathsf{B}\mathbb{G}_m}^\otimes)$,
\[
\Map_{\Mod_{A}(\Mod_{\mathsf{B}\mathbb{G}_m}^\otimes)}(A(i),N) \to \Map_{\Mod^\otimes_{[\Spec \overline{A}/\mathbb{G}_m]}}(h(A(i)),h(N))
\]
is an equivalence in $\mathcal{S}$.
Furthermore, it follows from the fact that $h$ is colimit-preserving that
for any $M,N\in \Mod_{A}(\Mod_{\mathsf{B}\mathbb{G}_m}^\otimes)$,
\[
\Map_{\Mod_{A}(\Mod_{\mathsf{B}\mathbb{G}_m}^\otimes)}(M,N) \to \Map_{\Mod^\otimes_{[\Spec \overline{A}/\mathbb{G}_m]}}(h(M),h(N))
\]
is an equivalence in $\mathcal{S}$.
It remains to show that $h$ is essentially surjective.
To this end, note that
$\Mod_{[\Spec \overline{A}/\mathbb{G}_m]}\simeq \Ind(\mathcal{E})$
where $\mathcal{E}$ is the smallest stable subcategory
which contains $\{\mathsf{1}_{[\Spec \overline{A}/\mathbb{G}_m]}(i)\}_{i\in \ZZ}$.
To see this, since $\mathsf{1}_{[\Spec \overline{A}/\mathbb{G}_m]}(i)$
are compact, thus (by \cite[Definition 3.7]{BFN})
it is enough to observe that the right orthogonal of
$\{\mathsf{1}_{[\Spec \overline{A}/\mathbb{G}_m]}(i)\}_{i\in \ZZ}$ is zero,
where $\mathsf{1}_{[\Spec \overline{A}/\mathbb{G}_m]}(i)=\pi^*(\mathsf{1}_{\mathsf{B}\mathbb{G}_m}(i))$.
The condition that
\[
\Map_{\Mod_{[\Spec \overline{A}/\mathbb{G}_m]}}(\mathsf{1}_{[\Spec \overline{A}/\mathbb{G}_m]}(i),N) \simeq \Map_{\Mod_{\mathsf{B}\mathbb{G}_m}}(\mathsf{1}_{\mathsf{B}\mathbb{G}_m}(i),\pi_*(N))=0
\]
for any $i\in \ZZ$ implies
that $\pi_*(N)=0$. Then since $\pi_*$ is conservative we deduce that
$N=0$, as desired.
Since the set $\{\mathsf{1}_{[\Spec \overline{A}/\mathbb{G}_m]}(i)\}_{i\in \ZZ}$
of compact objects
generates $\Mod_{[\Spec \overline{A}/\mathbb{G}_m]}$
(in the above sense), thus
$\Ind(h(\mathcal{D})) \simeq \Mod_{[\Spec \overline{A}/\mathbb{G}_m]}$
(see \cite[5.3.5.11]{HTT})
where $\mathcal{D}$ is the smallest stable subcategory in
$\Mod_{A}(\Mod_{\mathsf{B}\mathbb{G}_m}^\otimes)$
which contains $\{A(i)\}_{i\in \ZZ}$.
It follows that $h$ is essentially surjective, noting that
$h$ is colimit-preserving and fully faithful.
\QED

\subsection{Tannakization and Derived stack of mixed Tate motives}
Proposition~\ref{kimon?}, \ref{yoyu} and Lemma~\ref{ppp}
allow us to identify
the realization functor
$\mathsf{R}_T:\DTM^\otimes \to \Mod_{H\KKK}^\otimes$
with
\[
\rho^*:\Mod_{[\Spec \overline{Q}/\mathbb{G}_m]}^\otimes \to \Mod_{[\Spec \overline{R}/\mathbb{G}_m]}^\otimes
\]
induced by the morphism of derived stacks
$\rho:[\Spec \overline{R}/\mathbb{G}_m]\to [\Spec \overline{Q}/\mathbb{G}_m]$.
Here $\overline{R}$ is the image of $R$ in $\CAlg(\Mod_{H\mathbf{K}}^\otimes)$.

Observe that
$[\Spec \overline{R}/\mathbb{G}_m]\simeq \Spec H\KKK$.
To see this, note that by the property of the realization functor
the composite of left adjoint functors
\[
\Mod_{\mathsf{B}\mathbb{G}_m} \to\Mod_{H\mathbf{K},\ZZ} \to
\Mod_{Q}(\Mod_{H\mathbf{K},\ZZ}^\otimes)\simeq \DTM \to \Mod_{H\mathbf{K}}
\]
is equivalent to the forgetful functor
(since $\Mod_{\mathsf{B}\mathbb{G}_m}\to \Mod_{H\KKK}$ is $H\KKK$-linear,
the restriction to the full subcategory of the degree zero part of
$\Mod_{H\mathbf{K},\ZZ}$ is equivalent to the identity functor,
and moreover for any $i\in \ZZ$
the restriction to the degree $i$ part is equivalent to
the identity $\Mod_{H\KKK}\to \Mod_{H\KKK}$).
And its right adjoint functor sends $\mathbf{1}_{H\mathbf{K}}$
to the object $R$ of the form $(\ldots, \mathbf{1}_{H\mathbf{K}}, \mathbf{1}_{H\mathbf{K}}, \mathbf{1}_{H\mathbf{K}},\ldots)$ which belongs to $\CAlg(\Mod_{H\mathbf{K}}^\otimes)$.
By using adjunction maps and the fact that the above composite is symmetric monoidal,
we easily see that $R$ can be viewed as the coordinate ring of $\mathbb{G}_m$
endowed with the action of $\mathbb{G}_m$, determined by
the multiplication $\mathbb{G}_m\times \mathbb{G}_m\to \mathbb{G}_m$.
Hence $[\Spec \overline{R}/\mathbb{G}_m]\simeq \Spec H\KKK$.

We refer to $[\Spec \overline{Q}/\mathbb{G}_m]$
and $\rho:\Spec H\KKK\to [\Spec \overline{Q}/\mathbb{G}_m]$
as the derived stack of mixed Tate motives and the point determined by the mixed Weil cohomology $E$ respectively.

\begin{Theorem}
\label{tatemain}
Let $\mathsf{MTG}$ be 
the derived affine group scheme over $H\KKK$ which represents
the automorphism group functor of $\mathsf{R}_T:\DTM_{\vee}^\otimes\to \PMod_{H\KKK}^\otimes$, that is, the tannakization.
Then $\mathsf{MTG}$ is equivalent to the derived affine group scheme
arising from the \v{C}ech nerve of $\rho:\Spec H\KKK\to [\Spec \overline{Q}/\mathbb{G}_m]$.
\end{Theorem}

\Proof
Apply Corollary~\ref{main} to $\rho$.
\QED

\subsection{Cycle complex and $Q$}

We describe the ($\ZZ$-graded) complex $Q$ in terms of Bloch's cycle complexes. 
We here regard $Q$ as the object
in the $\infty$-category $\Mod_{H\KKK,\ZZ}$.
(The results of this subsection will not be used in other Sections
and the reader may skip.)

For this purpose, we need an explicit right adjoint functor
$r:\DM_{\ZZ}\to \Mod_{H\KKK,\ZZ}$ of $l:\Mod_{H\KKK,\ZZ}\to \DM_{\ZZ}$.
To this end, recall the Quillen adjoint pair
\[
\mathsf{1}\otimes(-):\textup{Comp}(\mathcal{A})\rightleftarrows \textup{DM}^{eff}(k):\Gamma
\]
where the right-hand side is the model category in \cite[Example 4.12]{CD1} (cf. Section 6.1) and the left adjoint functor carries
a complex $M$ to the tensor product $\mathsf{1}\otimes M$
with the (cofibrant) unit $\mathsf{1}$ of $\textup{DM}^{eff}(k)$.
Here the tensor product $\mathsf{1}\otimes M$ is considered to be
a complex of sheaves with transfers $U\mapsto L(\Spec k)(U)\otimes_{\KKK} M$.
The right adjoint functor sends a complex of Nisnevich sheaves 
with transfers $P$ to the complex $\Gamma(P)$ of sections at $\Spec k$.
Let $F$ be a Nisnevich sheaf with transfers.
Let $\triangle^\bullet$ be the cosimplicial scheme
where $\triangle^n=\Spec k[x_0,\ldots,x_n]/(\Sigma_{i=0}^{n}x_i=0)$
and the $j$-th face $\triangle^n\hookrightarrow \triangle^{n+1}$
is determined by $x_j=0$ (see e.g. \cite{MVW}).
We then have the Suslin complex $C_*(F)$ in $\textup{DM}^{eff}(k)$, that
is the complex of sheaves with transfers,
defined by $X\mapsto F(\triangle^\bullet\times_kX)$ (take the Moore complex).

\begin{Lemma}
\label{cycle0}
Let $F$ be a Nisnevich sheaf with transfers.
Let $F\to F'$ be a fibrant replacement in $DM^{eff}(k)$.
Then the global section $F'(\Spec k)$
is quasi-isomorphic to $C_*(F)(\Spec k)$.
\end{Lemma}

\Proof
It is well-known
that the natural morphism $F\to C_*(F)$ is a weak equivalence
in $\textup{DM}^{eff}(k)$ (cf. \cite[14.4]{MVW}).
Let $C_*(F)\to C_*(F)'$ be a fibrant replacement.
Then $F'\to C_*(F)'$ (induced by the functorial fibrant replacement)
is a weak equivalence.
According to \cite[2.19, 13.8]{MVW}, cohomology sheaves of $C_*(F)$ are homotopy invariant.
By \cite[13.8, 14.8]{MVW} and the definition of $\mathbb{A}^1$-local objects
\cite[4.12]{CD1}, $C_*(F)$, $C_*(F)'$ and $F'$
are $\mathbb{A}^1$-local.
Thus both $C_*(F)\to C_*(F)'$ and $F'\to C_*(F)'$
induce isomorphisms of cohomology sheaves.
Therefore, taking the Nisnevich topology of $\Spec k$ into account, we deduce
that $C_*(F)(\Spec k)$ is quasi-isomorphic to $F'(\Spec k)$.
\QED

For an equidimensional scheme $X$ over $k$, we denote by $z^n(X,*)$ the Bloch's cycle complex of $X$ (cf. e.g. \cite[Lecture 17]{MVW}).

\begin{Corollary}
\label{cycle1}
Let $n\ge 0$.
The total right Quillen derived functor
$\mathbb{R}\Gamma$ sends $\KKK(n)$ to a complex which is quasi-isomorphic
to $z^n(\Spec k,*)[-2n]$.
\end{Corollary}

\Proof
The comparison theorems \cite[16.7, 19.8]{MVW} together with
Lemma~\ref{cycle0}
imply that
$\mathbb{R}\Gamma(\KKK(n))$ is quasi-isomorphic
to $z^n(\mathbb{A}^n,*)[-2n]$, where
$\mathbb{A}^n$ is the $n$-dimensional affine
space. The homotopy invariance of higher Chow groups
(cf. \cite[17.4 (4)]{MVW}) shows that
$z^n(\mathbb{A}^n,*)[-2n]$ is quasi-isomorphic to
$z^n(\Spec k,*)[-2n]$.
\QED

\begin{Remark}
\label{cycle2}
Let $n$ be a negative integer. 
Then every morphism from $\KKK$ to $\KKK(n)[i]$ in
$\DM$ is null-homotopic for any $i\in \ZZ$.
Thus by adjunction, the right adjoint functor of the canonical functor
$\Mod_{H\KKK}\to \DM$ carries $\KKK(n)$ to zero
in $\Mod_{H\KKK}$.
\end{Remark}

\begin{Proposition}
\label{cyclemain}
Let $Q_n \in \Mod_{H\KKK}$ denote the complex of the $n$-th degree of $Q\in \Mod_{H\KKK,\ZZ}$ (it is not the homological degree).
Then $Q_n$ is equivalent to $z^n(\Spec k,*)[-2n]$
for any $n\ge0$, and $Q_n\simeq 0$ for $n<0$.
\end{Proposition}

\Proof
Recall that $Q$ is the image of
\[
\KKK(*):=(\ldots,\KKK(-1) ,\KKK(0),\KKK(1),\ldots)
\]
by
$r:\DM_{\ZZ}\to \Mod_{H\KKK,\ZZ}$ (we adopt the notation in Section 6.2).
The natural functor $\Sigma^{\infty}:\DM^{eff}\to \DM$ is fully faithful
by Voevodsky's cancellation theorem, and thus the right adjoint $\Omega^{\infty}:\DM\to \DM^{eff}$
sends $\KKK(i)$ to $\KKK(i)$ for $i\ge0$.
Now our claim follows from Corollary~\ref{cycle1} and Remark~\ref{cycle2}.
\QED

\section{Mixed Tate motives assuming Beilinson-Soul\'e vanishing conjecture}

In this Section, we adopt the notation in Section 6.
Contrary to the previous Section, in this Section we will assume Beilinson-Soul\'e vanishing conjecture for
the base field $k$; the motivic cohomology
\[
H^{n,i}(\Spec k, \KKK)
\]
is zero for $n\le 0$, $i>0$.
Here $H^{n,i}(\Spec k, \KKK)$ denotes the motivic cohomology
(following the notation in \cite[Definition 3.4]{MVW}).
What we need is that
this condition
imply
that $Q$ is cohomologically connective,
that is, $\pi_n(\overline{Q})=0$ for $n>0$,
and $\pi_0(\overline{Q})=\KKK$.
For example, Beilinson-Soul\'e vanishing conjecture holds
when $k$ is a number field.
The goal of this Section is to prove Theorem~\ref{comparison}
which relates our tannakization $\mathsf{MTG}$ of $\DTM_\vee^\otimes$
with the Galois group of mixed Tate motives constructed by
Bloch-Kriz \cite{BK}, Kriz-May \cite{KM}, Levine
\cite{Lev} (each group scheme is known to be equivalent to one another)
under this vanishing conjecture.

\subsection{Motivic $t$-structure on $\DTM$}
Under Beilinson-Soul\'e vanishing conjecture, one can
define motivic $t$-structure on $\DTM$, as proved by
Levine \cite{Lev} and Kriz-May \cite{KM}.
We will construct a $t$-structure in our setting
(we do not claim any originality).

We fix our convention on $t$-structures.
Let $\mathcal{C}$ be a stable $\infty$-category.
A $t$-structure on $\mathcal{C}$ is a
$t$-structure on the triangulated category $\textup{h}(\mathcal{C})$
(the homotopy category is naturally endowed with the structure of
triangulated category, see \cite[Chapter 1]{HA}).
That is to say, a pair of full subcategories $(\mathcal{C}_{\ge0},\mathcal{C}_{\le0})$ of $\mathcal{C}$ such that
\begin{itemize}

\item $\mathcal{C}_{\ge0}[1]\subset \mathcal{C}_{\ge0}$ and $\mathcal{C}_{\le0}[-1]\subset \mathcal{C}_{\le0}$,

\item for $X\in \mathcal{C}_{\ge0}$ and $Y\in \mathcal{C}_{\le0}$, the hom group $\Hom_{\textup{h}(\mathcal{C})}(X, Y[-1])$ is zero,

\item for $X\in \mathcal{C}$, there exists a distinguished triangle
\[
X'\longrightarrow X \longrightarrow X''
\]
in $\textup{h}(\mathcal{C})$
such that $X'\in \mathcal{C}_{\ge0}$ and $X''\in \mathcal{C}_{\le0}[-1]$.
\end{itemize}
We here assume that full subcategories are stable under equivalences.
We use homological indexing. Our reference on $t$-structure
is \cite{HA} and \cite{KS}.
We shall write $\mathcal{C}_{\ge n}$ and $\mathcal{C}_{\le n}$ for
$\mathcal{C}_{\ge 0}[n]$ and $\mathcal{C}_{\le 0}[n]$ respectively.
We denote by $\tau_{\ge n}$ the right adjoint to $\mathcal{C}_{\ge n}\subset \mathcal{C}$.
Similarly, we denote by $\tau_{\le n}$ the left adjoint to $\mathcal{C}_{\le n}\subset \mathcal{C}$.

\vspace{2mm}

Let $\mathsf{R}_T:\DTM \to \Mod_{H\KKK}$
be the realization functor of a fixed mixed Weil theory $E$.
Let $(\Mod_{H\KKK,\ge0}, \Mod_{H\KKK,\le0})$ be the standard $t$-structure
of $\Mod_{H\KKK}$ such that
$X$ belongs to $\Mod_{H\KKK,\ge0}$ (resp. $\Mod_{H\KKK,\le 0}$)
exactly when
the homotopy group $\pi_n(X)$ of the underlying spectra
is zero for $n<0$ (resp. $n>0$).

\begin{Proposition}
\label{motivic-t}
Let 
\[
\DTM_{\vee,\ge0}:=\mathsf{R}_T^{-1}(\Mod_{H\KKK,\ge0})\cap \DTM_\vee\ \ 
\textup{and}\ \ \DTM_{\vee,\le0}:=\mathsf{R}_T^{-1}(\Mod_{H\KKK,\le0})\cap \DTM_\vee.
\]
Then the pair $(\DTM_{\vee,\ge0}, \DTM_{\vee,\le0})$
is a bounded $t$-structure on $\DTM_\vee$.
(Of course, the realization functor is $t$-exact.)
\end{Proposition}

\Proof
Since $\mathsf{R}_T$ is exact, $\DTM_{\vee,\ge0}[1]\subset \DTM_{\vee,\ge0}$
and $\DTM_{\vee,\le0}[-1]\subset \DTM_{\vee,\le0}$.

We next claim that the realization functor induces
a conservative functor $\DTM_\vee\to \Mod_{H\KKK}$.
(Recall again that
an exact functor $p:\mathcal{K}\to \mathcal{L}$ between stable $\infty$-categories is said
to be conservative if for any $K\in \mathcal{K}$,
$p(K)\simeq 0$ implies that $K\simeq 0$.)
Note that the realization functor $\DTM\simeq \Mod_{[\Spec \overline{Q}/\mathbb{G}_m]}\stackrel{\rho^*}{\to} \Mod_{H\KKK}$
is induced by $\rho:\Spec H\KKK\to [\Spec \overline{Q}/\mathbb{G}_m]$
(see Section 6.4).
The morphism $\rho$ extends to $\overline{\rho}:\Spec H\KKK\to \Spec \overline{Q}$.
Thus the realization functor is decomposed into
\[
\DTM\simeq \Mod_{[\Spec \overline{Q}/\mathbb{G}_m]}\to \Mod_{\Spec\overline{Q}} \stackrel{\overline{\rho}^*}{\to} \Mod_{H\KKK}.
\]
By the definition, the pullback of the projection
$\Mod_{[\Spec \overline{Q}/\mathbb{G}_m]}\to \Mod_{\Spec\overline{Q}}$
is conservative.
The stable $\infty$-category $\Mod_{\overline{Q}}$ admits a $t$-structure
$(\Mod_{\overline{Q},\ge0},\Mod_{\overline{Q},\le0})$
such that $X$ in $\Mod_{\overline{Q}}$ belongs to $\Mod_{\overline{Q},\le0}$
if and only if $\pi_n(X)=0$ for $n>0$ (see, \cite[VIII, 4.5.4]{DAGn}).
According to \cite[VIII, 4.1.11]{DAGn}, the composite
$\bigcup_{n\in \ZZ}\Mod_{\overline{Q},\le n}\to \Mod_{H\KKK}$
is conservative.
Observe that
every object $X\in \PMod_{\overline{Q}}$ lies in $\bigcup_{n\in \ZZ}\Mod_{\overline{Q},\le n}$.
To see this, note that $\PMod_{\overline{Q}}$ is the
smallest stable subcategory which contains $\overline{Q}$ and is closed
under retracts. Since $\overline{Q}$ belongs to
$\bigcup_{n\in \ZZ}\Mod_{\overline{Q},\le n}$
and $\bigcup_{n\in \ZZ}\Mod_{\overline{Q},\le n}$ is closed under
retracts, we see that
$\PMod_{\overline{Q}}\subset \bigcup_{n\in \ZZ}\Mod_{\overline{Q},\le n}$.
Therefore the
composite $\DTM_\vee\simeq \PMod_{[\Spec \overline{Q}/\mathbb{G}_m]}\to \Mod_{H\KKK}$
is conservative.
By using this fact, we verify the second condition of the definition of $t$-structure.

It remains to show the third condition of $t$-structure.
For this purpose, note first that
if $Z\subset \Mod_{[\Spec \overline{Q}/\mathbb{G}_m]}$ denotes the inverse image of $\bigcup_{n\in \ZZ}\Mod_{\overline{Q},\le n}$  and  $f:Z \to \Mod_{H\KKK}$ denotes the restriction of the realization functor, we have
$f^{-1}(\PMod_{H\KKK})=\PMod_{[\Spec \overline{Q}/\mathbb{G}_m]}$.
Clearly, $f^{-1}(\PMod_{H\KKK})\supset\PMod_{[\Spec \overline{Q}/\mathbb{G}_m]}$ since the realization functor is symmetric monoidal.
An object in $\Mod_{[\Spec \overline{Q}/\mathbb{G}_m]}$
is dualizable if and only if its image in $\Mod_{\overline{Q}}$ is dualizable.
Thus it is enough to show that
$g^{-1}(\PMod_{H\KKK})=\PMod_{\overline{Q}}$ where 
$g:\bigcup_{n\in \ZZ}\Mod_{\overline{Q},\ge n}\to \Mod_{H\KKK}$.
According to \cite[VIII 4.5.2 (7)]{HA}, we have the natural symmetric monoidal
fully faithful functor
$\bigcup_{n\in \ZZ}\Mod_{\overline{Q},\le n} \to \lim_{\overline{Q}\to B}\Mod_{B}$
where $B$ run over connective commutative ring spectra under $\overline{Q}$.
An object $M\in \lim_{\overline{Q}\to B}\Mod_{B}$ belongs to
its essential image if and only if the image $M(H\KKK)$ of $M$ in $\Mod_{H\KKK}$ under the natural projection
has trivial homotopy groups $\pi_m(M(H\KKK))=0$ for sufficiently large $m>>0$.
Note that every morphism $\overline{Q}\to B$ factors through
$\overline{Q}\to H\KKK$ since $\overline{Q}$ is cohomologically connected.
Consequently, we deduce that $g^{-1}(\PMod_{H\KKK})\simeq \lim_{\overline{Q}\to B}\PMod_{B}$.
Thus all objects in $g^{-1}(\PMod_{H\KKK})$ are dualizable.
It follows that $g^{-1}(\PMod_{H\KKK})=\PMod_{\overline{Q}}$.
Next consider
\[
\Mod_{[\Spec \overline{Q}/\mathbb{G}_m],\ge 0}:= \Mod_{[\Spec \overline{Q}/\mathbb{G}_m]}\times_{\Mod_{\overline{Q}}}\Mod_{\overline{Q},\ge0}.
\]
Then this category is presentable, by \cite[5.5.3.13]{HTT}.
Define $\Mod_{[\Spec \overline{Q}/\mathbb{G}_m],\le 0}$
by replacing $\ge0$ on the right-hand side by $\le0$.
Then the comonad of $\Mod_{[\Spec \overline{Q}/\mathbb{G}_m]}\rightleftarrows \Mod_{\overline{Q}}$ is given by $M \mapsto M\otimes_{H\KKK}H\KKK[t^{\pm}]$
(it is checked by using the right adjointability; Lemma~\ref{bbb}).
Therefore we can apply \cite[VII 6.20]{DAGn} to deduce that
\[
(\Mod_{[\Spec \overline{Q}/\mathbb{G}_m],\ge 0},\Mod_{[\Spec \overline{Q}/\mathbb{G}_m],\le 0})
\]
is a $t$-structure.
Note that since $\Mod_{\overline{Q}}\to \Mod_{H\KKK}$ is $t$-exact
(it follows from \cite[VIII, 4.1.10, 4.5.4 (2)]{DAGn}),
$\Mod_{[\Spec \overline{Q}/\mathbb{G}_m]} \to \Mod_{H\KKK}$ is also
$t$-exact.
We now claim that $\PMod_{[\Spec \overline{Q}/\mathbb{G}_m]}$
is stable under the truncations $\tau_{\ge0},\tau_{\le0}$.
Let $M\in \PMod_{[\Spec \overline{Q}/\mathbb{G}_m]}$.
Then $\tau_{\ge 0}M$ and $\tau_{\le0}M$ are contained in $Z$.
Thus, to prove that $\tau_{\ge 0}M$ and $\tau_{\le0}M$ belong to
$\PMod_{[\Spec \overline{Q}/\mathbb{G}_m]}$, it will suffice to
prove that $g(\tau_{\ge 0}M)$ and $g(\tau_{\le0}M)$ belong to $\PMod_{H\KKK}$.
Let $H_i=\tau_{\ge i}\circ \tau_{\le i}=\tau_{\le i}\circ\tau_{\ge i}$
(this notation slightly differs from the standard one). Using $t$-exactness, we have
\begin{eqnarray*}
H_i(g(\tau_{\ge 0}M))&=& g(H_i\circ \tau_{\ge0}M) \\
&=& g(\tau_{\le i}\circ \tau_{\ge i}\circ \tau_{\ge 0}M) \\
&=& g(H_i(M)) \\
&=& H_i(g(M))
\end{eqnarray*}
for $i\ge 0$.
It follows that $H_i(g(\tau_{\ge 0}M))[-i]$ is equivalent to
a finite dimensional $\KKK$-vector space,
and the set
\[
\{i\in \ZZ|\ H_i(g(\tau_{\ge 0}M))[-i]\neq 0\}
\]
is finite.
This implies that $g(\tau_{\ge 0}M)$ lies in $\PMod_{H\KKK}$.
Similarly, $g(\tau_{\le 0}M)$ lies in $\PMod_{H\KKK}$.
Therefore
for any $M\in \PMod_{[\Spec \overline{Q}/\mathbb{G}_m]}$
we have the distinguished triangle (in the level of homotopy category)
\[
\tau_{\ge 0}M \longrightarrow M \longrightarrow \tau_{\le -1}M
\]
such that $\mathsf{R}_T(\tau_{\ge 0}M)\in \Mod_{H\KKK,\ge0}$
and $\mathsf{R}_T(\tau_{\le -1}M)\in \Mod_{H\KKK,\le0}[-1]$, as desired.

Finally, this $t$-structure is clearly bounded.
\QED

\begin{Remark}
The definition of $t$-structure in Proposition~\ref{motivic-t} is
compatible with the definition of motivic $t$-structure
on the triangulated category of (all) mixed motives developed by Hanamura \cite{Ha3}
(up to an anti-equivalence).
In loc. cit., the expected motivic $t$-structure is constructed using
Grothendieck's standard conjectures, Bloch-Beilinson-Murre conjecture and Beilinson-Soul\'e
vanishing conjecture for smooth projective varieties.

In Proposition~\ref{motivic-t}, 
by the extension of coefficients
$\mathbb{Q}\to \KKK$ we can replace $\KKK$ by $\mathbb{Q}$.
\end{Remark}

We refer to $(\DTM_{\vee,\ge0},\DTM_{\vee,\le0})$ as
motivic $t$-structure on $\DTM_\vee$.
We let $\DTM_\vee^{\heartsuit}:=\DTM_{\vee,\ge0}\cap\DTM_{\vee,\le0}$ be
the heart.
At first sight, it depends on the choice of our realization functor.
But the mapping space
$\Map(\Spec H\KKK, \Spec \overline{Q})$ is connected
since $\overline{Q}$ is cohomologically connected (cf. \cite[VIII, 4.1.7]{DAGn}). Therefore $\rho^*:\Mod_{[\Spec \overline{Q}/\mathbb{G}_m]}^\otimes\to \Mod_{H\KKK}^\otimes$ is unique up to equivalence.

\vspace{2mm}

As a by-product of the proof, we have

\begin{Corollary}
\label{byproduct}
Adopt the notation used in the proof of Proposition~\ref{motivic-t}.
The realization functor induces a conservative functor
$f:\bigcup_{n\in \ZZ}\Mod_{[\Spec\overline{Q}/\mathbb{G}_m],\le n} \to \Mod_{H\KKK}$. In particular, $\DTM_\vee\to \PMod_{H\KKK}$ is conservative.
Moreover, $f^{-1}(\PMod_{H\KKK})$ coincides with $\DTM_\vee$.
\end{Corollary}

Recall $\DTM$ is compactly generated. Namely, we have
the natural equivalence $\Ind(\DTM_\circ)\simeq \Ind(\DTM_\vee)\simeq \DTM$.

\begin{Corollary}
Let $\DTM_{\ge0}:=\Ind(\DTM_{\vee,\ge0})$ and $\DTM_{\le0}:=\Ind(\DTM_{\vee,\le0})$.
Then the pair $(\DTM_{\ge 0},\DTM_{\le 0})$ is an accessible right complete $t$-structure on $\DTM$.
\end{Corollary}

\Proof
It follows from Proposition~\ref{motivic-t}, \cite[VIII, 5.4.1]{DAGn} and \cite[1.4.4.13]{HA}.
\QED

Let $(\Mod_{H\KKK}^\heartsuit)^\otimes$ be the symmetric monoidal
abelian category such that the underlying category is $\Mod_{H\KKK,\ge0}\cap\Mod_{H\KKK,\le0}$ and its symmetric monoidal structure is induced by
that of $\Mod_{H\KKK}^\otimes$.
It is (the nerve of) the symmetric monoidal category of $\KKK$-vector spaces. 
For an affine group scheme $G$ over $\KKK$
(which can be viewed as a derived affine group scheme over $H\KKK$),
we let $\textup{Rep}(G)^\otimes$ be the symmetric monoidal full subcategory
$z^{-1}((\Mod_{H\KKK}^\heartsuit)^\otimes)$ of $\Mod_{\mathsf{B}G}^\otimes$
where $z:\Mod_{\mathsf{B}G}^\otimes\to \Mod_{H\KKK}^\otimes$ is the natural projection determined by $\Spec H\KKK\to \mathsf{B}G$.
We denote by $\textup{Rep}(G)_{\vee}^\otimes$ the symmetric monoidal
full subcategory
of $\textup{Rep}(G)^\otimes$ which consists of dualizable objects.
Applying the classical
Tannaka duality by Saavedra, Deligne-Milne, Deligne \cite{Sa}, \cite{DM}, \cite{D}
to the faithful symmetric monoidal exact functor of abelain categories
$(\DTM_{\vee}^{\heartsuit})^\otimes \to (\Mod_{H\KKK}^\heartsuit)^\otimes$
induced by the realization functor, we have

\begin{Corollary}
There are an affine group scheme $MTG$ over $\KKK$
and an equivalence $(\DTM_\vee^\heartsuit)^\otimes\stackrel{\sim}{\to} \textup{Rep}(MTG)_{\vee}^\otimes$ of symmetric monoidal $\infty$-categories.
\end{Corollary}

We here give a symmetric monoidal
equivalence between the abelian category $\DTM^\heartsuit_\vee$ and
the abelian category $\mathbf{TM}_k$ which is constructed via
the axiomatic formulation in \cite{Lev}.
Let $i$ be an integer.
Let $W_{\ge i}\DTM_{gm}\subset\DTM_{gm}$ (resp. $W_{\le i}\DTM_{gm}\subset\DTM_{gm}$
be the smallest stable subcategory generated by $\mathbf{K}(n)$ for $-2n\ge i$
(resp. $\mathbf{K}(n)$ for $-2n\le i$).
Then according to \cite[Lemma 1.2]{Lev},
the pair $(W_{\ge i}\DTM_{gm},W_{\le i}\DTM_{gm})$
is a $t$-structure.
Let $\textup{Gr}^W_i:\DTM_{gm}\to W_i\DTM_{gm}:=W_{\ge i}\DTM_{gm}\cap W_{\le i}\DTM_{gm}$ be the functor $H_0$ with respect to this $t$-structure.
When $i$ is even, the $\infty$-category $W_i\DTM_{gm}$ is equivalent to
the full subcategory $\textup{h}(\PMod_{H\KKK})$ of $\textup{h}(\Mod_{H\KKK})$ spanned by
bounded complexes of $\KKK$-vector spaces whose (co)homology
are finite dimensional.
This equivalence is given by the exact functor $\textup{h}(\PMod_{H\KKK})\to W_i\DTM_{gm}$
which carries $\KKK[r]$ to $\KKK(-i/2)[r]$.
If $i$ is odd, $W_i\DTM_{gm}$ is zero.
It gives rise to a natural symmetric monoidal exact functor
$\textup{Gr}:\textup{h}(\DTM_{gm})\to \textup{h}(\Mod_{H\KKK,\ZZ})$, which sends $X$ to $\{\textup{Gr}^W_i(X)\}_{i\in \ZZ}$, of homotopy categories
(which are furthermore triangulated categories). The triangulated category
$\textup{h}(\Mod_{H\KKK,\ZZ})\simeq \Pi_{\ZZ}\textup{h}(\Mod_{H\KKK})$
has the standard $t$-structure determined by the product of
pair $(\Mod_{H\KKK,\ge0},\Mod_{H\KKK,\le0})$.
We denote it by $(\textup{h}(\Mod_{H\KKK,\ZZ})_{\ge0},\textup{h}(\Mod_{H\KKK,\ZZ})_{\le0})$. Let $\DTM_{gm,\ge0}:=\textup{Gr}^{-1}(\textup{h}(\Mod_{H\KKK,\ZZ})_{\ge0})$ and $\DTM_{gm,\le0}:=\textup{Gr}^{-1}(\textup{h}(\Mod_{H\KKK,\ZZ})_{\le0})$. Then by \cite[Theorem 1.4]{Lev},
we have:

\begin{Lemma}[\cite{Lev}]
The pair $(\DTM_{gm,\ge0},\DTM_{gm,\le0})$ is a bounded $t$-structure, and $\textup{Gr}$ is $t$-exact and conservative.
\end{Lemma}
Let $\mathbf{TM}_k$ be its heart.

\begin{Lemma}
The realization functor $\mathsf{R}_{gm}:\DTM_{gm}\to \Mod_{H\KKK}$ (induced by $\mathsf{R}_T:\DTM\to \Mod_{H\KKK}$) is $t$-exact.
\end{Lemma}

\Proof
We will show that the essential image of $\DTM_{gm,\le0}$ is contained in
$\Mod_{H\KKK,\le0}$.
The dual case is similar.
Let $X\in \DTM_{gm,\le0}$.
Let $m$ be the cardinal of the set of integers $i$ such that
$H_i(X)[-i]$ is not zero (recall our (nonstandard) notation $H_i=\tau_{\le i}\circ \tau_{\ge i}$).
We proceed by induction on $m$.
If $m=0$, we conclude that $X\simeq 0$
(since the $t$-structure on $\DTM_{gm}$ is bounded).
Hence this case is clear.
By \cite[Theorem 1.4 (iii)]{Lev}
we see that the essential image of $\mathbf{TM}_k$ is contained in
$\Mod_{H\KKK}^{\heartsuit}$.
Hence the case $m=1$ follows.
Suppose that our claim holds for $m\le n$.
To prove the case when $m=n+1$, consider
the distinguished triangle
\[
H_{i}(X)\to X \to \tau_{\le i-1}X
\]
where $i$ is the largest number such that $H_i(X)[-i]\neq 0$.
Note that the functor $\DTM_{gm}\to \Mod_{H\KKK}$ is exact,
and the images of $H_{i}(X)$ and $\tau_{\le i-1}X$
is contained in $\Mod_{H\KKK,\le0}$. Thus
we conclude that the image of $X$ is also contained in $\Mod_{H\KKK,\le0}$.
\QED

\begin{Lemma}
\label{gacchiri}
$\DTM_{gm,\ge 0}$ (resp. $\DTM_{gm,\le 0}$) is the inverse image of
$\Mod_{H\KKK,\ge0}$ (resp. $\Mod_{H\KKK,\le0}$)
under $\mathsf{R}_{gm}:\DTM_{gm}\to \Mod_{H\KKK}$.
\end{Lemma}

\Proof
We will treat the case $\DTM_{gm,\le 0}$. Another case is similar.
We have already prove that
$\mathsf{R}_{gm}$ is $t$-exact in the previous Lemma.
It will suffice to show that if $X$ does not belong to
$\DTM_{gm,\le0}$, then
$\mathsf{R}_{gm}(X)$ does not lies in $\Mod_{H\KKK,\le0}$.
For such $X$, there exists $i\ge 1$ such that $H_i(X)\neq 0$.
According to Corollary~\ref{byproduct},
$\mathsf{R}_{gm}$ is conservative.
Combined with the $t$-exactness, we deduce that
$H_i(\mathsf{R}_{gm}(X))[-i]\neq 0$.
This implies that $\mathsf{R}_{gm}(X)$ is not in $\Mod_{H\KKK,\le0}$,
as required.
\QED

By Lemma~\ref{gacchiri}, we have a $t$-exact fully faithful functor
$\DTM_{gm}\to \DTM_{\vee}$, and it
induces a natural fully faithful functor
$\mathbf{TM}_k\to \DTM_\vee^{\heartsuit}$ between (nerves of) symmetric monoidal abelian categories.

\begin{Proposition}
The natural inclusion $\mathbf{TM}_k\to \DTM_\vee^{\heartsuit}$
is an equivalence.
\end{Proposition}

\Proof
Since $\mathbf{TM}_k$ is (the nerve of) an abelian category, and in particular it is
idempotent complete, thus it is enough to prove that
$\mathbf{TM}_k\to \DTM_\vee^{\heartsuit}$
is an idempotent completion.
Recall that $\DTM_{gm}\to \DTM_{\vee}$ is an idempotent completion.
Let $X\in \mathbf{TM}_k$.
The direct summand of $X$ (which automatically belongs to $\DTM_\vee$)
lies in $\DTM_\vee^\heartsuit$ by the definition of $t$-structure of $\DTM_{\vee}$.
Conversely, if $Y\in \DTM_\vee^\heartsuit$, then there
exists $X\in \DTM_{gm}$ such that $Y$ is equivalent to
a direct summand of
$X$. Then $Y$ is a direct summand of $H_0(X)\in \mathbf{TM}_k$
(note that we here use the $t$-exactness of $\DTM_{gm}\to \DTM_\vee$).
Consequently, $\mathbf{TM}_k\to \DTM_\vee^{\heartsuit}$ is an idempotent completion.
\QED

\begin{Corollary}
The Tannaka dual of $\mathbf{TM}_k$ (endowed with the realization functor)
is equivalent to $MTG$.
\end{Corollary}

\begin{Warning}
In \cite{Lev}, one works over rational coefficients.
In this note, we work over $\mathbf{K}$.
Therefore $MTG$ is the base change of the Tannaka dual of the abelian
category of mixed Tate motives in \cite{Lev} over $\QQ$ to $\mathbf{K}$.
\end{Warning}

\subsection{Completion and locally dimensional $\infty$-category}

Let $\DTM^\otimes \to \oDTM^\otimes$ be the left completion of $\DTM^\otimes$
with respect to the $t$-structure $(\DTM_{\ge 0},\DTM_{\le 0})$
(we refer the reader to \cite[1.2.1.17]{HA} and \cite[VIII, 4.6.17]{DAGn}
for the notions of left completeness and left completion).
It is symmetric monoidal, $t$-exact and colimit-preserving.
Here, the $\infty$-category $\oDTM$ is the limit of the diagram indexed by $\ZZ$
\[
\cdots \to \DTM_{\le n+1} \stackrel{\tau_{\le n}}{\to} \DTM_{\le n} \stackrel{\tau_{\le n-1}}{\to} \DTM_{\le n-1} \stackrel{\tau_{\le n-2}}{\to} \cdots
\]
of $\infty$-categories.
Note that according to \cite[3.3.3]{HTT}
the $\infty$-category $\oDTM$ can be identified with
the full subcategory of $\Fun(\NNNN(\ZZ),\DTM)$
spanned by functors $\phi:\NNNN(\ZZ)\to \DTM$ such that
\begin{itemize}
\item for any $n\in \ZZ$, $\phi([n])$ belongs to $\DTM_{\le -n}$,

\item for any $m\le n\in \ZZ$, the associated map $\phi([m])\to \phi([n])$
gives an equivalence $\tau_{\le -n}\phi([m])\to \phi([n])$.
\end{itemize}
Let $\oDTM_{\ge0}$ (resp. $\oDTM_{\le0}$) be the full subcategory of $\oDTM$
spanned by $\phi:\NNNN(\ZZ)\to \DTM$ such that $\phi([n])$ belongs to
$\DTM_{\ge0}$ (resp. $\DTM_{\le0}$) for each $n\in \ZZ$.
The functor $\DTM\to \oDTM$ induces an equivalence $\DTM_{\le0}\to \oDTM_{\le0}$.
The pair $(\oDTM_{\ge0}, \oDTM_{\le0})$ is an accessible, left complete and right complete $t$-structure of $\oDTM$.

\begin{Proposition}
\label{sayyes}
The followings hold.
\begin{enumerate}
\renewcommand{\labelenumi}{(\roman{enumi})}
\item $\oDTM_{\le0}$ is closed under filtered colimits.

\item The unit $\mathsf{1}$ belongs to the heart $\oDTM^{\heartsuit}:=\oDTM_{\ge0}\cap\oDTM_{\le0}$.

\item $\oDTM_{\ge 0}$ and $\oDTM_{\le0}$ are closed under the tensor product
$\oDTM \times \oDTM\to \oDTM$.

\item The unit $\mathsf{1}$ is compact in $\DTM_{\le n}$ for each $n\ge0$.

\item There exists a full subcategory $\oDTM^{\heartsuit}_{\textup{fd}}$ of $\oDTM^{\heartsuit}$ such that every object in $\oDTM^{\heartsuit}_{\textup{fd}}$
has the dual in $\oDTM^{\heartsuit}_{\textup{fd}}$, and $\oDTM^{\heartsuit}_{\textup{fd}}$ generates $\oDTM^{\heartsuit}$ under filtered colimits. 

\item $\pi_0(\Map_{\oDTM}(\mathsf{1},\mathsf{1}))=\mathbf{K}$.

\item For any $X\in \oDTM^{\heartsuit}_{\textup{fd}}$, the composite
\[
\mathsf{1}\to X\otimes X^\vee \to \mathsf{1}
\]
of the coevaluation map and the evaluation map corresponds to 
a nonnegative integer $\dim(X)\in \ZZ\subset \mathbf{K}$.
\end{enumerate}
\end{Proposition}

\Proof
By our construction and $\oDTM_{\le 0}=\DTM_{\le0}$, (i) is clear.
Since the unit of $\DTM$ lies in $\DTM^{\heartsuit}:=\DTM_{\ge0}\cap \DTM_{\le 0}$,
(ii) follows.

Next we will prove (iii).
Taking account into
the definition of $\DTM_{\vee,\le0}$ and $\DTM_{\vee,\ge0}$
and the conservatibity of $\mathsf{R}_T:\DTM_{\vee}\to \Mod_{H\mathbf{K}}$,
we see that $\DTM_{\vee,\le0}$ and $\DTM_{\vee,\ge0}$
are stable under tensor operation.
By $\Ind(\DTM_{\vee,\le0})=\oDTM_{\le 0}$
and the
tensor operation presrves colimits in each variable, we deduce that
$\oDTM_{\le 0}$ is stable under tensor operation.
Since $\DTM_{\ge0}$ is stable under tensor operation
of $\DTM$, by definition we also see that $\oDTM_{\ge 0}$
is stable under tensor operation.

The unit $\mathsf{1}$ is compact in $\DTM$,
and so is in $\DTM_{\le n}$ for any $n\in \ZZ$.
Noting that $\oDTM_{\le n}=\DTM_{\le n}$,
we have (iv).

To prove (v), note first that $\DTM\to \oDTM$ induces
equivalences
$\bigcup_{n\in \ZZ}\DTM_{\le n}\to\bigcup_{n\in \ZZ}\oDTM_{\le n}$
and $\DTM^\heartsuit \simeq \oDTM^\heartsuit$.
In particular, $\DTM_\vee\to \oDTM$ is fully faithful.
Let $X\in \oDTM^\heartsuit=\DTM^\heartsuit$.
Then $X$ is the filtered colimit
of a diagram $I\to \DTM_{\vee,\ge0}$ in $\DTM$ (or in $\oDTM$);
$\textup{colim}_{\lambda\in I}X_\lambda \simeq X$. (Recall
$\DTM^\heartsuit\subset \Ind(\DTM_{\vee,\ge0})$.)
Note that $X_{\lambda}\in \DTM_\vee$ and by definition
$\DTM_{\vee,\ge0}$, $\DTM_{\vee,\le0}$ and
its heart are stable under the tensor operation.
The heart is stable under taking dual objects.
It follows
that $\tau_{\le0}(X_{\lambda})=H_0(X_\lambda)$ is dualizable, that is, it belongs to
$\DTM_{\textup{fd}}^\heartsuit:=\DTM_\vee\cap \DTM^\heartsuit$
and 
the dual of $H_0(X_\lambda)$ lies in $\DTM_{\textup{fd}}^\heartsuit$.
Sience $\tau_{\le0}$ is a left adjoint, thus the natural morphism
$\textup{colim}_\lambda \tau_{\le0}(X_\lambda)\to \tau_{\le0}(\textup{colim}_\lambda X_\lambda)$ is an equivalence.
This shows that $\DTM_{\textup{fd}}^\heartsuit$
generates $\DTM^\heartsuit=\oDTM^\heartsuit$ under filtered colimits.

We remark that $H^{0,0}(\Spec k, \KKK)=\KKK$. Hence (vi) holds.
Finally, we will prove (vii). For any $X\in \oDTM^{\heartsuit}_{\textup{fd}}$, the element in $\KKK$ corresponding to the composite $\mathsf{1}\to X\otimes X^\vee \to \mathsf{1}$ is equal to the element in $\KKK$ corresponding to
$\mathsf{R}_T(\mathsf{1})\to \mathsf{R}_T(X)\otimes \mathsf{R}_T(X)^\vee\to \mathsf{R}_T(\mathsf{1})$. The latter element is nothing but the dimension of
$\mathsf{R}_T(X)$, which lies in $\ZZ$.
\QED

\begin{Remark}
Let $\mathcal{C}^\otimes$ be the symmetric monoidal
stable subcategory of $\DM^\otimes_\vee$ which is closed under
taking retracts and dual objects.
Suppose that $\mathcal{C}^\otimes$ admits a non-degenerate
$t$-structure $(\mathcal{C}_{\ge0},\mathcal{C}_{\le0})$
such that
\begin{itemize}
\item the realization functor $\mathcal{C}^\otimes\subset \DM^\otimes_{\vee} \to \Mod_{H\mathbf{K}}^\otimes$ is $t$-exact,

\item both $\mathcal{C}_{\ge0}$ and $\mathcal{C}_{\le0}$
are stable under the tensor operation $\mathcal{C}\times \mathcal{C}\to \mathcal{C}$ respectively.

\end{itemize}
As observed in \cite[1.3]{Bei}, its heart $\mathcal{C}_{\ge0}\cap \mathcal{C}_{\le0}$ is a tannakian category equipped with the realzaition functor as a fiber functor, and the realization functor $\mathcal{C}\to \Mod_{H\mathbf{K}}$
is conservative.
Let $\widehat{\mathcal{C}}_{\ge0}$ (resp. $\widehat{\mathcal{C}}_{\le0}$)
the left completion of $\Ind(\mathcal{C}_{\ge0})$
(resp. $\Ind(\mathcal{C}_{\le0})$).
Then as above the pair $(\widehat{\mathcal{C}}_{\ge0},\widehat{\mathcal{C}}_{\le0})$
is an accessible, both right complete and left complete $t$-structure 
on the left completion $\widehat{\mathcal{C}}$ of $\Ind{\mathcal{C}}$
(with respect to $(\Ind(\mathcal{C}_{\ge0}),\Ind(\mathcal{C}_{\le0}))$).
The argument of the above proof shows that the analogous assertions in Proposition~\ref{sayyes} hold also for
$(\widehat{\mathcal{C}}_{\ge0},\widehat{\mathcal{C}}_{\le0})$.
(Consequently, analogues of Corollary~\ref{loc} and Proposition~\ref{lurie}
also hold.)
\end{Remark}

\begin{Corollary}
\label{loc}
The symmetric monoidal $\infty$-category $\oDTM^\otimes$ endowed with the
$t$-structure $(\oDTM_{\ge0},\oDTM_{\le0})$ is a locally dimensional
$\infty$-category in the sense of \cite[VIII, 5.6]{DAGn}.
\end{Corollary}

To state the next result which follows from the theory
of locally dimensional
$\infty$-categories, we prepare some notation.
We say that a commutative ring spectrum $S$ is discrete if
$\pi_i(S)=0$ for $i\neq0$.
This property is equivalent to the property that
there exists a (usual) commutative ring $R$ such that
$HR\simeq S$ in $\CAlg$.
Let $\CAlg^{\textup{dis}}$ be the $\infty$-category of discrete commutative ring spectra.
The $\infty$-category $\CAlg^{\textup{dis}}$
is equivalent to the nerve of the category of (usual) commutative
rings (via Eilenberg-MacLane spectra).
Let $\mathfrak{S}:\CAlg^{\textup{dis}}\to \widehat{\SSS}$
be the functor which carries $A\in\CAlg^{\textup{dis}}$ to the space
$\Map_{\CAlg(\wCat^{\textup{L,st}})}(\oDTM^\otimes,\Mod_{A}^\otimes)$
(which can be constructed by $\Theta$ in Section 3.1 and Yoneda embedding).
Let $\xi:\CAlg^{\textup{dis}}\to \widehat{\SSS}$
be the functor which carries $A\in\CAlg^{\textup{dis}}$ to the space
$\Map_{\CAlg(\wCat^{\textup{L,st}})}(\Mod_{H\KKK}^\otimes,\Mod_{A}^\otimes)$.
Since there exists a natural equivalence
\[
\Map_{\CAlg(\wCat^{\textup{L,st}})}(\Mod_{H\KKK}^\otimes,\Mod_{A}^\otimes)\simeq \Map_{\CAlg}(H\KKK,A)
\]
(cf. \cite[Section 5]{FI}, \cite[6.3.5.18]{HA}), $\xi$ is corepresented by $H\KKK$.
We here write $\Spec H\KKK$ for $\xi$.
There is a sequence of functors $\Mod_{H\KKK}^\otimes\to \oDTM^\otimes\to \Mod_{H\KKK}^\otimes$ whose composite is equivalent to the identity.
Therefore we have $\Spec H\KKK\stackrel{\eta}{\to} \mathfrak{S}\to \Spec H\KKK$
whose composite is the identity.
Let $V:\CAlg^{\textup{dis}}\to \widehat{\SSS}$ be a functor
equipped with $V\to \Spec H\KKK$.
To $f:H\KKK\to A$ in $\CAlg^{\textup{dis}}_{H\KKK}:=(\CAlg^{\textup{dis}})_{H\KKK/}$ we associate $\{f\}\times_{\Spec H\KKK(A)}V(A)$.
It yields the functor $V_0:\CAlg^{\textup{dis}}_{H\KKK}\to \widehat{\SSS}$.
The morphism $\eta:\Spec H\KKK\to \mathfrak{S}$ induces
$\eta_0:(\Spec H\KKK)_0\to \mathfrak{S}_0$. Note that $(\Spec H\KKK)_0$ is
equivalent to the constant functor taking the value $\Delta^0$, that is, the final
object.

The following result is essentially proved by Lurie in the context of
locally dimensional $\infty$-categories (see \cite[VIII, 5.2.12, 5.6.1, 5.6.19 and their proofs]{DAGn}). We here state only the version in view of
Corollary~\ref{loc}, which fits in with
our need.

\begin{Proposition}[\cite{DAGn}]
\label{lurie}
Let $\textup{Grp}^{\textup{dis}}$ be the nerve of the category of (usual)
groups.
Consider the functor $\pi_1(\mathfrak{S}_0,\eta_0):\CAlg_{H\KKK}^{\textup{dis}}\to \textup{Grp}^{\textup{dis}}$ which is given by
$A\mapsto \pi_1(\mathfrak{S}_0(A),\eta_0)$.
Then $\pi_1(\mathfrak{S}_0,\eta_0)$ is represented by $MTG$, that is,
the Tannaka dual of $(\DTM_\vee^\heartsuit)^\otimes$.
\end{Proposition}

\subsection{Comparison theorem}

\begin{Definition}
\label{under}
Let $G:\CAlg_{H\KKK}\to \Grp(\mathcal{S})$ be a derived affine
group scheme over $H\KKK$.
Let $\pi_0:\Grp(\mathcal{S})\to \Grp^{\textup{dis}}$ be the truncation functor given by
$G\mapsto \pi_0(G)$.
If the composition
\[
G':\CAlg_{H\KKK}^{\textup{dis}} \hookrightarrow \CAlg_{H\KKK}\stackrel{G}{\to} \Grp(\mathcal{S}) \stackrel{\pi_0}{\to} \Grp^{\textup{dis}}
\]
is represented by an affine group scheme $G_0$ over $\KKK$,
we say that $G_0$ is an excellent coarse moduli space of $G$.
If there is an affine group scheme $G_0$ (considered as
$\CAlg_{H\KKK}^{\textup{dis}}\to \Grp(\mathcal{S})$)
and a morphism $G|_{\CAlg_{H\KKK}^{\textup{dis}}}\to G_0$ that is universal among morphisms into affine group scheme over $\KKK$,
we say that $G_0$ is a coarse moduli space of $G$.
We remark that an excellent coarse moduli space is a coarse moduli space.
\end{Definition}

\begin{Theorem}
\label{comparison}
Let $\mathsf{MTG}$ denote the tannakization of $\mathsf{R}_T:\DTM_\vee^\otimes\to \PMod_{H\KKK}^\otimes$ (cf. Theorem~\ref{tatemain}).
Then $MTG$ is an excellent coarse moduli space of $\mathsf{MTG}$.
\end{Theorem}

\Proof
For $A\in \CAlg^{\textup{dis}}$,
we set $\Mod_{A,\ge0}=\{X\in \Mod_A|\ \pi_i(X)=0\ \textup{for}\ i<0\}$ and
$\Mod_{A,\le0}=\{X\in \Mod_A|\ \pi_i(X)=0\ \textup{for}\ i>0\}$.
Then the pair $(\Mod_{A,\ge0},\Mod_{A,\le0})$
is an accessible, left and right complete $t$-structure.
Thus
we have
\begin{eqnarray*}
\Map_{\CAlg(\wCat^{\textup{L,st}})}^{\textup{rex}}(\oDTM^\otimes, \Mod_A^\otimes) &\simeq& \Map_{\CAlg(\wCat^{\textup{L,st}})}^{\textup{rex}}(\DTM^\otimes, \Mod_A^\otimes) \\
&\hookrightarrow& \Map_{\CAlg(\wCat)}(\DTM_\vee^\otimes, \Mod_A^\otimes)
\end{eqnarray*}
where $\Map^{\textup{rex}}$ indicates the full subcategory spanned by right $t$-exact functors, and 
the second arrow is fully faithful by Proposition~\ref{symKan}.
(The essential image consists of symmetric monoidal exact functors
which are right $t$-exact.)
Note that $\mathsf{R}_T:\DTM^\otimes\to \Mod_{H\KKK}^\otimes$ is $t$-exact,
and it belongs to $\Map_{\CAlg(\wCat^{\textup{L,st}})}^{\textup{rex}}(\DTM^\otimes, \Mod_A^\otimes)$.

Consider the automorphism group functor $\Aut(\mathsf{R}_T):\CAlg_{H\KKK}\to \Grp(\mathcal{S})$
of $\mathsf{R}_T:\DTM^\otimes_\vee \to \PMod_{H\KKK}^\otimes$
in $\CAlg(\sCat)$, (we abuse notation for $\mathsf{R}_T$).
According to Theorem~\ref{tatemain},
$\Aut(\mathsf{R}_T)$ is represented by $\mathsf{MTG}$.
On the other hand, using the above equivalence and unfolding the definition of $\pi_1(\mathfrak{S}_0,\eta_0)$
and $\Aut(\mathsf{R}_T)$,
we see that
the composite
\[
\CAlg_{H\KKK}^{\textup{dis}}\hookrightarrow \CAlg_{H\KKK}\stackrel{\Aut(\mathsf{R}_T)}{\longrightarrow} \Grp(\mathcal{S}) \stackrel{\pi_0}{\to} \Grp^{\textup{dis}}
\]
is equivalent to $\pi_1(\mathfrak{S}_0,\eta_0)$.
Combined with Proposition~\ref{lurie} we complete the proof.
\QED

\begin{Remark}
The truncation procedure given in \cite[Section 5, Truncated affine group schemes]{Tan} also allows us to construct the (usual)
affine group scheme $MTG'=\Spec H^0(\tau B)$ over $\mathbf{K}$ from $\mathsf{MTG}=\Spec B$ where $B$ is a commutative differential graded algebra.
Here we adopt the notation in \cite[Section 5]{Tan}.
This $MTG'$ coincides with the above $MTG$.
To observe this, we invoke $\tau$ in loc. cit.to
have $\Spec \tau B$, and we regard $\Spec \tau B$ as
a functor $\CAlg_{H\mathbf{K}}^{\textup{dis}}\to \Grp(\mathcal{S})$.
As functors $\CAlg_{H\mathbf{K}}^{\textup{dis}}\to \Grp(\mathcal{S})$,
$\Spec \tau B$ is equivalent to $\Spec B$.
Also, let us regard $MTG$ as a functor
$\CAlg_{H\mathbf{K}}^{\textup{dis}}\to \Grp(\mathcal{S})$.
Then we have a natural morphism $\Spec \tau B\to MTG$.
Also, the composition with $\Spec \tau B\to MTG$
induces an equivalence of spaces
\[
\Map_{\Fun(\CAlg_{H\mathbf{K}}^{\textup{dis}},\Grp(\mathcal{S}))}(MTG,F)\to \Map_{\Fun(\CAlg_{H\mathbf{K}}^{\textup{dis}},\Grp(\mathcal{S}))}(\Spec \tau B,F)
\]
for any
$F:\CAlg_{H\mathbf{K}}^{\textup{dis}}\to \Grp(\mathcal{S}^{\textup{dis}})\subset \Grp(\mathcal{S})$.
On the other hand, by the construction
in \cite{Tan} $\Spec \tau B\to MTG'$ is universal among morphisms
to usual affine group schemes over $\mathbf{K}$.
Hence we have $MTG\simeq MTG'$.
\end{Remark}

\section{Artin motives and
Absolute Galois group}

Let $G_k$ denote the absolute Galois group $\textup{Gal}(\bar{k}/k)$
with an algebraic closure $\bar{k}$ of a perfect field $k$.
For the sake of completeness,
we will construct a natural homomorphism
\[
\mathsf{MG}_E\longrightarrow G_k
\]
where $\mathsf{MG}_E$ is the derived motivic Galois group.
This represents the automorphism
functor of $\mathsf{R}_E:
\DM^\otimes_\vee\to \PMod_{H\mathbf{K}}^\otimes$.
Here $\DM^\otimes_\vee$ denotes the symmetric monoidal
full sucbategory of $\DM^\otimes$ spanned by dualizable objects.
To this end, we consider the full subcategory of $\mathsf{DM}$
which consists of Artin motives, and we will finish by proving that
its tannakization is the absolute Galois group (Proposition~\ref{toabsolute}).

Let $\Cor_{\mathbf{K},0}$ be the full subcategory of
$\Cor_{\mathbf{K}}$ spanned by smooth schemes $X$
which are \'etale over $\Spec k$.
We simply write $\Cor_0$ and $\Cor$ for $\Cor_{\mathbf{K},0}$
and $\Cor_{\mathbf{K}}$ respectively.
The classical Galois theory says
that the category of schemes which are \'etale over $k$
is equivalent to the category of finite $G_k$-sets.
Consequently, we easily see that there is a fully faithful functor
$\Cor_0\to \mathbf{K}[G_k]\textup{-Mod}$
which carries $X$ to the $\mathbf{K}$-vector space
generated by the set $X(\bar{k})$ endowed with action of $G_k$.
Here $\mathbf{K}[G_k]\textup{-Mod}$ denotes the category of
$\mathbf{K}[G_k]$-modules, i.e. abelian groups equipped with
the (left) actions of $\mathbf{K}[G_k]$.
The essential image consists of permutational representations
(see \cite[p. 216]{Vtri}).

Let $\iota:\Cor\to \Cor_0$ be the left adjoint of the inclusion
$\Cor_0 \hookrightarrow \Cor$. The funtor $\iota$
carries $X$ to the Zariski spectrum
of the integral closure of $k$ in $\Gamma(X)$.
Let 
$\PSh(\Cor_0)$
be the category of presheaves (with value of $\mathbf{K}$-vector
spaces) with transfers, that is the category of
$\mathbf{K}$-linear functors $(\Cor_0)^{op} \to \mathbf{K}\textup{-Vect}$
where 
$\mathbf{K}\textup{-Vect}$ is the category of $\mathbf{K}$-vector spaces.
Note that $\PSh(\Cor_0)$ contains $\Cor_0$ as a full subcategory by
enriched Yoneda lemma \cite{Ke}.
There is a symmetric monoidal structure on $\PSh(\Cor_0)$
which makes $\Cor_0\hookrightarrow \PSh(\Cor_0)$ symmetric monoidal
such that the tensor product
$\PSh(\Cor_0)\times \PSh(\Cor_0)\to \PSh(\Cor_0)$
preserves small colimits separately in each variable.
Such a symmetric monoidal structure is usually called Day convolution
\cite{Da}.
This exhibits $\PSh(\Cor_0)$ as a symmetric monoidal abelian category.
We define $\Sh(\Cor)$ to be the symmetric monoidal category of 
Nisnevich sheaves with transfer (see \cite{CD1}).
The composition with $\iota$ and the sheafification
induces a symmetric monoidal functor $\PSh(\Cor_0)\to \Sh(\Cor)$.
Hence it gives rise to a functor
\[
\textup{Comp}(\PSh(\Cor_0))\to \textup{Comp}(\Sh(\Cor)).
\]
Let us equip the category $\textup{Comp}(\Sh(\Cor))$ with the model structure
given in \cite[2.4]{CD1},
in which weak equivalences are quasi-isomorphisms.
We equip $\textup{Comp}(\PSh(\Cor_0))$
with the model structure in \cite[2.5]{CD1} 
by choosing the descent structure $(\mathcal{G},\mathcal{H})$ in
\cite[2.2]{CD1} as $\mathcal{G}$:= sheaves represented by objects
in $\Cor_0$, and $\mathcal{H}=\{0\}$.
Then by \cite[2.14]{CD1}
we see that the above functor is a left Quillen adjoint symmetric monoidal functor.
Hence we take their localization and have
the symmetric monoidal colimit-preserving functor
\[
\NNNN(\textup{Comp}(\PSh(\Cor_0))^c)^\otimes_\infty\to \NNNN(\textup{Comp}(\PSh(\Cor))^c)^\otimes_\infty
\]
of symmetric monoidal $\infty$-categories.
By the construction of $\mathsf{DM}$ (cf. \cite[7.15]{CD1})
there is the natural symmetric monoidal colimit-preserving
functor
\[\NNNN(\textup{Comp}(\Sh(\Cor))^c)^\otimes_\infty\to \mathsf{DM}^\otimes
\]
which is induced by the localization by $\mathbb{A}^1$-homotopy equivalence
and the stabilization by the Tate sphere.
Thus by composition we obtain the symmetric monoidal functor
\[
\mathfrak{A}: \mathsf{Q}^\otimes:=\NNNN(\textup{Comp}(\PSh(\Cor_0))^c)^\otimes_\infty\to \mathsf{DM}^\otimes\simeq \mathsf{Sp}_{\textup{Tate}}(\mathbf{HK})^\otimes.
\]
The image of the inclusion $\Cor_0\hookrightarrow \textup{Comp}(\PSh(\Cor_0))$
is contained in $\textup{Comp}(\PSh(\Cor_0)^c)$.
Let $\mathsf{Art}(k)$ be the smallest stable idempotent complete
subcategory
which contains its essential image.
Alternatively, if we let $A$ be the triangulated thick subcategory of
$\textup{h}(\mathsf{Q})$
generated by the essential image of $\Cor_0\to \textup{h}(\mathsf{Q})$, then
$\mathsf{Art}(k)\simeq \mathsf{Q}\times_{\textup{h}(\mathsf{Q})}A$.
Observe that by the elementary representation theory and the fully faithful embedding
$\Cor_0\subset G_k\textup{-Mod}$
the idempotent completion $\Cor_0^\sim$ of $\Cor_0$ (in $\PSh(\Cor_0)$)
can be identified with the abelian category of discrete representations
of $G_k$ (that is, actions $\rho:G_k\to \Aut(V)$
of $G_k$ on finite dimensional
$\mathbf{K}$-vector spaces $V$
such that $\rho$ factor through some finite quotient $G_k\to H$).
The abelian category $\Cor_0^\sim$ is semi-simple.
Hence it is easy to see that
the stable subcategory $\mathsf{Art}(k)$ of $\mathsf{Q}$ is spanned
by bounded complexes $C$ such that $C^n$ belongs to $\Cor_0^\sim$
for each $n\in \ZZ$ (indeed such complexes are cofibrant).
Note that the symmetric monoidal structure of $\mathsf{Q}^\otimes$
induces the symmetric monoidal structure of
$\mathsf{Art}(k)$.
According to \cite[3.4.1]{Vtri} and Voevodsky's calcellation theorem
together with \cite[Lemma 5.8]{Tan} we deduce:

\begin{Lemma}
The natural functor
$\mathsf{Art}(k)\to \mathsf{DM}$ is fully faithful.
\end{Lemma}

We identify $\mathsf{Art}(k)^\otimes$ with a symmetric monoidal full
subcategory of $\mathsf{DM}^\otimes$ and refer to it as the $\infty$-category of Artin motives.
We remark that $\mathsf{Art}(k)$ is contained in
the full subcategory of
$\mathsf{DM}$
spanned by compact objects.

We regard $G_k$ as a limit $\lim(\textup{Gal}(L/k))$ where $L$ run through
all finite Galois extensions $L$ of $k$.
Let $\textup{Gal}(L/k)\textup{-Perm}$ be the $\mathbf{K}$-linear category
of permutational representations.
We define $\PSh(\textup{Gal}(L/k)\textup{-Perm})$
to be the symmetric monoidal category of presheaves (with values of
$\mathbf{K}$-vector spaces)
on $\textup{Gal}(L/k)\textup{-Perm}$
in the same way as $\PSh(\Cor_0)$.
Let us equip the category
$\textup{Comp}(\PSh(\textup{Gal}(L/k)\textup{-Perm}))$
with the symmetric monoidal model structure given in \cite[2.5, 3.2]{CD1}
by choosing the descent structure
$(\mathcal{G},\mathcal{H})$
in \cite[2.2]{CD1} as
$\mathcal{G}:=$
sheaves represented by objects $\textup{Gal}(L/k)\textup{-Perm}$
and $\mathcal{H}=\{0\}$.
Let $\textup{Gal}(L/k)\textup{-Perm}^{\sim}$
be the idempotent completion of $\textup{Gal}(L/k)\textup{-Perm}$
 which can be
identified with the abelian category of finite dimensional
representations of $\textup{Gal}(L/k)$.
Let $\mathsf{A}_L$ be the stable subcategory of $\mathsf{B}_L:=\NNNN(\textup{Comp}(\PSh(\textup{Gal}(L/k)\textup{-Perm}))^c)_\infty$
spanned by bounded complexes $C$ such that
$C^n$ lies in $\textup{Gal}(L/k)\textup{-Perm}^{\sim}$
for each $n\in \ZZ$. The full subcategory $\mathsf{A}_L$ is a symmetric monoidal full subcategory of $\mathsf{B}_L^\otimes$ (spanned by
dualizable objects).
The quotient homomorphism $\pi_L:G_k\to \textup{Gal}(L/k)$
naturally induces a symmetric monoidal
functor $\PSh(\textup{Gal}(L/k)\textup{-Perm}) \to \PSh(\Cor_0)$ which is a left Kan extension of the natural functor $\xi:\textup{Gal}(L/k)\textup{-Perm}\hookrightarrow \Cor_0\hookrightarrow \PSh(\Cor_0)$. This left adjoint
is given by the formula
$M\mapsto \textup{colim}_{H_\sigma\to M}\xi(H_\sigma)$
where $H_\sigma$ is a presheaf represented by $\sigma\in \textup{Gal}(L/k)\textup{-Perm}$ and $H_\sigma\to M$ run through the overcategory $\textup{Gal}(L/k)\textup{-Perm}_{/M}$.
It gives rise to a symmetric monoidal
colimit-preserving functor
\[
\textup{Comp}(\PSh(\textup{Gal}(L/k)\textup{-Perm})) \to \textup{Comp}(\PSh(\Cor_0)).
\]
By the definition of their model structures
this is a left Quillen adjoint, and we
obtain a symmetric monoidal colimit-preserving functor
\[
\NNNN(\textup{Comp}(\PSh(\textup{Gal}(L/k)\textup{-Perm}))^c)^\otimes_\infty \to \mathsf{Q}^\otimes=\NNNN(\textup{Comp}(\PSh(\Cor_0))^c)^\otimes_\infty.
\]
Taking account of all finite Galois extensions $L$
we have
\[
f:\textup{colim}_L(\mathsf{A}^\otimes_L) \to \mathsf{Art}(k)^\otimes.
\]

\begin{Lemma}
\label{artincolimit}
The functor $f$ is an equivalence of
symmetric monoidal $\infty$-categories.
\end{Lemma}

\Proof
By \cite[3.2.3.1, 4.2.3.5]{HA} we can regard the underling $\infty$-category of
$\textup{colim}_L(\mathsf{A}^\otimes_L)$
as a colimit of the diagram of the underling $\infty$-categories
$\mathsf{A}_L$ in $\textup{Cat}_\infty$.
Moreover, the filtered colimit of stable $\infty$-categories
$\mathsf{A}_L$ is also stable \cite[1.1.4.6]{HA}.
Thus by \cite[Lemma 5.8]{Tan} it is enough to observe that
$f:\textup{colim}_L(\mathsf{A}_L) \to \mathsf{Art}(k)$
induces an equivalence of their homotopy categories.
Clearly, $f$ is essentially surjective.
By computing the hom sets in the homotopy category
we see that $f$ induces the fully faithful functor
$\textup{h}(\mathsf{A}_L)\to \textup{h}(\mathsf{Art}(k))$ of
homotopy categories for each $L$.
\QED

Let $\mathsf{R}':\Art(k)^\otimes\to \PMod_{H\mathbf{K}}^\otimes$ be the
composition of $\Art(k)^\otimes\to \mathsf{DM}^\otimes_\vee$
and the realization functor $\mathsf{R}:\mathsf{DM}^\otimes_\vee
\to \PMod_{H\mathbf{K}}^\otimes$ associated to a mixed Weil theory $E$.
We study the automorphism group of $\mathsf{R}'$.
We will show that it is represented by $G_k$.
Here for a finite Galois extension $L$,
we regard $\Gal(L/k)$ as the constant derived affine group scheme
over $H\mathbf{K}$
and we think of $G_k$ as
the limit of derived affine group schemes $\Gal(L/k)$.

Let $\textup{\'Et}/k$ be the category of \'etale schemes over $k$.
There is a natural functor $\textup{\'Et}/k\to \Cor_0$ determined by graphs.
Then we have the composition
\[
\textup{\'Et}/k\to \Cor_0\to \Art(k) \to \PMod_{H\mathbf{K}}
\]
where the second functor is the natural functor induced by $\Cor_0\to \textup{Comp}(\PSh(\Cor_0))$.
(We often omit to take the simplicial nerves of the ordinary categories.)
Note that the second functor is fully faithful.
The essential image is contained in the heart of $\PMod_{H\mathbf{K}}$
with respect to standard $t$-structure, that is, the category of $\mathbf{K}$-vector spaces (the standard $t$-structure is determined by
a pair of full subcategories: the first consists of spectra which are concentrated in non-negative degrees, and the second consists of spectra which are
concentrated in non-positive degrees).
Then this gives rise to $\textup{\'Et}/k\to \mathbf{K}\textup{-Vect}$.

Now suppose that the mixed Weil theory $E$ is 
either $l$-adic \'etale
cohomology theory or Betti cohomology,
see \cite{CD2} ($\mathbf{K}$ depends on the choice of a mixed Weil cohomology
theory).
Then $\textup{\'Et}/k\to \mathbf{K}\textup{-Vect}$
carries $X$ to the $\mathbf{K}$-vector space
generated by the set of $X(\bar{k})$ ($\bar{k}$ is the algebraic closure).
Applying \cite[6.5]{MVW} (after taking the dual vector spaces)
we see that there exists a unique extension 
$\Cor_0\to \mathbf{K}\textup{-Vect}$ of $\textup{\'Et}/k\to \mathbf{K}\textup{-Vect}$. Such a functor $\Cor_0\to \mathbf{K}\textup{-Vect}$
is given by $\Cor_0\simeq \mathbf{K}[G_k]\textup{-Perm} \to \mathbf{K}\textup{-Vect}$ where $\mathbf{K}[G_k]\textup{-Perm}$
denotes the category of permutational representions
and the second functor is the forgetful functor (it is also
symmetric monoidal). Consequently,
we see that the
restriction $\Cor_0^\sim\to \mathbf{K}\textup{-Vect}$ of $\Art(k)\to \PMod_{H\mathbf{K}}$ to
$\Cor_0^\sim$ (contained in $\Art(k)$
as the symmetric monoidal full subcategory)
is equivalent to the forgetful functor $\mathbf{K}[G_k]\textup{-rep}\to \mathbf{K}\textup{-Vect}$ as symmetric monoidal functors.
Here $\mathbf{K}[G_k]\textup{-rep}$ is the category of finite dimensional
discrete representations of $G_k$.
The stable $\infty$-category  $\mathsf{A}_L$
has the standard $t$-structure,
whose heart is $\textup{Gal}(L/k)\textup{-Perm}^\sim$.
Recall that this idempotent completion is equivalent to
the category of finite dimensional representation of $\textup{Gal}(L/k)$.
The composition
$\mathsf{A}_L^\otimes \to \Art(k)^\otimes \to \PMod_{H\mathbf{K}}^\otimes$
induces a ($\mathbf{K}$-linear) symmetric monoidal functor $(\textup{Gal}(L/k)\textup{-Perm}^\sim)^\otimes\to \mathbf{K}\textup{-Vect}^\otimes$ which we
can identify with the forgetful functor.
According to the main theorem in \cite[Section 5]{FI}
together with the classical Tannaka duality (cf. \cite{DM}),
$\mathsf{A}_L^\otimes \to\PMod_{H\mathbf{K}}^\otimes$
is equivalent to the forgetful functor $\PRep_{\Gal(L/k)}^\otimes\to \PMod_{H\mathbf{K}}^\otimes$.

\begin{Proposition}
\label{toabsolute}
The absolute Galois group $G_k$
is the tannakization of $\mathsf{Art}^\otimes(k)\to \PMod_{H\mathbf{K}}^\otimes$.
\end{Proposition}

\Proof
By Lemma~\ref{artincolimit},
we are reduced to showing that the tannakization of
the forgetful functor
$\mathsf{A}_L^\otimes\to \PMod_{H\mathbf{K}}^\otimes$
is the constant finite group scheme $\textup{Gal}(L/k)$ over $H\mathbf{K}$.
Our claim follows from Corollary~\ref{main}.
\QED

$\bullet$ The author is partly supported by Grant-in-aid for Scientific Reseach 23840003,
Japan Society for the Promotion of Science. He also thanks
Research Institute for Mathematical Sciences for the hospitality.

\end{document}